\documentclass[reqno,12pt]{amsart}
\usepackage[T1]{fontenc}
\usepackage{amsfonts}
\usepackage{pifont}
\usepackage{amsmath}
\usepackage{amsthm}
\usepackage{amssymb}
\usepackage{graphicx}
\usepackage{subcaption}
\usepackage{url}
\usepackage{hyperref}
\usepackage{geometry}


\newtheorem{theorem}[subsection]{Theorem}
\newtheorem{definition}{Definition}[section]
\newtheorem{lemma}{Lemma}[section]
\newtheorem{proposition}{Proposition}[section]
\newtheorem{corollary}{Corollary}[section]
\newtheorem{remark}{Remark}[section]

\title[A doubly nonlinear elliptic problem]{A doubly nonlinear elliptic problem with variable exponents, homogeneous Neumann boundary conditions and generalized logistic source}
\author{Bogdan Maxim}
\email{maxim.bogdan.n6h@student.ucv.ro}
\date{\today}

\begin{document}
	\maketitle

\begin{center} \footnotesize{Department of Mathematics, University of Craiova, Al. I. Cuza Street, no. 13, 200585, Craiova, Romania}
\end{center}

\begin{abstract}
	We study a class of doubly nonlinear elliptic problems with variable exponents and homogeneous Neumann boundary conditions arising in the time discretization, by Rothe's method, of an associated doubly nonlinear parabolic problem. A main feature of the present work is that the structural assumptions imposed on the source term are weakened: in particular, no local Lipschitz continuity with respect to the state variable is required. Instead, we work under a weaker monotonicity-type assumption involving the source term and the nonlinear function $b$. A key ingredient in the analysis is the continuity of the associated Nemytskii operators between suitable variable-exponent Lebesgue spaces.
	
	We first establish existence and uniqueness of weak solutions for the auxiliary elliptic problem arising in the Rothe scheme. We then introduce an order-preserving nonlinear solution operator and prove its $L^2$-continuity. Using monotone iteration, we obtain the existence of minimal and maximal weak solutions for the stationary problem, both on the full admissible interval $[0,\delta_0]$ and on invariant subintervals $[\varepsilon,\delta]$. We also provide sufficient conditions for uniqueness of steady states. These results establish the elliptic framework needed for the subsequent study of the associated doubly nonlinear parabolic problem by Rothe's method.
\end{abstract}

\bigskip


\smallskip


\tableofcontents
	
\section{Introduction}\label{sec1}

Given a mathematical model we are interested mainly in those solutions that are in a known bounded interval $[\varepsilon,\delta]$. For example in image processing we are interested in solutions between $0$ (black) and $1$ (white). In population models with a logistic source we are interested in solutions between $0$ and some $K$ (the maximal density of the population supported by the environment), and when we take into account the \textbf{Allee effect} we get interested in solutions between some $k>0$ (population density has to be at least $k$ in order for the population to avoid extinction) and $K$. 

\bigskip

We want to study of the solutions that lie in a given interval $[\varepsilon,\delta]$ for the following very general abstract doubly nonlinear parabolic problem involving variable exponents with a very wide spectrum of applications:

\begin{equation}\tag{$DNP$}\label{eqpg}
	\begin{cases}\dfrac{\partial b(x,u(t,x))}{\partial t}-\operatorname{div}\mathbf{a}(x,\nabla u)=f\big (x,u(t,x)\big ), & (t,x)\in (0,\infty)\times\Omega\\[3mm] \mathbf{a}(x,\nabla u)\cdot\nu=0, & (t,x)\in (0,\infty)\times\partial\Omega\\[3mm] u(0,x)=u_0(x)\in [\varepsilon,\delta], & x\in\Omega\end{cases}
\end{equation}

This article constitutes the first part of the announced study. Specifically, we establish all the foundational results required to address problem \eqref{eqpg} by means of Rothe’s method. Moreover, we present qualitative results regarding the existence and uniqueness of stationary solutions to problem \eqref{eqpg}. To be more precise we will study the following two problems on $[0,\delta_0]$ and then give a method to find information about the solutions that lie in a subinterval $[\varepsilon,\delta]\subset [0,\delta_0]$\footnote{See Theorem \ref{theorem91}.}:

\begin{equation}\tag{$DE_{\lambda}$}\label{eqedglambda}
	\begin{cases}-\operatorname{div}\mathbf{a}\big (x,\nabla V(x)\big )+\lambda b\big (x,V(x)\big )=g(x), & x\in\Omega\\[3mm] \mathbf{a}(x,\nabla V)\cdot\nu=0, & x\in \partial\Omega\\[3mm] 0\leq V(x)\leq \delta_0, & x\in\Omega\end{cases}
\end{equation}

\noindent and:

\begin{equation}\tag{$E$}\label{eqedg}
	\begin{cases}-\operatorname{div}\mathbf{a}\big (x,\nabla U(x)\big )=f\big (x,U(x)\big ), & x\in\Omega\\[3mm] \mathbf{a}(x,\nabla U)\cdot\nu =0, & x\in \partial\Omega\\[3mm] 0\leq U(x)\leq \delta_0, & x\in\Omega\end{cases}
\end{equation}

 An important advancement was made here by \textbf{dropping the assumption that the source $f$ is locally-Lipschitz}\footnote{See Proposition \ref{proposignom} \textbf{(6)}.} while maintaining a weaker Leray-Lions structural condition, introduced by the author in \cite[hypothesis \textbf{(H7)}]{max2}. Only a weak monotonicity assumption is assumed for the source (see \textbf{(H13)}). The key ingredient was the continuity property of the Nemytskii operator between two Lebesgue spaces with variable exponents.\footnote{See Theorem \ref{athnem} from the Appendix.} The most important results obtained in the present paper are: Theorem \ref{thmepsilon}, Theorem \ref{3thmrot}, Proposition \ref{propK} \textbf{(4)}, Theorem \ref{theorem81} and Theorem \ref{theorem91}.

The hypothesis under which we work are the mainly the same from the previous two papers of the author \cite{max2} and \cite{max3}. This paper generalizez the results obtained in \cite{max2} for the doubly nonlinear case.

\section{Hypotheses and notations}

\noindent We consider the following hypotheses and notations:

\begin{enumerate}
	\item[\textbf{(H1)}] $\Omega\subset\mathbb{R}^N,\ N\geq 2$ is an open, bounded and connected Lipschitz domain.

	\item[\textbf{(H2)}] $p:\overline{\Omega}\to (1,\infty)$ is a log-H\"{o}lder continuous variable exponent with the property that: $p^-:=\displaystyle\min_{x\in\overline{\Omega}}\ p(x)\geq\dfrac{2N}{N+2}$ or a continuous exponent with $p^->\dfrac{2N}{N+2}$. In both situations we may write $W^{1,p(x)}(\Omega)\hookrightarrow L^2(\Omega)$.\footnote{See \cite[Theorem 2.3]{Fan2}.}
	
	\bigskip
	
	\item[$\bullet$] Denote $p^{-}=\displaystyle\min_{x\in\overline{\Omega}}\ p(x) >1$ and $p^+=\displaystyle\max_{x\in\overline{\Omega}}\  p(x)<\infty$. Let also $p'(x)=\dfrac{p(x)}{p(x)-1}$ be the conjugate variable exponent of $p(x)$.

	\bigskip
	
	\item[\textbf{(H3)}] $\Psi:\overline{\Omega}\times (0,\infty)\to (0,\infty)$ with $\Psi(\cdot, s)$ measurable for each $s\in (0,\infty)$ and $\Psi(x,\cdot)\in\operatorname{AC}_{\text{loc}}\big ((0,\infty)\big )$ for a.e. $x\in\Omega$.

	\item[\textbf{(H4)}] $\lim\limits_{s\to 0^+} \Psi(x,s)s=0$ for a.e. $x\in\Omega$.
	
	\item[\textbf{(H5)}] For a.e. $x\in\Omega$ we have that $(0,\infty)\ni s\longmapsto \Psi(x,s)s\in (0,\infty)$ is a strictly increasing function.
	
	\bigskip

	\item[$\bullet$] We take $\mathbf{a}:\overline{\Omega}\times\mathbb{R}^N\to\mathbb{R}^N,\ \mathbf{a}(x,\xi)=\Psi(x,|\xi|)\xi$ if $\xi\neq \mathbf{0}$ and $\mathbf{a}(x,\textbf{0})=\textbf{0}$. 
	\bigskip
	
	\item[$\bullet$] Define $\Phi:\overline{\Omega}\times\mathbb{R}\to [0,\infty),\ \Phi(x,s)=\begin{cases} \Psi(x,|s|)|s|, & s\neq 0\\ 0, & s=0\end{cases}$.
	
	\bigskip
	
	\item[\textbf{(H6)}] $L^{p(x)}(\Omega)\ni v\longmapsto \Phi\big(\cdot,v(\cdot)\big )\in L^{p'(x)}(\Omega)$. This is the same as saying that the Nemytskii operator of $\Phi$ is defined as follows $\mathcal{N}_{\Phi}:L^{p(x)}(\Omega)\to L^{p'(x)}(\Omega)$. An other equivalent version would be: $L^{p(x)}(\Omega)^N\ni\mathbf{v}\longmapsto\mathbf{a}\big (\cdot,\mathbf{v}(\cdot)\big )\in L^{p'(x)}(\Omega)^N$.
	
	\bigskip
	
	\item[$\bullet$] $A:\Omega\times\mathbb{R}^N\to [0,\infty),\ A(x,\xi):=\displaystyle\int_{0}^{|\xi|}\Phi(x,s)\ ds$. Note that $A$ is well-defined since $\Phi(x,\cdot)$ is continuous on $[0,\infty)$ for a.e. $x\in\Omega$.

	\item[$\bullet$] $\mathcal{A}:W^{1,p(x)}(\Omega)\to [0,\infty)$ given by $\mathcal{A}(V)=\displaystyle\int_{\Omega} A(x,\nabla V(x))\ dx$.
	
	\bigskip
	
	\item[\textbf{(H7)}] For any sequence $(V_n)_{n\geq 1}\subset W^{1,p(x)}(\Omega)$ with $\displaystyle\int_{\Omega}|\nabla V_n|^{p(x)}\ dx\to+\infty$ we have that: $\lim\limits_{n\to\infty} \displaystyle\int_{\Omega} A(x,\nabla V_n) =+\infty$.
	
	\bigskip
	
	\item[\textbf{(H8)}] $f:\Omega\times \mathbb{R}\to\mathbb{R}$ is a 
	\textbf{generalized logistic source}, meaning that there is a constant $0<\delta_0<\infty$ such that $ f(x,0)\geq 0$ and $f(x,\delta_0)\leq 0$ for a.e. $x\in\Omega$. 
	
%	\item[\textbf{(H9)}] $f$ is \textbf{locally uniformly Lipschitz on $(0,\delta_0]$}, i.e. for each $0<\varepsilon<\delta_0$ there is some $\gamma_{\varepsilon}>0$ such that $|f(x,s_1)-f(x,s_2)|\leq\gamma_{\varepsilon} |s_1-s_2|$ for any $s_1,s_2\in [\varepsilon,\delta_0]$ and for a.e. $x\in\Omega$.
	
	\item[\textbf{(H9)}] The restriction $f:\Omega\times [0,\delta_0]\to\mathbb{R}$ is a Carath\'{e}odory function, i.e. $[0,\delta_0]\ni s\mapsto f(x,s)$ is continuous for a.e. $x\in\Omega$ and $\Omega\ni x\mapsto f(x,s)$ is measurable for any $s\in [0,\delta_0]$. 
	
	\item[\textbf{(H10)}] $b:\Omega\times [0,\delta_0]\to \mathbb{R}$ is a Carath\'{e}odory function, i.e. $[0,\delta_0]\ni s\mapsto b(x,s)$ is continuous for a.e. $x\in\Omega$ and $\Omega\ni x\mapsto b(x,s)$ is measurable for any $s\in [0,\delta_0]$.
	
	\begin{remark}\label{measurabilityfigueiredo} First, from Lemma 5.2.1 given in \cite[page 153]{Alip} we know that $b:\Omega\times [0,\delta_0]\to\mathbb{R}$ and $f:\Omega\times [0,\delta_0]\to\mathbb{R}$ are both \textbf{jointly measurable}, being Carath\'{e}odory functions. Moreover for any \textbf{measurable} function $V:\Omega\to [0,\delta_0]$ (i.e. $V\in\mathcal{U})$ we have from \cite[Theorem 2.1, page 9]{deFigueiredo1989Ekeland} that $b(\cdot,V):\Omega\to\mathbb{R}$ and $f(\cdot,V):\Omega\to\mathbb{R}$ are also \textbf{measurable}.
		
	\end{remark}
	
	\item[\textbf{(H11)}] $b(\cdot,0),b(\cdot,\delta_0)\in L^2(\Omega)$.
	
	\item[\textbf{(H12)}] $[0,\delta_0]\ni s\mapsto b(x,s)$ is a strictly increasing function for a.e. $x\in\Omega$.
	
	\begin{remark}\label{remarkmodulusb} For a.e. $x\in\Omega$ the following inequality holds:
		\begin{equation}
			|b(x,s)|\leq \max\{|b(x,0)|,|b(x,\delta_0)|\},\ \forall\ s\in [0,\delta_0].
		\end{equation}
	\end{remark}
	
	\item[\textbf{(H13)}] There is a constant $\lambda_0>0$ such that the function $[0,\delta_0]\ni s\mapsto f(x,s)+\lambda_0 b(x,s)$ is strictly increasing for a.e. $x\in\Omega$.
	
	\begin{remark}\label{remlambdaincrease}
		We may modify the value of $\lambda_0$ by taking a bigger value when needed. This is because, from \textnormal{\textbf{(H12)}}, the function $[0,\delta_0]\ni s\mapsto f(x,s)+\lambda b(x,s)=f(x,s)+\lambda_0 b(x,s)+(\lambda-\lambda_0)b(x,s)$ is also strictly increasing (being a sum of strictly increasing functions), for a.e. $x\in\Omega$ and for every $\lambda\geq \lambda_0$.
	\end{remark}

	\bigskip
	\bigskip
	
	\noindent In what follows we will introduce some extra hypothesis that shall be useful when dealing with the uniqueness problem for \eqref{eqedg}:
	
	\bigskip
	
	\item[\textbf{(EH)}] There is some constant $\min\{p^{-},2\}\geq\alpha>1$ such that for a.e. $x\in\Omega$ the function $(0,\infty)\ni s\mapsto\dfrac{\Phi(x,s)}{s^{\alpha-1}}$ is \textbf{increasing} and the function $(0,\delta_0^{\alpha}]\ni s\longmapsto \dfrac{f(x,\sqrt[\alpha]{s})}{\sqrt[\alpha]{s^{\alpha-1}}}$ is \textbf{decreasing}.

	\item[\textbf{(EH$_\Phi$)}] For $b(\cdot,0),\ b(\cdot,\delta_0)\in L^{\infty}(\Omega)$\footnote{This will force $f\in L^{\infty}\big ([0,\delta_0]\times\Omega\big )$ as we shall see in Proposition \ref{prop35} \textbf{(1)}.}, there is some $\alpha\in (1,2)$ with $p^-\geq\alpha$ such that for a.e. $x\in\Omega$ the function $(0,\infty)\ni s\mapsto\dfrac{\Phi(x,s)}{s^{\alpha-1}}$ is \textbf{strictly increasing} and the function $(0,\delta_0^{\alpha}]\ni s\longmapsto \dfrac{f(x,\sqrt[\alpha]{s})}{\sqrt[\alpha]{s^{\alpha-1}}}$ is \textbf{decreasing}.
	
	\item[\textbf{(EH$_f$)}] For $b(\cdot,0),\ b(\cdot,\delta_0)\in L^{\infty}(\Omega)$, there is some $\alpha\in (1,2)$ with $p^-\geq\alpha$ such that for a.e. $x\in\Omega$ the function $(0,\delta_0^{\alpha}]\ni s\longmapsto \dfrac{f(x,\sqrt[\alpha]{s})}{\sqrt[\alpha]{s^{\alpha-1}}}$ is \textbf{strictly decreasing} and the function $(0,\infty)\ni s\mapsto\dfrac{\Phi(x,s)}{s^{\alpha-1}}$ is \textbf{increasing}.

		\bigskip
		\bigskip
		
		\noindent In this article we will use the following notations:
		
		\bigskip
	
	\item[$\bullet$] $\mathcal{U}:=\big\{U\in L^{\infty}(\Omega)\ |\ 0\leq U\leq \delta_0\ \text{a.e. on}\ \Omega \big \}$ and for any $0\leq\varepsilon\leq \delta\leq\delta_0$ we set:  $\mathcal{U}_{[\varepsilon,\delta]}=\big\{U\in L^{\infty}(\Omega)\ |\ \varepsilon\leq U\leq \delta\ \text{a.e. on}\ \Omega \big \}$
	
	\item[$\bullet$] $X=\mathcal{U}\cap W^{1,p(x)}(\Omega)$ is the space in which we want the solution of the problem \eqref{eqedg} to be.
	
	\item[$\bullet$] $\mathcal{M}_{\lambda}:=\big \{g\in L^{2}(\Omega)\ |\ \lambda b(x,0)\leq g(x)\leq \lambda b(x,\delta_0)\ \text{for a.e.}\ x\in\Omega\big \}$ and for any $0\leq\varepsilon\leq\delta\leq\delta_0$ let $\mathcal{M}^{[\varepsilon,\delta]}_{\lambda}:=\big \{g\in L^{2}(\Omega)\ |\ \lambda b(x,\varepsilon)\leq g(x)\leq \lambda b(x,\delta)\ \text{for a.e.}\ x\in\Omega\big \}$, whenever $\lambda>0$.
	
	\begin{remark}\label{remark22} Note that if $\lambda\geq \lambda_0$ then for any $U\in\mathcal{U}$ the function $\Omega\ni x\mapsto g(x):=f(x,U(x))+\lambda b(x,U(x))\in\mathcal{M}_{\lambda}$, because for a.e. $x\in\Omega$ we have that:
		
	\begin{equation}
		\lambda b(x,0)\leq f(x,0)+\lambda b(x,0)\leq g(x)\leq f(x,\delta_0)+\lambda b(x,\delta_0)\leq \lambda b(x,\delta_0).
	\end{equation}
	
	\noindent In the same manner if for $0\leq\varepsilon\leq\delta\leq\delta_0$ we have that $f(x,\varepsilon)\geq 0$ and $f(x,\delta)\leq 0$ for a.e. $x\in\Omega$ then for each $U\in\mathcal{U}_{[\varepsilon,\delta]}$ we get that $\Omega\ni x\mapsto g(x):=f(x,U(x))+\lambda b(x,U(x))\in\mathcal{M}^{[\varepsilon,\delta]}_{\lambda}$. Indeed, for a.e. $x\in\Omega$ we have that:
	
	\begin{equation}
		\lambda b(x,\varepsilon)\leq f(x,\varepsilon)+\lambda b(x,\varepsilon)\leq g(x)\leq f(x,\delta)+\lambda b(x,\delta)\leq \lambda b(x,\delta).
	\end{equation}
	\end{remark}
	
	\item[$\bullet$] $g_{\lambda}:\Omega\times [0,\delta_0]\to\mathbb{R},\ g_{\lambda}(x,s):=f(x,s)+\lambda b(x,s)$, for any $\lambda>0$.
	
	\begin{remark}
		For any $\lambda\geq\lambda_0$ we have that $g_{\lambda}$ is a Carath\'{e}odory function, $[0,\delta_0]\ni s\mapsto g_{\lambda}(x,s)$ is strictly increasing for a.e. $x\in\Omega$.
	\end{remark}
	
	\item[$\bullet$] $\overline{f}:\Omega\times\mathbb{R}\to\mathbb{R}$, $\overline{f}(x,s)=\begin{cases} f(x,0)-\dfrac{\lambda_0}{2}s, & s\in (-\infty,0)\\[3mm] f(x,s), & s\in [0,\delta_0] \\[3mm] f(x,\delta_0)-\tilde{\lambda}_0(s-\delta_0), & s\in (\delta_0,\infty) \end{cases}$, where $0<\tilde{\lambda}_0<\lambda_0$ is any fixed constant.
	
	\item[$\bullet$] $F:\Omega\times [0,\delta_0]\to\mathbb{R}$, $F(x,s)=\displaystyle\int_{0}^s f(x,\tau)\ d\tau$.
	
	\item[$\bullet$] $\mathcal{F}:X\to\mathbb{R},\ \mathcal{F}(V)=\displaystyle\int_{\Omega} F\big (x,V(x)\big )\ dx$.
	
	\item[$\bullet$] $\overline{F}:\Omega\times\mathbb{R}\to\mathbb{R}$,
	$\overline{F}(x,s)=\displaystyle\int_{0}^s\overline{f}(x,\tau)\ d\tau$
	
	$=\begin{cases} f(x,0)s-\dfrac{\lambda_0}{4} s^2, & s<0 \\[3mm] F(x,s)=\displaystyle\int_{0}^s f(x,\tau)\ d\tau, &  s\in [0,\delta_0] \\[3mm] F(x,\delta_0)+f(x,\delta_0)(s-\delta_0)-\dfrac{\tilde{\lambda}_0}{2}(s-\delta_0)^2, & s>\delta_0 \end{cases}$.

	\item[$\bullet$] $\overline{\mathcal{F}}:W^{1,p(x)}(\Omega)\to\mathbb{R},\ \overline{\mathcal{F}}(V)=\displaystyle\int_{\Omega} \overline{F}\big (x,V(x)\big )\ dx$.
	
	\item[$\bullet$] $\overline{b}:\Omega\times\mathbb{R}\to\mathbb{R},\  \overline{b}(x,s):=\begin{cases}b(x,0)+s, & s<0\\[3mm] b(x,s), & s\in [0,\delta_0]\\[3mm] b(x,\delta_0)+s-\delta_0, & s>\delta_0 \end{cases}$.
	
	\item[$\bullet$] $B:\Omega\times [0,\delta_0]\to\mathbb{R}$, $B(x,s)=\displaystyle\int_{0}^s b(x,\tau)\ d\tau$.
	
	\item[$\bullet$] $\mathcal{B}:X\to\mathbb{R},\ \mathcal{B}(V)=\displaystyle\int_{\Omega} B\big (x,V(x)\big )\ dx$.
	
	\item[$\bullet$] $\overline{B}:\Omega\times \mathbb{R}\to\mathbb{R}$, $\overline{B}(x,s)=\displaystyle\int_{0}^s \overline{b}(x,\tau)\ d\tau$
	
	$=\begin{cases} b(x,0)s+\dfrac{s^2}{2}, & s<0 \\[3mm] B(x,s)=\displaystyle\int_{0}^s b(x,\tau)\ d\tau, &  s\in [0,\delta_0] \\[3mm] B(x,\delta_0)+b(x,\delta_0)(s-\delta_0)+\dfrac{(s-\delta_0)^2}{2}, & s>\delta_0 \end{cases}$.
	
	\item[$\bullet$] $\overline{\mathcal{B}}:W^{1,p(x)}(\Omega)\to\mathbb{R},\ \overline{\mathcal{B}}(V)=\displaystyle\int_{\Omega} \overline{B}\big (x,V(x)\big )\ dx$.
	
	\item[$\bullet$] $\overline{\mathcal{C}}:W^{1,p(x)}(\Omega)\to\mathbb{R},\ \overline{\mathcal{C}}(V)=\displaystyle\int_{\Omega} \dfrac{|V(x)|^{p(x)}}{p(x)}\ dx$.

	\bigskip
	\end{enumerate}

	\begin{remark}
		Here we use a more general version of the usual Leray-Lions structural condition:
		
		\begin{equation}\label{eqstruct}
			A(x,\xi)\geq \frak{c}_1 |\xi|^{p(x)}-\frak{c}_2,\ \text{for a.e.}\ x\in\Omega,\ \forall\ \xi\in\mathbb{R}^N,
		\end{equation} 
		
		\noindent for some $\frak{c}_1>0$ and $\frak{c}_2\geq 0$.  It is important to emphasize that there are classes of functions that satisfy \textnormal{\textbf{(H7)}} but do not satisfy \eqref{eqstruct}. Take a look at the following two propositions.
	\end{remark}
	
	\noindent The following proposition gives a characterization of \textbf{(H7)}.
	
	\medskip
	
	\begin{proposition} \textnormal{\textbf{(1)}} Consider the \textbf{coercivity profile} $\gamma:[0,\infty)\to [0,\infty)$ defined by
		
		\begin{equation}
			\gamma(s):=\inf\bigl\{\mathcal{A}(V)\ |\ V\in W^{1,p(x)}(\Omega)\ \text{and}\ \rho_{p(x)}\bigl(|\nabla V|\bigr)\geq s\bigr\},
		\end{equation}
		
		\noindent then $\gamma$ is a well-defined, nondecreasing function with $\lim\limits_{s\to\infty} \gamma(s)=+\infty$. Moreover the following inequality holds:
		
		\begin{equation}\label{bogmax1}
			\mathcal{A}(V)\geq\gamma\left(\rho_{p(x)}\bigl(|\nabla V|\bigr)\right),\ \forall\ V\in W^{1,p(x)}(\Omega).
		\end{equation}	
		
		\noindent \textnormal{\textbf{(2)}} Conversely, if there is a function $\gamma:[0,\infty)\to [0,\infty)$ that is nondecreasing with $\lim\limits_{s\to \infty} \gamma(s)=+\infty$ and \eqref{bogmax1} holds, then \textnormal{\textbf{(H7)}} is true.
	\end{proposition}
	
	\begin{proof} \textbf{(1)} $\bullet$ First we show that $\gamma$ is well-defined. In fact, we will show that for any $s\geq 0$ the set $E_s:=\bigl\{V\in W^{1,p(x)}(\Omega)\ |\ \rho_{p(x)}\bigl(|\nabla V|\bigr)\geq s\bigr\}$ is not empty.
		
		\medskip
		
		\noindent Indeed, if we fix any $\varphi_0\in C^{\infty}_{c}(\Omega)$ that is not a constant function, we get that $|\nabla\varphi_0(x)|>0$ on a subset of $\Omega$ with strictly positive measure. Therefore $\rho_{p(x)}(|\nabla\varphi_0|)=\displaystyle\int_{\Omega} |\nabla\varphi_0(x)|^{p(x)}\ dx>0$. Note that for any $\theta\geq 1$ we have for any $x\in\Omega$ that $|\nabla (\theta\varphi_0)(x)|^{p(x)}=\theta^{p(x)}|\nabla\varphi_0(x)|^{p(x)}\geq \theta^{p^-}|\nabla\varphi_0(x)|^{p(x)}$. From here we deduce that $\rho_{p(x)}\bigl(|\nabla (\theta\varphi_0)|\bigr)\geq \theta^{p^-}\rho_{p(x)}\bigl(|\nabla\varphi_0|\bigr)$ for any $\theta\geq 1$. But for any $\theta\geq\max\left\{1,\left [\dfrac{s}{\rho_{p(x)}\bigl(|\nabla\varphi_0|\bigr)}\right ]^{\frac{1}{p^-}}\right\}$ one has that $\theta\varphi_0\in C^{\infty}_c(\Omega)\subset W^{1,p(x)}(\Omega)$ and $\rho_{p(x)}\bigl(|\nabla (\theta\varphi_0)|\bigr)\geq s$, from which we deduce that $E_s\neq\emptyset$. But then $\bigl\{\mathcal{A}(V)\ |\ V\in W^{1,p(x)}(\Omega)\ \text{and}\ \rho_{p(x)}\bigl(|\nabla V|\bigr)\geq s\bigr\}\neq \emptyset$, and henceforth $0\leq \gamma(s)\leq \mathcal{A}(\theta\varphi_0)<\infty$. This shows that $\gamma:[0,\infty)\to [0,\infty)$ is well-defined.
		
		\medskip
		
		\noindent $\bullet$ Let any $0\leq s_1\leq s_2<\infty$. Remark that $E_{s_2}\subset E_{s_1}$. This happens because if $V\in E_{s_2}$ then $\rho_{p(x)}(|\nabla V|)\geq s_2\geq s_1$, which means that $V\in E_{s_1}$. Taking the infimum of the same functional over a smaller set can only increase the
		infimum. Hence $\gamma(s_1)=\displaystyle\inf_{V\in E_{s_1}}\mathcal{A}(V)\leq \displaystyle\inf_{V\in E_{s_2}}\mathcal{A}(V)=\gamma(s_2)$. 
		
		\medskip
		
		\noindent $\bullet$ Since $\gamma$ is a nondecreasing function the following limit exists: $\lim\limits_{s\to\infty} \gamma(s)=\ell\in [0,\infty]$. Suppose that $0\leq \ell<\infty$. We have that $\gamma(s)\leq \ell,\ \forall\ s\in [0,\infty)$. From the definition of the infimum, for every $n\geq 1$ we get that there is some $V_n\in E_n\subset W^{1,p(x)}(\Omega)$ such that: 
		
		\begin{equation}\label{capasuna1}
			0\leq \mathcal{A}(V_n)\leq \dfrac{1}{n}+\inf_{V\in E_n}\mathcal{A}(V)=\dfrac{1}{n}+\gamma(n)\leq\ell+\dfrac{1}{n}\leq\ell+1,\ \forall\ n\geq 1.
		\end{equation}
		
		\noindent But, from $V_n\in E_n$ for every $n\geq 1$, we get that $\rho_{p(x)}(|\nabla V_n|)=\displaystyle\int_{\Omega} |\nabla V_n(x)|^{p(x)}\ dx\geq n,\ \forall\ n\geq 1$. Therefore we have that $\lim\limits_{n\to\infty} \displaystyle\int_{\Omega} |\nabla V_n(x)|^{p(x)}\ dx=+\infty$. So, from \textbf{(H7)}, we get that $\lim\limits_{n\to\infty} \mathcal{A}(V_n)=+\infty$, which contradicts \eqref{capasuna1}. Thereby, we have shown that $\ell=\infty$, i.e. $\lim\limits_{s\to\infty} \gamma(s)=\infty$.
		
		\medskip
		
		\noindent $\bullet$ Let any $V\in W^{1,p(x)}(\Omega)$ and denote $s_V:=\rho_{p(x)}(|\nabla V|)$. Since $\rho_{p(x)}(|\nabla V|)=s_V\geq s_V$ it follows that $V\in E_{s_V}$. Thence, from the definition of the infimum we get that:
		
		\begin{equation}
			\gamma\bigl(\rho_{p(x)}(|\nabla V|) \bigr)=\gamma(s_V)=\inf_{\tilde{V}\in E_{s_V}}\mathcal{A}(\tilde{V})\leq \mathcal{A}(V).
		\end{equation}
		
	\bigskip
	
	\noindent\textbf{(2)} Just remark that for every sequence $(V_n)_{n\geq 1}\subset W^{1,p(x)}(\Omega)$ with $\rho_{p(x)}\bigl(|\nabla V_n|\bigr)\longrightarrow\infty$ we have that:
	
	\begin{equation}
		\mathcal{A}(V_n)\geq \gamma\left(\rho_{p(x)}\bigl(|\nabla V_n|\bigr)\right)\stackrel{n\to\infty}{\longrightarrow}\infty.
	\end{equation}
		
	\end{proof}
		
		\begin{proposition} For:
		\begin{equation}\Phi(x,s)=\max\{h(x,s),a(x,s)|s|^{p(x)-1}-\tilde{a}(x)\},
		\end{equation}
		
		\noindent where:
		
		\begin{enumerate}
			\item[\ding{202}] $a,h:\Omega\times\mathbb{R}\to [0,\infty)$ are Carath\'{e}odory functions (i.e. $a(\cdot,s),h(\cdot,s)$ are measurable for all $s\in [0,\infty)$ and $[0,\infty)\ni s\mapsto a(x,s)$, $[0,\infty)\ni s\mapsto h(x,s)$ are continuous for a.e. $x\in\Omega$) with the extra property that: $a(x,\cdot),h(x,\cdot)\in\operatorname{AC}_{\textnormal{loc}}\big ((0,\infty) \big )$ for a.e. $x\in\Omega$;
			
			\item[\ding{203}] $h(\cdot,v(\cdot))\in L^{p'(x)}(\Omega)$ for any $v\in L^{p(x)}(\Omega)$;
			
			\item[\ding{204}] $\tilde{a}\in L^{p'(x)}(\Omega)^+\setminus L^{\infty}(\Omega)$ ($\tilde{a}$ is \textbf{unbounded});
			
			\item[\ding{205}] $\lim\limits_{s\to\infty} \dfrac{s^{p(x)-1}}{h(x,s)}=\infty$, uniformly with respect to $x\in\Omega$;
			
			\item[\ding{206}] $\underset{\Omega\times (0,\infty)}{\operatorname{ess\ inf}}\ a:=a_0>0$, $a\in L^{\infty}(\Omega\times (0,\infty))$;
			
			\item[\ding{207}] $[0,\infty)\ni s\mapsto h(x,s)$ and $[0,\infty)\ni s\mapsto a(x,s)s^{p(x)-1}$ are strictly increasing for a.e. $x\in\Omega$;
			\item[\ding{208}] $h(x,0)=0$ for a.e. $x\in\Omega$,
		\end{enumerate} 
		
		\noindent we have that \textnormal{\textbf{(H3), (H4), (H5), (H6)}} and \textnormal{\textbf{(H7)}} hold but \eqref{eqstruct} did not hold.
		\end{proposition}
		
		\begin{proof} For any $s\in (0,\infty)$ we have that $\Omega\ni x\mapsto\Phi(x,s)$ is measurable. Therefore $\Omega\ni x\mapsto\Psi(x,s)=\dfrac{\Phi(x,s)}{s}$ is also measurable $\forall\ s\in (0,\infty)$. Now for a fixed $x\in\Omega$ and any $[a,b]\subset (0,\infty)$ the function $[a,b]\ni s\mapsto \Phi(x,s)$ is absolutely continuous, from \ding{202}. This is true since the product of two absolutely continuous functions on a compact interval is also absolutely continuous and the maximum between two absolutely continuous functions is also absolutely continuous. In the same manner $[a,b]\ni s\mapsto\Psi(x,s)=\Phi(x,s)\cdot\dfrac{1}{s}$ is absolutely continuous. This proves \textbf{(H3)}. Note that $\lim\limits_{s\searrow 0}\Psi(x,s)s=\lim\limits_{s\searrow 0} \Phi(x,s)=\max\{h(x,0),0-\tilde{a}(x)\}=\max\{0,-\tilde{a}(x)\}=0$ as \textbf{(H4)} requires. Hypothesis \textbf{(H5)} follows directly from \ding{207}.
			
	\noindent Set any $v\in L^{p(x)}(\Omega)\Rightarrow |v|^{p(x)-1}\in L^{\frac{p(x)}{p(x)-1}}(\Omega)=L^{p'(x)}(\Omega)$. From \ding{203} we have that $h\big (\cdot,v(\cdot)\big )\in L^{p'(x)}(\Omega)$. Also $|a(x,v(x))|v(x)|^{p(x)-1}-\tilde{a}(x)|\leq \Vert a\Vert_{L^{\infty}(\Omega\times (0,\infty))}|v|^{p(x)-1}+\tilde{a}(x)\in L^{p'(x)}(\Omega)$. Therefore $\Phi\big (\cdot,v(\cdot)\big )\in L^{p'(x)}(\Omega)$ and \textbf{(H6)} is true.
	
	\noindent Now it's time to show that \textbf{(H7)} holds. Let any $(V_n)_{n\geq 1}\subset W^{1,p(x)}(\Omega)$ with $\displaystyle\int_{\Omega}|\nabla V_n|^{p(x)}\ dx\to\infty$. Thence:
	
	\begin{align*}
	\displaystyle\int_{\Omega} A\big (x,\nabla V_n(x)\big )\ dx &=\int_{\Omega}\int_{0}^{|\nabla V_n(x)|} \Phi(x,s)\ ds\ dx\geq\int_{\Omega}\int_{0}^{|\nabla V_n(x)|} a(x,s)|s|^{p(x)-1}-\tilde{a}(x)\ ds\ dx \\
	&\geq \int_{\Omega}\int_{0}^{|\nabla V_n(x)|} a_0|s|^{p(x)-1}-\tilde{a}(x)\ ds\ dx=\int_{\Omega} a_0\dfrac{|\nabla V_n(x)|^{p(x)}}{p(x)}\ -\tilde{a}(x)|\nabla V_n(x)|\ dx\\
\text{(H\"{o}lder ineq.)}\ \ \ \ 	&\geq \dfrac{a_0}{p^+}\int_{\Omega} |\nabla V_n(x)|^{p(x)}\ dx-2\Vert\tilde{a}\Vert_{L^{p'(x)}(\Omega)}\Vert |\nabla V_n|\Vert_{L^{p(x)}(\Omega)}\\
(\forall\ n\geq n_0) \ \ \ \ &\geq \dfrac{a_0}{p^+}\rho_{p(x)}(|\nabla V_n|) -2\Vert\tilde{a}\Vert_{L^{p'(x)}(\Omega)}\rho_{p(x)}(|\nabla V_n|)^{\frac{1}{p^{-}}}\\
&=\rho_{p(x)}(|\nabla V_n|)^{\frac{1}{p^{-}}}\left (\dfrac{a_0}{p^+}\rho_{p(x)}(|\nabla V_n|)^{1-\frac{1}{p^{-}}}-2\Vert\tilde{a}\Vert_{L^{p'(x)}(\Omega)} \right )\longrightarrow\infty,
	\end{align*}
	
\noindent as $n\to\infty$. It was essential that $p^{-}>1$. We have also used the inequality $\Vert |\nabla V_n|\Vert_{L^{p(x)}(\Omega)}\leq \rho_{p(x)}(|\nabla V_n|)^{\frac{1}{p^-}}$ which holds when $\Vert |\nabla V_n|\Vert_{L^{p(x)}(\Omega)}>1$. But since $\rho_{p(x)}(|\nabla V_n|)\to\infty$ implies that $\Vert |\nabla V_n|\Vert_{L^{p(x)}(\Omega)}\to\infty$ we deduce that there is some $n_0\geq 1$ such that $\Vert |\nabla V_n|\Vert_{L^{p(x)}(\Omega)}>1$ for each $n\geq n_0$, as needed.

\noindent In the final part of the proof we show that \eqref{eqstruct} cannot hold. Suppose that there are some $\frak{c}_1>0$ and $\frak{c}_2\geq 0$ such that $A(x,\xi)\geq \frak{c}_1 |\xi|^{p(x)}-\frak{c}_2$ for a.e. $x\in\Omega$ and for all $\xi\in\mathbb{R}^N$. So from the monotonicity of $\Phi(x,\cdot)$ we get that for a.e. $x\in\Omega$:

\begin{equation}
	\Phi(x,|\xi|)|\xi|\geq \int_{0}^{|\xi|} \Phi(x,\tau)\ d\tau =A(x,\xi)\geq \frak{c}_1 |\xi|^{p(x)}-\frak{c}_2,\ \forall\ \xi\in\mathbb{R}^N.
\end{equation}

\noindent By setting $|\xi|=s$ we get that:

\begin{equation}\label{eqstructphi}
\Phi(x,s)s\geq \frak{c}_1 s^{p(x)}-\frak{c}_2\ \text{for any}\ s\in [0,\infty)\ \text{and for a.e.}\ x\in\Omega.
\end{equation}

\noindent From \ding{205} we obtain that $\lim\limits_{s\to\infty} \dfrac{\frak{c}_1 s^{p(x)}-\frak{c}_2}{h(x,s)s}=+\infty$ uniformly with respect to $x\in\Omega$. Thus, there exists some $s_0\geq 0$ such that: 

\begin{equation}\label{eqstructcontr}
	\frak{c}_1 s^{p(x)}-\frak{c}_2>h(x,s)\cdot s,\ \forall\ s\geq s_0,\ \text{for a.e.}\ x\in\Omega.
\end{equation}

\noindent For any fixed $s\geq s_0$, set $\Omega_s=\{x\in\Omega\ |\ \Phi(x,s)=h(x,s)\}$. Suppose that $|\Omega_s|=0$. Then $h(x,s)<a(x,s)s^{p(x)}-\tilde{a}(x)$ for a.e. $x\in\Omega$. This shows that: $\tilde{a}(x)<a(x,s)s^{p(x)}-h(x,s)\leq \Vert a\Vert_{L^{\infty}(\Omega\times (0,\infty))}\max\{s^{p^+},s^{p^-}\}<\infty$ for a.e. $x\in\Omega$. Hence $\tilde{a}\in L^{\infty}(\Omega)$, which contradicts \ding{204}. In conclusion $|\Omega_s|>0$ and \eqref{eqstructcontr} combined with \eqref{eqstructphi} gives us:

\begin{equation}
	\frak{c}_1 s^{p(x)}-\frak{c}_2>h(x,s)s=\Phi(x,s)s\geq\frak{c}_1 s^{p(x)}-\frak{c}_2,\ \text{for a.e.}\ x\in\Omega_s
\end{equation}

\noindent which represents the desired contradiction. So \eqref{eqstruct} cannot hold no matter what constants $\frak{c}_1>0$ and $\frak{c}_2\geq 0$ we choose.

		\begin{proposition}[\textbf{The signomial source}]\label{proposignom} For the source function $f:\Omega\times [0,\infty)\to\mathbb{R}$ given by:
		
		\begin{equation}
			f(x,s)=\sum_{i=1}^n a_i(x)s^{q_i(x)}-\sum_{j=1}^m b_j(x)s^{r_j(x)}+c(x),
		\end{equation}
		
		\noindent where:
		
		\begin{enumerate}
			\item[\ding{202}] $a_1,a_2,\hdots, a_n, b_1,b_2,\hdots, b_m, c\in L^{\infty}(\Omega)^+$, where $n,m\geq 1$;
			
			\item[\ding{203}] $q_1,q_2,\hdots,q_n:\Omega\to (0,\alpha-1)$  and $r_1,r_2,\hdots, r_m:\Omega\to (\alpha-1,\infty)$ are measurable exponents, for some $\alpha\in (1,2)$;
			
			\item[\ding{204}] $\underset{x\in\Omega}{\operatorname{ess\ inf}}\ a_1:=a_0>0$ and $\underset{x\in\Omega}{\operatorname{ess\ inf}}\ b_1:=b_0>0$;
			
			\item[\ding{205}] $\displaystyle\min_{i\in\overline{1,n}}\underset{x\in\Omega}{\operatorname{ess\ inf}} \big (r_1-q_i\big )\geq p_0>0$ and $\displaystyle\min_{j\in\overline{1,m}}\underset{x\in\Omega}{\operatorname{ess\ inf}} \big (r_j-q_1\big )\geq p_0>0$,
		\end{enumerate}
		
		\noindent we have that:
		
		\begin{enumerate}
			\item[\textnormal{\textbf{(1)}}] $f$ is a Carath\'{e}odory function.
			
			\item[\textnormal{\textbf{(2)}}] There is some $\delta_0>0$ with $\underset{x\in\Omega}{\operatorname{ess\ sup}} f(x,\delta)<0$ for each $\delta\in [\delta_0,\infty)$.
			
			\item[\textnormal{\textbf{(3)}}] There is some $\varepsilon_0>0$ such that $\underset{x\in\Omega}{\operatorname{ess\ inf}} f(x,\varepsilon)>0$ for each $\varepsilon\in (0,\varepsilon_0]$.
			
			\item[\textnormal{\textbf{(4)}}] For any $\varepsilon\in (0,\varepsilon_0]$ and any $b:\Omega\times [0,\delta_0]\to\mathbb{R}$ satisfying \textnormal{\textbf{(H10)}}, \textnormal{\textbf{(H11)}} and \textnormal{\textbf{(H12)}} with the extra property that there is a constant $\ell_{\varepsilon}>0$ such that:
			
			\begin{equation}
				|b(x,s_2)-b(x,s_1)|\geq \ell_{\varepsilon}|s_2-s_1|,\ \forall\ s_1,s_2\in [\varepsilon,\delta_0],\ \text{for a.e.}\ x\in\Omega,
			\end{equation}
			
			\noindent there is a constant $\lambda_\varepsilon>0$ such that the function $[0,\delta_0]\ni s\mapsto f_{[\varepsilon,\delta_0]}(x,s)+\lambda_{\varepsilon} b(x,s)$ is strictly increasing for a.e. $x\in\Omega$, where $f_{[\varepsilon,\delta_0]}:\Omega\times [0,\delta_0]\to\mathbb{R}$ is defined as:
			
			\begin{equation}
				f_{[\varepsilon,\delta_0]}(x,s)=\begin{cases}f(x,s), & s\in [\varepsilon,\delta_0]\\[3mm] f(x,\varepsilon), & s\in [0,\varepsilon) \end{cases}.
			\end{equation}
			
			\item[\textnormal{\textbf{(5)}}] For a.e. $x\in\Omega$ the function $(0,\delta_0^{\alpha}]\ni s\longmapsto \dfrac{f(x,\sqrt[\alpha]{s})}{\sqrt[\alpha]{s^{\alpha-1}}}$ is \textbf{strictly decreasing}.
			
			\item[\textnormal{\textbf{(6)}}] For a.e. $x\in\Omega$ the function $f(x,\cdot)$ is not locally-Lipschitz at $0$. 
		\end{enumerate}
	\end{proposition}
	
	\begin{proof}\textbf{(1)} It is obvious that for any $s\in [0,\infty)$ the function $\Omega\ni x\mapsto f(x,s)$ is measurable. So it is to notice that for a.e. $x\in\Omega$ the function $[0,\infty)\ni s\mapsto f(x,s)$ is continuous. Thus $f$ is a Carath\'{e}odory function.
		
		\noindent\textbf{(2)}  If we set $M:=\displaystyle\max_{i\in\overline{1,n}} \Vert a_i\Vert_{L^{\infty}(\Omega)}$ then for any $s\geq\max\left\{ \left (\dfrac{nM}{b_0}\right )^{\frac{1}{p_0}},1\right\}$ and for a.e. $x\in\Omega$ we will have that:
		
		\begin{align*}
			f(x,s) &\leq \left (\sum_{i=1}^n \Vert a_i\Vert_{L^{\infty}(\Omega)} s^{q_i(x)}\right )-b_0 s^{r_1(x)}+\Vert c\Vert_{L^{\infty}(\Omega)}\\
			&\leq M\left (\sum_{i=1}^n  s^{q_i(x)}\right )-b_0 s^{r_1(x)}+\Vert c\Vert_{L^{\infty}(\Omega)}\\
			&=\Vert c\Vert_{L^{\infty}(\Omega)}+\sum_{i=1}^n\left [ M s^{q_i(x)}-\dfrac{b_0}{n}s^{r_1(x)}\right ]\\
			&=\Vert c\Vert_{L^{\infty}(\Omega)}+\sum_{i=1}^n s^{q_i(x)}\left ( M-\dfrac{b_0}{n}s^{r_1(x)-q_i(x)} \right )\\
		(s\geq 1)\ \ \ \ \ \ \ \	&\leq \Vert c\Vert_{L^{\infty}(\Omega)}+\sum_{i=1}^n s^{q_i(x)}\underbrace{\left ( M-\dfrac{b_0}{n}s^{p_0} \right )}_{\leq 0}\\
			&\leq \Vert c\Vert_{L^{\infty}(\Omega)}+\sum_{i=1}^n M-\dfrac{b_0}{n}s^{p_0}=\Vert c\Vert_{L^{\infty}(\Omega)}+Mn-b_0s^{p_0}\longrightarrow -\infty,
		\end{align*}
		
		\noindent uniformly with respect to $x\in\Omega$, as $s\to\infty$. So there is some $\delta_0>0$ such that for any $\delta\geq\delta_0$ and for a.e. $x\in\Omega$ we will have that $f(x,\delta)\leq -1$. The claim follows.
		
	\end{proof}
	
	\noindent\textbf{(3)} Take any $\varepsilon_0<\min\left \{\left (\dfrac{a_0}{mM}\right )^{\frac{1}{p_0}}, 1\right \}$ and denote $M:=\displaystyle\max_{j\in\overline{1,m}} \Vert b_j\Vert_{L^{\infty}(\Omega)}$. We have that for a.e. $x\in\Omega$ and any $\varepsilon\in (0,\varepsilon_0]$ the following inequalities hold:
	
	\begin{align*}
		f(x,s)&\geq a_1(x)\varepsilon^{q_1(x)}-\sum_{j=1}^m \Vert b_j\Vert_{L^{\infty}(\Omega)}\varepsilon^{r_j(x)}\geq a_0\varepsilon^{q_1(x)}-\sum_{j=1}^m M\varepsilon^{r_j(x)}\\
		&=\sum_{j=1}^m \dfrac{a_0}{m}\varepsilon^{q_1(x)}- M\varepsilon^{r_j(x)}=\sum_{j=1}^m \varepsilon^{q_1(x)}\left (\dfrac{a_0}{m}- M\varepsilon^{r_j(x)-q_1(x)}\right )\\
(\varepsilon<1)\ \ \ \ \ \ \ \ 		&\geq \sum_{j=1}^m\varepsilon^{q_1(x)}\underbrace{\left (\dfrac{a_0}{m}- M\varepsilon^{p_0}\right )}_{>0}\geq \sum_{j=1}^m\varepsilon^{\alpha-1}\left (\dfrac{a_0}{m}- M\varepsilon^{p_0}\right )\\
&=\varepsilon^{\alpha-1}\left (a_0- Mm\varepsilon^{p_0}\right )>0.
	\end{align*}
	
	\noindent\textbf{(4)} For any $\varepsilon\in (0,\varepsilon_0]$ it is straightforward to see (using \ding{202}) that one can find a constant $\gamma_{\varepsilon}>0$ such that $\left |\dfrac{\partial f}{\partial s}(x,s) \right |\leq \gamma_{\varepsilon}$ for all $s\in [\varepsilon,\delta_0]$ and for a.e. $x\in\Omega$. This shows that the function $s\mapsto f(x,s)+\gamma_{\varepsilon} s$ is increasing on $[\varepsilon,\delta_0]$ for a.e. $x\in \Omega$. In conclusion for a fixed $\lambda_{\varepsilon}>\dfrac{\gamma_{\varepsilon}}{\ell_{\varepsilon}}$ and any $\delta_0\geq s_2>s_1\geq\varepsilon$ we will get that:
	
	\begin{align*}
		&f(x,s_2)+\lambda_{\varepsilon} b(x,s_2)-f(x,s_1)-\lambda_{\varepsilon} b(x,s_1)=f(x,s_2)+\gamma_{\varepsilon} s_2 - f(x,s_1)-\gamma_{\varepsilon} s_1+\\
		 &+\lambda_{\varepsilon}\underbrace{\bigl[b(x,s_2)-b(x,s_1)\bigr]}_{\geq 0}-\gamma_{\varepsilon}(s_2-s_1)\\
		\geq &\underbrace{\big [f(x,s_2)+\gamma_{\varepsilon} s_2 - f(x,s_1)-\gamma_{\varepsilon} s_1\big ]}_{\geq 0}+\underbrace{(\lambda_{\varepsilon}\ell_{\varepsilon}-\gamma_{\varepsilon} )(s_2-s_1)}_{>0}>0,
	\end{align*}
	
	\noindent for a.e. $x\in\Omega$. This shows that $[\varepsilon,\delta_0]\ni s\mapsto f(x,s)+\lambda_{\varepsilon}b(x,s)$ is strictly increasing for a.e. $x\in\Omega$. It is obvious that $[0,\varepsilon]\ni s\mapsto f(x,\varepsilon)+\lambda_{\varepsilon} b(x,s)$ is also striclty increasing for a.e. $x\in\Omega$. The conclusion follows.
	
	\noindent\textbf{(5)} We will show that for $\tau:=\sqrt[\alpha]{s}$ the function $(0,\delta_0]\ni \tau\mapsto\dfrac{f(x,\tau)}{\tau^{\alpha-1}}$ is striclty decreasing for a.e. $x\in\Omega$. Indeed,

	\begin{equation}
	\dfrac{f(x,\tau)}{\tau^{\alpha-1}}=\sum_{i=1}^n a_i(x)\tau^{q_i(x)-(\alpha-1)}-\sum_{j=1}^m b_j(x)\tau^{r_j(x)-(\alpha-1)}+c(x)\tau^{1-\alpha},
	\end{equation}
	
	\noindent is a sum of strictly decreasing functions, since $q_i(x)<\alpha-1$, $r_j(x)>\alpha-1$ for a.e. $x\in\Omega$ and $\alpha>1$. It is essential here that \ding{204} holds.
	
		\noindent\textbf{(6)} We will show that for a.e. $x\in \Omega$: $\lim\limits_{s\to 0^+}\dfrac{\partial f}{\partial s}(x,s)=+\infty$. Indeed,
		
		\begin{align*}
			\dfrac{\partial f}{\partial s}(x,s)&=\sum_{i=1}^n a_i(x)q_i(x)s^{q_i(x)-1}-\sum_{j=1}^mb_j(x)r_j(x)s^{r_j(x)-1}\geq a_1(x)q_1(x)s^{q_1(x)-1}-\sum_{j=1}^mb_j(x)r_j(x)s^{r_j(x)-1}\\
			&\geq a_0q_1(x)s^{q_1(x)-1}-\sum_{j=1}^mb_j(x)r_j(x)s^{r_j(x)-1}\\
			&=\sum_{j=1}^m s^{\overbrace{q_1(x)-1}^{<0}}\left(\dfrac{a_0}{m}q_1(x)-b_j(x)r_j(x) \underbrace{s^{\overbrace{r_j(x)-q_1(x)}^{>0}}}_{\to 0,\ \text{as}\ s\to 0^+}\right )\stackrel{s\to 0^+}{\longrightarrow} +\infty.
		\end{align*}
	
		\end{proof}
		
		\section{Preliminary results}
		
		\begin{proposition}\label{4prop1} The following basic and well-known properties hold:
			\begin{enumerate}
				\item[\textnormal{\textbf{(1)}}] For any $x\in\overline{\Omega}$ and any $\xi_1,\xi_2\in\mathbb{R}^N$ we have that:
				
				\begin{equation}
					\big (\mathbf{a}(x,\xi_1)-\mathbf{a}(x,\xi_2)\big )\cdot (\xi_1-\xi_2)\geq 0,
				\end{equation}
				
				\noindent with equality iff $\xi_1=\xi_2$.
				
				\item[\textnormal{\textbf{(2)}}] For a.e. $x\in\Omega$ and for all $\xi\in\mathbb{R}^N$ we have that $\nabla_{\xi} A(x,\xi)=\mathbf{a}(x,\xi)$.

				%This equality holds for all $(x,\xi)\in\overline{\Omega}\times\mathbb{R}^N$ if $[0,\infty)\ni s\mapsto\Phi(x,s)$ is continuous for each $x\in\overline{\Omega}$.

				\item[\textnormal{\textbf{(3)}}] For a.e. $x\in\Omega$ one has:
				
				\begin{equation}\label{4eqAa}
					\mathbf{a}(x,\xi_2)\cdot (\xi_2-\xi_1)\geq A(x,\xi_2)-A(x,\xi_1)\geq \mathbf{a}(x,\xi_1)\cdot (\xi_2-\xi_1),\ \forall\ \xi_1,\xi_2\in\mathbb{R}^N.
				\end{equation}
				
				\noindent Equality occurs only when $\xi_1=\xi_2$.
				
				\item[\textnormal{\textbf{(4)}}] For a.e. $x\in\Omega$ we have that $A(x,\cdot)$ is an even function that is also \textbf{strictly convex}. 
				
				\item[\textnormal{\textbf{(5)}}] $\mathcal{N}_{\Phi}$ is a continuous and bounded operator and moreover there is some non-negative function $a\in L^{p'(x)}(\Omega)$ and a constant $b\geq 0$ such that: 
				
				\begin{equation}
					\Phi(x,s)\leq a(x)+b|s|^{p(x)-1},\ \forall\ s\in \mathbb{R},\ \text{for a.e.}\ x\in\Omega.
				\end{equation}
				\noindent Thus the following inequality is true:
				
				\begin{equation}
					|\mathbf{a}(x,v)|\leq a(x)+b|v|^{p(x)-1},\ \forall\ v\in\mathbb{R}^N,\ \text{for a.e.}\ x\in\Omega.
				\end{equation}

				\item[\textnormal{\textbf{(6)}}] $\mathcal{A}$ is well-defined, $\mathcal{A}\in C^1\big (W^{1,p(x)}(\Omega)\big )$ and:
				
				\begin{equation}
					\langle \mathcal{A}'(u), \phi\rangle=\int_{\Omega} \mathbf{a}(x,\nabla u(x))\cdot \nabla\phi(x)\ dx,\ \forall\ u,\phi\in W^{1,p(x)}(\Omega).
				\end{equation}
			\end{enumerate}
		\end{proposition}
		
		\begin{proof} \textbf{(1)} Fix some $x\in\overline{\Omega}$ and take any $\xi_1,\xi_2\in\mathbb{R}^N$. At that point, we distinguish three cases:
			
			\begin{enumerate}
				\item[\textbf{(i)}] $\xi_1=\xi_2=\mathbf{0}$. Here there is nothing to prove.
				
				\item[\textbf{(ii)}] $\xi_1\neq \mathbf{0}$ and $\xi_2=\mathbf{0}$. We need to show that $\mathbf{a}(x,\xi_1)\cdot\xi_1>0\ \Leftrightarrow\ \Psi(x,|\xi_1|)\xi_1\cdot\xi_1>0\ \Leftrightarrow\ \Psi(x,|\xi_1|)|\xi_1|^2>0$ which is clearly true from \textbf{(H3)}.
				
				\item[\textbf{(iii)}] $\xi_1,\xi_2\neq\mathbf{0}$. We have that:

				\begin{align*}
					\big (\mathbf{a}(x,\xi_1)-\mathbf{a}(x,\xi_2)\big )\cdot (\xi_1-\xi_2)&=\Psi(x,|\xi_1|)|\xi_1|^2+\Psi(x,|\xi_2|)|\xi_2|^2-\big [\Psi(x,|\xi_1|)+\Psi(x,|\xi_2|) \big ]\xi_1\cdot\xi_2\\
					\text{(Cauchy ineq.)}	&\geq \Psi(x,|\xi_1|)|\xi_1|^2+\Psi(x,|\xi_2|)|\xi_2|^2-\big [\Psi(x,|\xi_1|)+\Psi(x,|\xi_2|) \big ]\cdot |\xi_1|\cdot|\xi_2|\\
					&=\big (|\xi_1|-|\xi_2| \big )\cdot \big (\Psi(x,|\xi_1|)|\xi_1|-\Psi(x,|\xi_2|)|\xi_2| \big )\geq 0\\
					&=\big (|\xi_1|-|\xi_2| \big )\cdot \big (\Phi(x,|\xi_1|)-\Phi(x,|\xi_2|) \big )\geq 0,
				\end{align*}
				
				\noindent because $\Phi(x,\cdot)$ is a strictly increasing function.
				
				\noindent Equality takes place only when $|\xi_1|=|\xi_2|$ and $\xi_1\cdot\xi_2=|\xi_1|\cdot |\xi_2|$ (i.e. there is some $\lambda\geq 0$ with $\xi_1=\lambda\xi_2$ or $\xi_2=\lambda\xi_1$). So we have that $\xi_1=\xi_2$ as claimed.
				
			\end{enumerate}
			
			\noindent\textbf{(2)} Let $x\in\overline{\Omega}$ such that $\Phi(x,\cdot)$ is continuous on $[0,\infty)$. If $\xi\neq 0$ we get from the \textit{chain rule} and the \textit{fundamental theorem of calculus} that:
			
			\begin{equation}
				\nabla_{\xi} A(x,\xi)=\Phi(x,|\xi|)\cdot \nabla_{\xi}|\xi|=\Psi(x,|\xi|)\cdot|\xi|\cdot\dfrac{\xi}{|\xi|}=\Psi(x,|\xi|)\xi=\mathbf{a}(x,\xi). 
			\end{equation}
			
			\noindent For $\xi=0$ we get for each $k\in\overline{1,N}$ that:
			
			\begin{equation}
				\dfrac{\partial A}{\partial\xi_k}(x,\mathbf{0})=\lim\limits_{\xi_k\to 0} \dfrac{A(x,\xi_k\mathbf{e}_k)-A(x,\mathbf{0})}{\xi_k}=\lim\limits_{\xi_k\to 0}\dfrac{1}{\xi_k}\int_{0}^{|\xi_k|} \Phi(x,s)\ ds=0.
			\end{equation}
			
			\noindent We have used the \textit{squeezing principle} and the monotonicity of $\Phi(x,\cdot)$: $\left | \dfrac{1}{\xi_k}\displaystyle\int_{0}^{|\xi_k|} \Phi(x,s)\ ds \right |\leq \dfrac{1}{|\xi_k|}\cdot |\xi_k|\cdot |\Phi(x,|\xi_k|)|=|\Phi(x,|\xi_k|)\stackrel{\xi_k\to 0}{\to} \Phi(x,0)=0$. Thus $\nabla_{\xi} A(x,\mathbf{0})=\mathbf{0}=\mathbf{a}(x,\mathbf{0})$.
			
			\bigskip
			
			\noindent\textbf{(3)} We prove just the right part of the inequality. The left part follows from the right one by simply switching $\xi_1$ and $\xi_2$. Suppose that $|\xi_2|\geq |\xi_1|$. If $\xi_1=\mathbf{0}$ we have that $\mathbf{a}(x,\xi_1)=\mathbf{0}$ and the inequality became $\displaystyle\int_0^{|\xi_2|}\Phi(x,s)\ ds\geq 0$ which is true, and the equality takes place only when $\xi_2=\mathbf{0}$.

			\noindent For $|\xi_2|\geq |\xi_1|>0$ our inequality is equivalent to: $\displaystyle\int_{|\xi_1|}^{|\xi_2|} \Phi(x,s)\ ds\geq \Psi(x,|\xi_1|)\xi_1\cdot\xi_2-\Psi(x,|\xi_1|)|\xi_1|^2$. This is true since from the monotonicity of $\Phi(x,\cdot)$ we derive that $\displaystyle\int_{|\xi_1|}^{|\xi_2|} \Phi(x,s)\ ds\geq \displaystyle\int_{|\xi_1|}^{|\xi_2|} \Phi(x,|\xi_1|)\ ds=\big (|\xi_2|-|\xi_1| \big )\Psi(x,|\xi_1|)|\xi_1|=\Psi(x,|\xi_1|)|\xi_1|\cdot|\xi_2|-\Psi(x,|\xi|)|\xi_1|^2\geq \Psi(x,|\xi_1|)\xi_1\cdot\xi_2-\Psi(x,|\xi_1|)|\xi_1|^2$ from Cauchy inequality. Since $\Phi(x,\cdot)$ is strictly increasing the first inequality becomes equality only when $|\xi_1|=|\xi_2|$. The second inequality becomes equality iff $\xi_1\cdot\xi_2=|\xi_1|\cdot |\xi_2|$, i.e. $\xi_1=\lambda\xi_2$ or $\xi_2=\lambda\xi_1$ for some $\lambda>0$. Since $|\xi_1|=|\xi_2|$ it follows that $\lambda=1$ which means that $\xi_1=\xi_2$.
			
			\bigskip
			
			\noindent For $|\xi_1|\geq |\xi_2|$ and $\xi_2=\mathbf{0}$ the inequality is as follows: $-A(x,\xi_1)\geq -\mathbf{a}(x,\xi_1)\cdot \xi_1\ \Longleftrightarrow\ \Phi(x,|\xi_1|)|\xi_1|=\mathbf{a}(x,\xi_1)\cdot \xi_1\geq \displaystyle\int_{0}^{|\xi_1|}\Phi(x,s)\ ds$ which is true from the strict monotonicity of $\Phi(x,\cdot)$. Indeed: $\displaystyle\int_{0}^{|\xi_1|}\Phi(x,s)\ ds\leq \displaystyle\int_{0}^{|\xi_1|}\Phi(x,|\xi_1|)\ ds=\Phi(x,|\xi_1|)|\xi_1|$. Equality is attained only if $\xi_1=\mathbf{0}=\xi_2$.
			
			\bigskip
			
			\noindent For $|\xi_1|\geq |\xi_2|>0$ we need to show that: $A(x,\xi_1)-A(x,\xi_2)\leq \mathbf{a}(x,\xi_1)\cdot (\xi_1-\xi_2)$. The same arguments apply here:
			
			\begin{align*}
				A(x,\xi_1)-A(x,\xi_2)&=\int_{|\xi_2|}^{|\xi_1|}\Phi(x,s)\ ds\leq \int_{|\xi_2|}^{|\xi_1|}\Phi(x,|\xi_1|)\ ds=\Phi(x,|\xi_1|)|\xi_1|-\Phi(x,|\xi_1|)|\xi_2|\\
				\text{(Cauchy ineq.)}\ \ \ \ \ \	 &=\mathbf{a}(x,\xi_1)\cdot\xi_1-\Psi(x,|\xi_1|)|\xi_1|\cdot|\xi_2|\leq \mathbf{a}(x,\xi_1)\cdot\xi_1-\Psi(x,|\xi_1|)\xi_1\cdot\xi_2\\
				&=\mathbf{a}(x,\xi_1)\cdot\xi_1-\mathbf{a}(x,\xi_1)\cdot\xi_2=\mathbf{a}(x,\xi_1)\cdot (\xi_1-\xi_2).
			\end{align*}
			
			\noindent Equality takes place only when $|\xi_1|=|\xi_2|$ and $\xi_1\cdot\xi_2=|\xi_1|\cdot |\xi_2|$ which means that $\xi_1=\xi_2$.
			
			\bigskip
			
			\noindent\textbf{(4)} $A(x,-\xi)=A(x,\xi)$ for any $\xi\in\mathbb{R}^N$. To show that $A(x,\cdot)$ is strictly convex we take $\xi_1,\xi_2\in\mathbb{R}^N$, $\lambda\in (0,1)$ and then using \textbf{(3)} we obtain that:
			
			\begin{equation}
				\begin{cases} \lambda A(x,\xi_1)-\lambda A(x,\lambda\xi_1+(1-\lambda)\xi_2)\geq \lambda(1-\lambda)\mathbf{a}(x,\lambda\xi_1+(1-\lambda)\xi_2)\cdot(\xi_1-\xi_2)\\[3mm] (1-\lambda) A(x,\xi_2)-(1-\lambda)A(x,\lambda\xi_1+(1-\lambda)\xi_2)\geq -\lambda(1-\lambda)\mathbf{a}(x,\lambda\xi_1+(1-\lambda)\xi_2)\cdot (\xi_1-\xi_2)\end{cases}
			\end{equation}
			
			\noindent Adding these two relations will result in $\lambda A(x,\xi_1)+(1-\lambda)A(x,\xi_2)\geq A(x,\lambda\xi_1+(1-\lambda)\xi_2)$. Equality takes place only when $\xi_1=\xi_2=\lambda\xi_1+(1-\lambda)\xi_2$.
			
			\bigskip
			
			\noindent\textbf{(5)} The proof is a direct application of Theorem \ref{athnem}.
			
			\bigskip
			
			\noindent\textbf{(6)} First we need to explain why is $\mathcal{A}(u)$ a real number for each $u\in W^{1,p(x)}(\Omega)$. Since $A(x,\nabla u(x))\geq 0$ for a.e. $x\in\Omega$ we conclude that the integral $\displaystyle\int_{\Omega}A(x,\nabla u(x))\ dx$ exists, but it can be $\infty$. Now, using the fact that $\Phi(x,\cdot)$ is increasing on $[0,\infty)$ for a.e. $x\in\Omega$ we get for a.e. $x\in\Omega$  that:
			
			\begin{equation}
				0\leq A(x,\nabla u(x))=\int_{0}^{|\nabla u(x)|}\Phi(x,s)\ ds\leq \int_{0}^{|\nabla u(x)|}\Phi(x,|\nabla u(x)|)\ ds=\Phi(x,|\nabla u(x)|)|\nabla u(x)|.
			\end{equation}
			
			\noindent Since $u\in W^{1,p(x)}(\Omega)$ we have that $|\nabla u|\in L^{p(x)}(\Omega)$, and from \textbf{(H6)} we deduce that $\Phi(\cdot,|\nabla u|)\in L^{p'(x)}(\Omega)$. Finally, using \textit{H\"{o}lder inequality} we deduce that $\Phi(\cdot,|\nabla u|)|\nabla u|\in L^{1}(\Omega)$ as needed. So:
			
			\begin{equation}
				0\leq \mathcal{A}(u)\leq \int_{\Omega}\Phi(x,|\nabla u(x)|)|\nabla u(x)|\ dx<\infty.
			\end{equation}
			
			\bigskip
			
			\noindent Secondly, since $\phi\in W^{1,p(x)}(\Omega)$ we get that $|\nabla\phi|\in L^{p(x)}(\Omega)$. Thus:
			
			\begin{align*}\left | \int_{\Omega} \mathbf{a}(x,\nabla u)\cdot\nabla\phi\ dx\right |&\leq \int_{\Omega} |\mathbf{a}(x,|\nabla u|)|\cdot |\nabla\phi|\ dx = \int_{\Omega} |\Phi(x,|\nabla u|)|\cdot |\nabla\phi|\ dx \\
				\text{(H\"{o}lder)}\ \ \ \ \ \ \ 	&\leq 2\Vert \Phi(x,|\nabla u|)\Vert_{L^{p'(x)}(\Omega)}\cdot \Vert |\nabla\phi(x)|\Vert_{L^{p(x)}(\Omega)}<\infty.
			\end{align*}

			\noindent So the mapping $x\longmapsto a(x,\nabla u(x))\cdot\nabla\phi(x)$ is in $L^1(\Omega)$ for each $\phi\in W^{1,p(x)}(\Omega)$. Therefore we can define the linear operator $L_u:W^{1,p(x)}(\Omega)\to\mathbb{R}$ by:
			
			\begin{equation}
				L_u(\phi)=\int_{\Omega} \mathbf{a}(x,\nabla u(x))\cdot\nabla\phi(x)\ dx.
			\end{equation}
			
			\noindent Next, we prove that $L_u\in W^{1,p(x)}(\Omega)^*$, i.e. it is a bounded linear operator. Indeed, for any $\phi\in W^{1,p(x)}(\Omega)$ we have that:
			
			\begin{align*}
				|L_u(\phi)|&=\left |\int_{\Omega}\mathbf{a}(x,\nabla u(x))\cdot\nabla\phi(x)\ dx \right |\leq \int_{\Omega} |\mathbf{a}(x,\nabla u(x))|\cdot |\nabla\phi(x)|\ dx\\
				&=\int_{\Omega} \Phi(x,|\nabla u(x)|)\cdot |\nabla\phi(x)|\ dx\\
				\text{(H\"{o}lder ineq.)}\ \ \ \ \ \ \		&\leq 2\Vert \Phi(\cdot,|\nabla u|)\Vert_{L^{p'(x)}(\Omega)}\cdot \Vert |\nabla\phi|\Vert_{L^{p(x)}(\Omega)}\leq  2\Vert \Phi(\cdot,|\nabla u|)\Vert_{L^{p'(x)}(\Omega)}\cdot \Vert \phi\Vert_{W^{1,p(x)}(\Omega)}.
			\end{align*}
			
			\noindent Observe now that, for any sequence $(\phi_n)_{n\geq 1}\subset W^{1,p(x)}(\Omega)$ with $\lim\limits_{n\to\infty} \Vert \phi\Vert_{W^{1,p(x)}(\Omega)}=0$ we may write:
			
			\begin{align*}
				|\mathcal{A}(u+\phi_n)-\mathcal{A}(u)-L_{u}(\phi_n)|&=\left |\int_{\Omega} A(x,\nabla u+\nabla\phi_n)-A(x,\nabla u)-\mathbf{a}(x,\nabla u)\cdot \nabla\phi_n\ dx \right |\\
				&\leq\int_{\Omega} |A(x,\nabla u+\nabla\phi_n)-A(x,\nabla u)-\mathbf{a}(x,\nabla u)\cdot \nabla\phi_n|\ dx\\
				\text{\eqref{4eqAa}}\ \ \ \ \ \ \ 	&=\int_{\Omega} A(x,\nabla u+\nabla\phi_n)-A(x,\nabla u)-\mathbf{a}(x,\nabla u)\cdot \nabla\phi_n\ dx\\
				\text{\eqref{4eqAa}}\ \ \ \ \ \ \ 	&\leq \int_{\Omega} \big (\mathbf{a}(x,\nabla u+\nabla\phi_n)-\mathbf{a}(x,\nabla u)\big )\cdot \nabla\phi_n\ dx\\
				\text{(Cauchy ineq. and \textbf{(H6)})}\ \ \ \ \ \ \ &\leq  \int_{\Omega} \underbrace{\big |\mathbf{a}(\cdot,\nabla u+\nabla\phi_n)-\mathbf{a}(\cdot,\nabla u)\big |}_{\in L^{p'(x)}(\Omega)}\cdot |\nabla\phi_n|\ dx\\
				\text{(H\"{o}lder ineq.)}\ \ \ \ \ \ \ &\leq 2\big \Vert |\mathbf{a}(\cdot,\nabla u+\nabla\phi_n)-\mathbf{a}(\cdot,\nabla u)|\big \Vert_{L^{p'(x)}(\Omega)}\cdot \big \Vert |\nabla \phi_n|\big \Vert_{L^{p(x)}(\Omega)}\\
				&\leq 2\big \Vert |\mathbf{a}(\cdot,\nabla u+\nabla\phi_n)-\mathbf{a}(\cdot,\nabla u)|\big \Vert_{L^{p'(x)}(\Omega)}\cdot \Vert \phi_n \Vert_{W^{1,p(x)}(\Omega)}.
			\end{align*}
			
			\noindent Consider now $v:=\nabla u\in L^{p(x)}(\Omega)^N$ and define $\Phi_v:\overline{\Omega}\times\mathbb{R}^N\to\mathbb{R},\ \Phi_{v}(x,h)=|\mathbf{a}(x,v+h)-\mathbf{a}(x,v)|,\ \forall\ x\in\overline{\Omega},\ \forall\ h\in\mathbb{R}^N$ which is clearly a Carath\'{e}odory function. The Nemytskii operator associated to $\Phi_v$ is $\mathcal{N}_{\Phi_v}:L^{p(x)}(\Omega)^N\to L^{p'(x)}(\Omega)$. Indeed, since for all $h\in L^{p(x)}(\Omega)^N$ we have that:
			
			\begin{equation*}
				\mathcal{N}_{\Phi_{v}}(h)=|\mathbf{a}(\cdot,v+h)-\mathbf{a}(\cdot,v)|\leq |\mathbf{a}(\cdot,v+h)|+|\mathbf{a}(\cdot,v)|=\mathcal{N}_{\Phi}(|v+h|)+\mathcal{N}_{\Phi}(|v|)\in L^{p'(x)}(\Omega).
			\end{equation*}
			\noindent we deduce from Theorem \ref{athleb} \textbf{(10)} in the Appendix that $\mathcal{N}_{\Phi_{v}}(h)\in L^{p'(x)}(\Omega)$. So, using Remark \ref{arenem}, we finally get that $\mathcal{N}_{\Phi_v}$ is continuous and bounded. Now, for $h_n:=\nabla\phi_n\in L^{p(x)}(\Omega)^N$ we know that $\Vert h_n\Vert_{L^{p(x)}(\Omega)^N}=\Vert |h_n|\Vert_{L^{p(x)}(\Omega)}=\Vert |\nabla\phi_n|\Vert_{L^{p(x)}(\Omega)}\longrightarrow 0$. Thus

			\begin{equation*}\big \Vert |\mathbf{a}(\cdot,\nabla u+\nabla\phi_n)-\mathbf{a}(\cdot,\nabla u)|\big \Vert_{L^{p'(x)}(\Omega)}=\Vert\mathcal{N}_{\Phi_v}(h_n)-\mathcal{N}_{\Phi_v}(0)\Vert_{L^{p'(x)}(\Omega)}\longrightarrow 0.
			\end{equation*}
			
			\noindent Consequently, putting it all together we obtain that:
			
			\begin{equation*}
				0\leq\lim\limits_{n\to\infty} \dfrac{|\mathcal{A}(u+\phi_n)-\mathcal{A}(u)-L_{u}(\phi_n)|}{\Vert\phi_n\Vert_{W^{1,p(x)}(\Omega)}}\leq 2 \lim\limits_{n\to\infty} \big \Vert |\mathbf{a}(\cdot,\nabla u+\nabla\phi_n)-\mathbf{a}(\cdot,\nabla u)|\big \Vert_{L^{p'(x)}(\Omega)}=0.
			\end{equation*}
			
			\noindent This proves that $\mathcal{A}$ is \textit{Fr\'{e}chet differentiable} on $W^{1,p(x)}(\Omega)$ and $\langle\mathcal{A}'(u),\phi\rangle=L_u(\phi),\ \forall\ \phi\in W^{1,p(x)}(\Omega)$.
			
			\begin{remark}\label{4rema} We have proved at \textnormal{\textbf{(6)}} that if $v_n\longrightarrow v$ in $L^{p(x)}(\Omega)^N$ then $\mathbf{a}(\cdot, v_n)\longrightarrow \mathbf{a}(\cdot,v)$ in $L^{p'(x)}(\Omega)^N$.
			\end{remark}
			
			\bigskip
			
			\noindent It remains to show that the (nonlinear) operator $\mathcal{A}':W^{1,p(x)}(\Omega)\to W^{1,p(x)}(\Omega)^*,\ \mathcal{A}'(u):=L_u,\ \forall\ u\in W^{1,p(x)}(\Omega)$ is continuous. First consider any $u\in W^{1,p(x)}(\Omega)$. Then:
			
			\begin{align*}
				\Vert \mathcal{A}'(u)\Vert_{W^{1,p(x)}(\Omega)^*}&=\sup_{\phi\in W^{1,p(x)}(\Omega)\setminus\{0\}}\dfrac{\big |\langle\mathcal{A}'(u),\phi\rangle\big |}{\Vert\phi\Vert_{W^{1,p(x)}(\Omega)}}\leq \sup_{\phi\in W^{1,p(x)}(\Omega)\setminus\{0\}}\dfrac{\displaystyle\int_{\Omega}\big |\mathbf{a}(x,\nabla u)\big |\cdot\big |\nabla\phi\big |\ dx}{\Vert\phi\Vert_{W^{1,p(x)}(\Omega)}}\\
				\text{(H\"{o}lder ineq.)}\ \ \ \ \ \ \ &\leq \sup_{\phi\in W^{1,p(x)}(\Omega)\setminus\{0\}}\dfrac{2\big \Vert|\mathbf{a}(\cdot,\nabla u)|\big\Vert_{L^{p'(x)}(\Omega)} \cdot\big \Vert|\nabla\phi|\big\Vert_{L^{p(x)}(\Omega)}}{\Vert\phi\Vert_{W^{1,p(x)}(\Omega)}}\leq 2\big \Vert|\mathbf{a}(\cdot,\nabla u)|\big\Vert_{L^{p'(x)}(\Omega)}.
			\end{align*}
			
			\noindent Let us now take a sequence $(u_n)_{n\geq 1}\subset W^{1,p(x)}(\Omega)$ with $u_n\longrightarrow u$ in $W^{1,p(x)}(\Omega)$. Observe that:
			
			\begin{align*}
				\big |\langle\mathcal{A}'(u_n)-\mathcal{A}'(u),\phi\rangle \big |&=\left |\int_{\Omega} \big (\mathbf{a}(x,\nabla u_n)-\mathbf{a}(x,\nabla u) \big )\cdot \nabla \phi\ dx \right |\leq \int_{\Omega} \big |\mathbf{a}(x,\nabla u_n)-\mathbf{a}(x,\nabla u) \big |\cdot |\nabla\phi|\ dx\\
				\text{(H\"{o}lder ineq.)}\ \ \ \ \ \ \ &\leq 2\big \Vert |\mathbf{a}(\cdot,\nabla u_n)-\mathbf{a}(\cdot,\nabla u)|\big \Vert_{L^{p'(x)}(\Omega)}\cdot \big\Vert |\nabla\phi |\big \Vert_{L^{p(x)}(\Omega)}\\
				&\leq 2\big \Vert |\mathbf{a}(\cdot,\nabla u_n)-\mathbf{a}(\cdot,\nabla u)|\big \Vert_{L^{p'(x)}(\Omega)}\cdot \Vert \phi \Vert_{W^{1,p(x)}(\Omega)}.
			\end{align*}

			\noindent We have that $\nabla u_n\longrightarrow \nabla u$ in $L^{p(x)}(\Omega)^N$ and from Remark \ref{4rema} we get that:
			
			\begin{align*}
				\big \Vert\mathcal{A}'(u_n)-\mathcal{A}'(u)\Vert_{W^{1,p(x)}(\Omega)^*}&=\sup_{\phi\in W^{1,p(x)}(\Omega)\setminus\{0\}}\dfrac{|\langle\mathcal{A}'(u_n)-\mathcal{A}'(u),\phi\rangle \big |}{\Vert \phi \Vert_{W^{1,p(x)}(\Omega)}}\\
				&\leq 2\big \Vert |\mathbf{a}(\cdot,\nabla u_n)-\mathbf{a}(\cdot,\nabla u)|\big \Vert_{L^{p'(x)}(\Omega)}\stackrel{n\to\infty}{\longrightarrow} 0.
			\end{align*}
			
			\noindent Thus $\mathcal{A}\in C^1\big (W^{1,p(x)}(\Omega) \big)$.
			
		\end{proof}
		
		\begin{proposition}\label{propprel1}
		The function $b:\Omega\times [0,\delta_0]\to\mathbb{R}$ has the following properties:
		\begin{enumerate}
			\item[\textnormal{\textbf{(1)}}] $b(\cdot,V(\cdot))\in L^2(\Omega)$ for every $V\in L^2(\Omega)\cap\mathcal{U}$. In particular if $V_n\to V$ in $L^2(\Omega)$ and $V_n,V\in\mathcal{U}\subset L^2(\Omega)$ for each $n\geq 1$, then $b(\cdot,V_n)\to b(\cdot,V)$ in $L^2(\Omega)$.
			
			\item[\textnormal{\textbf{(2)}}] If $b(\cdot,0),\ b(\cdot,\delta_0)\in L^{\infty}(\Omega)$ then for any $V\in\mathcal{U}$ we get that $b(\cdot, V(\cdot))\in L^{\infty}(\Omega)$.
			
			\item[\textnormal{\textbf{(3)}}] $b\in L^{\infty}\bigl(0,\delta_0;L^2(\Omega)\bigr)$ and in particular $b\in L^2\bigl((0,\delta_0)\times\Omega\bigr)$.
	
		\end{enumerate}
		\end{proposition}
		
		\begin{proof} \textbf{(1)} First, from Remark \ref{measurabilityfigueiredo} we get that $b(\cdot,V(\cdot))$ is a measurable function. Moreover, from Remark \ref{remarkmodulusb} we have for a.e. $x\in\Omega$ that
			
		$$|b(x,V(x))|\leq b(x,V(x))\leq \max\{|b(x,0),|b(x,\delta_0)|\}\in L^2(\Omega).$$ 
		
		\noindent Thus $b(\cdot,V(\cdot))\in L^2(\Omega)$. The last part follows directly from \cite[Theorem 19.1, page 155]{vain}.
			
		\noindent\textbf{(2)} From \textbf{(H12)}: $-\Vert b(\cdot,0)\Vert_{L^{\infty}(\Omega)}\leq b(x,0)\leq b(x,V(x))\leq b(x,\delta_0)\leq \Vert b(\cdot,\delta_0)\Vert_{L^{\infty}(\Omega)}$ for a.e. $x\in\Omega$. So $b(\cdot,V(\cdot))\in L^{\infty}(\Omega)$.
		
		\noindent\textbf{(3)} Using Remark \ref{measurabilityfigueiredo} we first have that $b$ is jointly measurable, i.e. $b:\Omega\times [0,\delta_0]\to\mathbb{R}$ is measurable. Next, from Remark \ref{remarkmodulusb}, we have for a.e. $x\in\Omega$ that:
		
		\begin{equation}
			|b(x,s)|\leq \max\bigl\{|b(x,0)|,|b(x,\delta_0)|\}\leq |b(x,0)|+|b(x,\delta_0)|,\ \forall\ s\in [0,\delta_0].
		\end{equation}
		
		\noindent Thus, from \textbf{(H11)}: $\Vert b(\cdot,s)\Vert_{L^2(\Omega)}\leq \Vert b(\cdot,0)\Vert_{L^2(\Omega)}+\Vert b(\cdot,\delta_0)\Vert_{L^2(\Omega)}$ for any $s\in [0,\delta_0]$. This shows that:
		
		\begin{equation}
			\begin{cases} \Vert b\Vert_{L^{\infty}(0,\delta_0;L^2(\Omega))}&=\underset{s\in [0,\delta_0]}{\operatorname{ess\ sup}}\ \Vert b(\cdot,s)\Vert_{L^2(\Omega)}\leq \Vert b(\cdot,0)\Vert_{L^2(\Omega)}+\Vert b(\cdot,\delta_0)\Vert_{L^2(\Omega)}<\infty\\[3mm] 
			\Vert b\Vert_{L^{2}((0,\delta_0)\times\Omega)}&\leq \Vert\max\bigl\{|b(\cdot,0)|,|b(\cdot,\delta_0)|\} \Vert_{L^2((0,\delta_0)\times\Omega)}\leq \sqrt{\delta_0}\cdot \bigl(\Vert b(\cdot,0)\Vert_{L^2(\Omega)}+\Vert b(\cdot,\delta_0)\Vert_{L^2(\Omega)}\bigr),\end{cases}
		\end{equation}
		
		\noindent i.e. $b\in L^{\infty}\bigl(0,\delta_0;L^2(\Omega)\bigr)$ and $b\in  L^2\bigl((0,\delta_0)\times\Omega\bigr)$.
			
		\end{proof}
		
		\begin{proposition}\label{prop32} The function $\overline{b}:\Omega\times\mathbb{R}\to\mathbb{R}$ has the following properties:
			
			\begin{enumerate}
				\item[\textnormal{\textbf{(1)}}] $\overline{b}$ is a Carath\'{e}odory function.
				
				\item[\textnormal{\textbf{(2)}}] For any $V\in L^2(\Omega)$ we have that $\overline{b}(\cdot,V(\cdot))\in L^2(\Omega)$, i.e. $\mathcal{N}_{\overline{b}}:L^2(\Omega)\to L^2(\Omega)$, where $\mathcal{N}_{\overline{b}}$ is the Nemytskii operator of $\overline{b}$. In particular if $V_n\to V$ in $L^2(\Omega)$ then $\overline{b}(\cdot,V_n)\to \overline{b}(\cdot,V)$ in $L^2(\Omega)$. Also, the following inequality holds for a.e. $x\in\Omega$ and any $s\in\mathbb{R}$:
				
				\begin{equation}
					|\overline{b}(x,s)|\leq \max\{|b(x,0)|,|b(x,\delta_0)|\}+|s|.
				\end{equation}
				
				\item[\textnormal{\textbf{(3)}}] If $b(\cdot,0),\ b(\cdot,\delta_0)\in L^{\infty}(\Omega)$, then for any $V\in L^{\infty}(\Omega)$ we have that $\mathcal{N}_{\overline{b}}(V)=\overline{b}\big (\cdot,V(\cdot)\big )\in L^{\infty}(\Omega)$.
				
				\item[\textnormal{\textbf{(4)}}] For a.e. $x\in\Omega$ the function $\mathbb{R}\ni s\mapsto \overline{b}(x,s)$ is strictly increasing.	
				
			\end{enumerate}
			
		\end{proposition}
		
		\begin{proof} \textbf{(1)} If $s<0$ then $\Omega\ni x\mapsto\overline{b}(x,s)=b(x,0)+s$ is measurable, knowing \textbf{(H10)}. If $s\in [0,\delta_0]$, $\Omega\ni x\mapsto\overline{b}(x,s)=b(x,s)$ is also measurable from \textbf{(H10)}. For $s>\delta_0$ then $\Omega\ni x\mapsto\overline{b}(x,s)=b(x,\delta_0)+s-\delta_0$ is clearly measurable from the same hypothesis \textbf{(H10)}. For a.e. $x\in\Omega$ it is easy to see that $\mathbb{R}\ni s\mapsto \overline{b}(x,s)$ is continuous, using again \textbf{(H10)}. Thus $\overline{b}$ is a Carath\'{e}odory function.
			
		\noindent\textbf{(2)} For any $s\in\mathbb{R}$ we have that:
		
		\begin{align}
		|\overline{b}(x,s)|&=\begin{cases} |b(x,0)+s|, & s<0\\[3mm] |b(x,s)|, & s\in [0,\delta_0] \\[3mm] |b(x,\delta_0)+s-\delta_0|, & s>\delta_0\end{cases}\leq \begin{cases} |b(x,0)|+|s|, & s<0\\[3mm] \max\{|b(x,0)|,|b(x,\delta_0)|\}, & s\in [0,\delta_0] \\[3mm] |b(x,\delta_0)|+|s|-\delta_0, & s>\delta_0\end{cases} \nonumber\\
		&\leq \max\{|b(x,0)|,|b(x,\delta_0)|\}+|s|\label{eqV}.
		\end{align}
			
		\noindent So for any $V\in L^2(\Omega)$ we have a.e. on $\Omega$ that $|\overline{b}(x,V(x))|\leq \max\{|b(x,0)|,|b(x,\delta_0)|\}+|V(x)|\in L^2(\Omega)$. Thence $\overline{b}(\cdot,V(\cdot))\in L^2(\Omega)$. 
		
		\noindent\textbf{(3)} Knowing that $V\in L^{\infty}(\Omega)$, from \textbf{(2)} we get a.e. on $\Omega$ that $|\overline{b}(x,V(x))|\leq \max\{|b(x,0)|,|b(x,\delta_0)|\}+|V(x)|\in L^\infty(\Omega)$. Thence $\overline{b}(\cdot,V(\cdot))\in L^{\infty}(\Omega)$.
		
		\noindent\textbf{(4)} For a.e. fixed $x\in\Omega$ we have that $\overline{b}(x,\cdot)$ is continuous (from \textbf{(1)}) and on each of its branches it is strictly increasing. The conclusion follows with ease.
		
		\end{proof}
		
			\begin{proposition}\label{prop33} The function $\overline{B}:\Omega\times\mathbb{R}\to\mathbb{R}$ has the following properties:
			
			\begin{enumerate}
				\item[\textnormal{\textbf{(1)}}] $\overline{B}$ is a Carath\'{e}odory function. In particular its restriction $B$ is a Carath\'{e}odory function.
				
				\item[\textnormal{\textbf{(2)}}] For a.e. $x\in\Omega$ we have that $\overline{B}(x,\cdot)\in C^1(\mathbb{R})$ and $\dfrac{\partial\overline{B}}{\partial s}(x,s)=\overline{b}(x,s)$ for any $s\in\mathbb{R}$.  In particular, for a.e. $x\in\Omega$ we have that $B(x,\cdot)\in C^1([0,\delta_0])$ and $\dfrac{\partial B}{\partial s}(x,s)=b(x,s)$ for any $s\in [0,\delta_0]$.
				
				\item[\textnormal{\textbf{(3)}}] For any $V\in L^2(\Omega)$ we have that $\overline{B}\big (\cdot,V(\cdot)\big )\in L^1(\Omega)$. In particular for a.e. $x\in\Omega$ and any $s\in\mathbb{R}$:
				
				\begin{equation}
					|\overline{B}(x,s)|\leq \max\{|b(x,0)|,|b(x,\delta_0)|\}|s|+\dfrac{1}{2}s^2+\dfrac{1}{2}\big (s-\delta_0\big )^2.
				\end{equation}
				
				\item[\textnormal{\textbf{(4)}}] For a.e. $x\in\Omega$ and for any $s\in\mathbb{R}$ we have that:
				
				\begin{equation}
					\overline{B}(x,s)\geq -\dfrac{\max\bigl\{b^2(x,0),b^2(x,\delta_0)\bigr\}}{2}-\delta_0|b(x,0)|.
				\end{equation}
				\item[\textnormal{\textbf{(5)}}] For a.e. $x\in\Omega$ the function $\overline{B}(x,\cdot):\mathbb{R}\to\mathbb{R}$ is strictly convex.
				
			\end{enumerate}
			
		\end{proposition}
		
		\begin{proof} \textbf{(1)} and \textbf{(2)} Since from Proposition \ref{propprel1} \textbf{(3)} we know that $b\in L^2\bigl((0,\delta_0)\times\Omega\bigr)$, we can use \textit{Fubini's theorem} given at \cite[Proposition 5.2.2, page 148]{Cohn} in order to deduce that for any fixed $s\in [0,\delta_0]$ the function $\Omega\ni x\mapsto\overline{B}(x,s)=\displaystyle\int_{0}^s\overline{b}(x,\tau)\ d\tau=\displaystyle\int_{0}^s b(x,\tau)\ d\tau$ is measurable. For $s<0$ the function $\Omega\ni x\mapsto b(x,0)+\dfrac{s^2}{2}$ is clearly measurable. Similarly for $s>\delta_0$ the function $\Omega\ni x\mapsto B(x,\delta_0)+b(x,\delta_0)(s-\delta_0)+\dfrac{(s-\delta_0)^2}{2}$ is also measurable. Thus for any $s\in\mathbb{R}$ we proved that $\overline{B}(\cdot,s)$ is measurable. Next, from the \textit{Fundamental theorem of calculus}, because $\mathbb{R}\ni s\mapsto\overline{b}(x,s)$ is continuous for a.e. $x\in\Omega$,  we get that $\overline{B}(x,\cdot)\in C^1(\mathbb{R})\subset C(\mathbb{R})$ and $\dfrac{\partial\overline{B}}{\partial s}(x,s)=\overline{b}(x,s)$ for any $s\in\mathbb{R}$.
			
		\bigskip
		
		\noindent\textbf{(3)} For a.e. $x\in\Omega$ we have that:
		
		\begin{align*}
			\left |\overline{B}\big (x, V(x)\big ) \right |&=\begin{cases}b(x,0)V(x)+\dfrac{V^2(x)}{2}, & V(x)<0\\[3mm] \displaystyle\int_{0}^{V(x)} b(x,\tau)\ d\tau, & V(x)\in [0,\delta_0]\\[3mm] \displaystyle\int_{0}^{\delta_0} b(x,\tau)\ d\tau+b(x,\delta_0)\big (V(x)-\delta_0\big )+\dfrac{1}{2}\big (V(x)-\delta_0\big )^2, & V(x)>\delta_0 \end{cases} \\[3mm]
		\text{\textbf{(H12)}}\ \ \ \ 	&\leq \begin{cases}|b(x,0)|\cdot|V(x)|+\dfrac{V^2(x)}{2}, & V(x)<0\\[3mm] V(x)b\big (x,V(x)\big )\leq V(x)|b(x,\delta_0)|, & V(x)\in [0,\delta_0]\\[3mm] \delta_0 b(x,\delta_0)+|b(x,\delta_0)|\big (V(x)-\delta_0\big )+\dfrac{1}{2}\big (V(x)-\delta_0\big )^2, & V(x)>\delta_0 \end{cases}\\[3mm]
	\text{\textbf{(H11)}}\ \ \ \ 	&\leq \underbrace{\max\{|b(x,0)|,|b(x,\delta_0)|\}}_{\in L^2(\Omega)}\underbrace{|V(x)|}_{\in L^2(\Omega)}+\dfrac{1}{2}V^2(x)+\dfrac{1}{2}\big (V(x)-\delta_0\big )^2\in L^1(\Omega),
		\end{align*}
			
		\noindent from Cauchy inequality. So $\overline{B}(\cdot,V(\cdot))\in L^2(\Omega)$. Choosing $V\equiv s\in\mathbb{R}$ gives us the desired inequality.
		
		\noindent \textbf{(4)} Indeed:
		
		\begin{align*}
			\overline{B}(x,s)&=\begin{cases}b(x,0)s+\dfrac{s^2}{2}, & s<0 \\[3mm] \displaystyle\int_{0}^s b(x,\tau)\ d\tau\geq \int_{0}^s b(x,0)\ d\tau=sb(x,0), &  s\in [0,\delta_0] \\[3mm] \displaystyle\int_{0}^{\delta_0} b(x,\tau)\ d\tau+b(x,\delta_0)(s-\delta_0)+\dfrac{(s-\delta_0)^2}{2}, & s>\delta_0 \end{cases}\\[3mm]
			&\geq \begin{cases}b(x,0)s+\dfrac{s^2}{2}, & s<0 \\[3mm] -\delta_0|b(x,0)|, &  s\in [0,\delta_0] \\[3mm] -\delta_0 |b(x,0)|+b(x,\delta_0)(s-\delta_0)+\dfrac{(s-\delta_0)^2}{2}, & s>\delta_0 \end{cases}\\
			&=\begin{cases}-\dfrac{b(x,0)^2}{2}+\dfrac{(s+b(x,0))^2}{2}, & s<0 \\[3mm] -\delta_0 |b(x,0)|, &  s\in [0,\delta_0] \\[3mm] -\delta_0 |b(x,0)|-\dfrac{b(x,\delta_0)^2}{2}+\dfrac{(s-\delta_0+b(x,\delta_0))^2}{2}, & s>\delta_0 \end{cases}\\[3mm]
			&\geq \begin{cases}-\dfrac{b(x,0)^2}{2}, & s<0 \\[3mm] -\delta_0 |b(x,0)|, &  s\in [0,\delta_0] \\[3mm] -\delta_0 |b(x,0)|-\dfrac{b(x,\delta_0)^2}{2}, & s>\delta_0 \end{cases}\geq -\dfrac{\max\{b^2(x,0),b^2(x,\delta_0)\}}{2}-\delta_0|b(x,0)|.
		\end{align*}
		
		\noindent\textbf{(5)} From \textbf{(2)} the derivative of this function is $\overline{b}(x,\cdot):\mathbb{R}\to\mathbb{R}$, which is strictly increasing, from Proposition \ref{prop32} \textbf{(4)}. The conclusion follows.
		\end{proof}
		
		\begin{proposition}\label{prop34}
			The functional $\overline{\mathcal{B}}:W^{1,p(x)}\to\mathbb{R}$ has the following properties:
			\begin{enumerate}
				
				\item[\textnormal{\textbf{(1)}}] $\overline{\mathcal{B}}$ is well-defined. In particular its restriction $\mathcal{B}$ is also well-defined.
				
				\item[\textnormal{\textbf{(2)}}] $\overline{\mathcal{B}}(V)=\mathcal{B}(V)$ for every $V\in X$.
					
				\item[\textnormal{\textbf{(3)}}] $\displaystyle\inf_{V\in W^{1,p(x)}(\Omega)} \overline{\mathcal{B}}(V)>-\infty$, i.e. $\overline{\mathcal{B}}$ is bounded from below.
					
				\item[\textnormal{\textbf{(4)}}] $\displaystyle\overline{\mathcal{B}}\in C^1\big (W^{1,p(x)}(\Omega)\big )$ and for any $V\in W^{1,p(x)}(\Omega)$:
				
				\begin{equation}
					\langle\overline{\mathcal{B}}'(V),\phi\rangle=\int_{\Omega} \overline{b}\big (x,V(x)\big )\phi(x)\ dx,\ \forall\  \phi\in W^{1,p(x)}(\Omega).
					\end{equation}
			\end{enumerate}
		\end{proposition}
		
		\begin{proof} \textbf{(1)} For any $V\in W^{1,p(x)}(\Omega)$, from \textbf{(H2)} we get that $V\in L^2(\Omega)$. From Proposition \ref{prop33} \textbf{(3)} we have that $\overline{B}(\cdot,V(\cdot))\in L^1(\Omega)$ and therefore $\overline{\mathcal{B}}(V)=\displaystyle\int_{\Omega}\overline{B}(x,V(x))\ dx\in\mathbb{R}$.
			
			\noindent\textbf{(2)} $\overline{\mathcal{B}}(V)=\displaystyle\int_{\Omega}\overline{B}(x,V(x))\ dx=\displaystyle\int_{\Omega} B(x,V(x))\ dx=\mathcal{B}(V)$, because $V(x)\in [0,\delta_0]$ for a.e. $x\in\Omega$.

			\noindent\textbf{(3)} For any $V\in W^{1,p(x)}(\Omega)$, from Proposition \ref{prop33} \textbf{(4)} we have that:
			
			\begin{align*}
				\overline{\mathcal{B}}(V)&=\displaystyle\int_{\Omega}\overline{B}(x,V(x))\ dx\geq -\dfrac{1}{2}\int_{\Omega}\max\bigl\{b^2(x,0),b^2(x,\delta_0)\bigr\}\ dx-\delta_0\int_{\Omega}|b(x,0)|\ dx\\
				&\geq -\dfrac{1}{2}\int_{\Omega}b^2(x,0)\ dx-\dfrac{1}{2}\int_{\Omega}b^2(x,\delta_0)\ dx -\delta_0\int_{\Omega}|b(x,0)|\ dx\\
		\text{(Cauchy ineq.)}\ \ \		&\geq -\dfrac{1}{2}\Vert b(\cdot,0)\Vert_{L^2(\Omega)}^2-\dfrac{1}{2}\Vert b(\cdot,\delta_0)\Vert_{L^2(\Omega)}^2-\delta_0\sqrt{|\Omega|}\Vert b(\cdot,0)\Vert_{L^2(\Omega)}>-\infty.
			\end{align*}

			\noindent\textbf{(4)} 	\noindent Next, we'll show that $\overline{\mathcal{B}}$ is G\^ateaux-differentiable on $W^{1,p(x)}(\Omega)$. Set a direction $\phi\in W^{1,p(x)}(\Omega)$ (so $\phi\in L^{2}(\Omega)$) and consider any sequence $(\varepsilon_n)_{n\geq 1}\subset\mathbb{R}^*$ that is convergent to $0$. Without loss of generality we can assume that $|\varepsilon_n|\leq 1,\ \forall\ n\geq 1$. 
			
			\noindent Consider for each $n\geq 1$ the function: $g_n:\Omega\to\mathbb{R},\ g_n(x)=\dfrac{\overline{B}(x,V(x)+\varepsilon_n\phi(x))-\overline{B}(x,V(x))}{\varepsilon_n}$. From Proposition \ref{prop33} \textbf{(2)} and the classical 1-dimensional \textit{mean value theorem} we have that $g_n(x)=\dfrac{\partial \overline{B}}{\partial s}(x,V(x)+s_n(x)\phi(x))\phi(x)=\overline{b}(x,V(x)+s_n(x)\phi(x))\phi(x)$ for some real number $s_n(x)$ with $|s_n(x)|\leq\varepsilon_n$. Using the continuity of $\overline{b}(x,\cdot)$ for a.e. $x\in\Omega$ (see Proposition \ref{prop32} \textbf{(1)}) we get that:
			
			\begin{equation}
				\lim\limits_{n\to\infty} g_n(x)=\lim\limits_{n\to\infty}\overline{b}(x,V(x)+s_n(x)\phi(x))\phi(x)\stackrel{|s_n(x)|\leq \varepsilon_n}{=}\overline{b}(x,V(x))\phi(x),\ \text{for a.e.}\ x\in\Omega.
			\end{equation}
			
			\noindent This proves that $g_n$ has a pointwise limit. The next step will be to show that $|g_n|$ is bounded by a function from $L^1(\Omega)$ for any $n\geq 1$. Indeed, from Proposition \ref{prop32} \textbf{(2)} we get that:
			
			\begin{align*}
				|g_n(x)|&=|\overline{b}(x,V(x)+s_n(x)\phi(x))|\cdot |\phi(x)|\leq \\
				&\leq \big [\max\bigl\{|b(x,0)|,|b(x,\delta_0)|\bigl\}+|V(x)+s_n(x)\phi(x)| \big ]\cdot |\phi(x)|\\
				&\leq \big [\max\bigl\{|b(x,0)|,|b(x,\delta_0)|\bigl\}+|V(x)|+|s_n(x)|\cdot|\phi(x)| \big ]\cdot |\phi(x)|\\
			 |s_n|\leq\varepsilon_n\leq 1, \textbf{(H11)}\ \ \ \	&\leq \underbrace{\big [\max\bigl\{|b(x,0)|,|b(x,\delta_0)|\bigl\}+|V(x)|+|\phi(x)| \big ]}_{\in L^2(\Omega)}\cdot \underbrace{|\phi(x)|}_{\in L^2(\Omega)}\in L^1(\Omega),\\
			\end{align*}
			
			\noindent because $V,\phi\in W^{1,p(x)}(\Omega)\hookrightarrow L^2(\Omega)$.
			
			\noindent From the \textit{Lebesgue dominated convergence theorem} we get that:
			
			\begin{align*}
				\lim\limits_{n\to \infty}\dfrac{\overline{\mathcal{B}}(V+\varepsilon_n\phi)-\overline{\mathcal{B}}(V)}{\varepsilon_n}&=\lim\limits_{n\to\infty} \int_{\Omega}\dfrac{\overline{B}(x,V(x)+\varepsilon_n\phi(x))-\overline{B}(x,V(x))}{\varepsilon_n}\ dx\\
				&=\lim\limits_{n\to\infty}\int_{\Omega} g_n(x)\ dx=\int_{\Omega} \overline{b}(x,V(x))\phi(x)\ dx.
			\end{align*}
			
			\noindent It is easy to remark that $\partial\overline{\mathcal{B}}(V):W^{1,p(x)}(\Omega)\to\mathbb{R},\ \partial_{\phi}\overline{\mathcal{B}}(V)=\displaystyle\int_{\Omega} \overline{b}(x,V(x))\phi(x)\ dx$ is a linear operator. We show now that it is in fact a bounded linear operator. Indeed:
			
			\begin{align*}
				\left |\partial_{\phi}\overline{\mathcal{B}}(V) \right |&=\left | \int_{\Omega} \overline{b}(x,V(x))\phi(x)\ dx \right |\leq \int_{\Omega} |\overline{b}(x,V(x))|\cdot |\phi(x)|\ dx\\
			\text{Proposition \ref{prop32}\ \bf{(2)}}\ \ \ \ \	&\leq \int_{\Omega} \left ( \max\bigl\{|b(x,0)|,|b(x,\delta_0)|\bigl\}+|V(x)| \right )|\phi(x)|\ dx\\
				\text{(Cauchy ineq.)}	\ \ \ \ \ &\leq \big\Vert \max\bigl\{|b(x,0)|,|b(x,\delta_0)|\bigl\}+|V(x)|\big \Vert_{L^2(\Omega)}\Vert \phi\Vert_{L^2(\Omega)}\\
				\big (W^{1,p(x)}(\Omega)\hookrightarrow L^{2}(\Omega)\big )\ \ \ \ \ \ \ &\lesssim \big\Vert \max\bigl\{|b(x,0)|,|b(x,\delta_0)|\bigl\}+|V(x)|\big \Vert_{L^2(\Omega)}\Vert \phi\Vert_{W^{1,p(x)}(\Omega)}.
			\end{align*}
			
			\noindent Henceforth: $\displaystyle\sup_{\phi\in W^{1,p(x)}(\Omega)\setminus\{0\}}\dfrac{\big |\partial_{\phi}\overline{\mathcal{B}}(V) \big |}{\Vert \phi\Vert_{W^{1,p(x)}(\Omega)}}<\infty$, which shows that $\partial\overline{\mathcal{B}}(V)\in\mathcal{L}(W^{1,p(x)}(\Omega);\mathbb{R})=W^{1,p(x)}(\Omega)^*$, for any $V\in W^{1,p(x)}(\Omega)$. This observation allows us to write that: $\partial\overline{\mathcal{B}}:W^{1,p(x)}(\Omega)\to W^{1,p(x)}(\Omega)^*$. The following step is to prove that $\partial\overline{\mathcal{B}}$ is a continuous mapping.

			\noindent Consider a sequence $(V_n)_{n\geq 1}\subset W^{1,p(x)}(\Omega)$ and some $V\in W^{1,p(x)}(\Omega)$ such that $V_n\to V$ in $W^{1,p(x)}(\Omega)$. Since $W^{1,p(x)}(\Omega)\hookrightarrow L^{2}(\Omega)$ we also have that $V_n\to V$ in $L^{2}(\Omega)$.

			\begin{align*}\big \Vert\partial\overline{\mathcal{B}}(V_n)-\partial\overline{\mathcal{B}}(V)\big\Vert_{W^{1,p(x)}(\Omega)^*}&=\sup_{\phi\in W^{1,p(x)}(\Omega)\setminus\{0\}}\dfrac{\big |\partial_{\phi}\overline{\mathcal{B}}(V_n)-\partial_{\phi}\overline{\mathcal{B}}(V)\big |}{\Vert\phi\Vert_{W^{1,p(x)}(\Omega)}}\\
				&=\sup_{\phi\in W^{1,p(x)}(\Omega)\setminus\{0\}}\dfrac{\left |\displaystyle\int_{\Omega} \big (\overline{b}(x,V_n(x))-\overline{b}(x,V(x)))\cdot\phi(x) dx\right |}{\Vert\phi\Vert_{W^{1,p(x)}(\Omega)}} \\
				&\leq \sup_{\phi\in W^{1,p(x)}(\Omega)\setminus\{0\}}\dfrac{\displaystyle\int_{\Omega} \big |\overline{b}(x,V_n(x))-\overline{b}(x,V(x))|\cdot|\phi(x)| dx}{\Vert\phi\Vert_{W^{1,p(x)}(\Omega)}}\\
				\text{(Cauchy ineq.)}\ \ \ \ \ \ \ \	&\leq \sup_{\phi\in W^{1,p(x)}(\Omega)\setminus\{0\}}\dfrac{\Vert \overline{b}(\cdot,V_n)-\overline{b}(\cdot,V) \Vert_{L^{2}(\Omega)}\Vert \phi\Vert_{L^{2}(\Omega)}}{\Vert\phi\Vert_{W^{1,p(x)}(\Omega)}}\\
				\big (W^{1,p(x)}(\Omega)\hookrightarrow L^{2}(\Omega)\big )\ \ \ \ \ \ \ &\lesssim \sup_{\phi\in W^{1,p(x)}(\Omega)\setminus\{0\}}\dfrac{\Vert \overline{b}(\cdot,V_n)-\overline{b}(\cdot,V)\Vert_{L^{2}(\Omega)}\Vert \phi\Vert_{W^{1,p(x)}(\Omega)}}{\Vert\phi\Vert_{W^{1,p(x)}(\Omega)}}\\
				&=\Vert \overline{b}(\cdot,V_n)-\overline{b}(\cdot,V)\Vert_{L^{2}(\Omega)}\stackrel{n\to\infty}{\longrightarrow} 0,
			\end{align*}

			\noindent from Proposition \ref{prop32} \textbf{(2)}.
			\noindent Finally, from Theorem \ref{athmgatfre}, we obtain that $\overline{\mathcal{B}}\in C^1\big (W^{1,p(x)}(\Omega)\big )$ and:
			
			\begin{equation*}
				\langle \overline{\mathcal{B}}'(V),\phi \rangle=\displaystyle\int_{\Omega} \overline{b}(x,V(x))\phi(x)\ dx,\ \forall\ V,\phi\in W^{1,p(x)}(\Omega).
			\end{equation*}

		\end{proof}
		
		\begin{proposition}\label{prop35} The source function $f:\Omega\times [0,\delta_0]\to\mathbb{R}$ has the following properties:
			
			\begin{enumerate}
				\item[\textnormal{\textbf{(1)}}] For any $s\in [0,\delta_0]$ and for a.e. $x\in\Omega$ the following inequality holds:
				
				\begin{equation}
					\label{oecuatie}
					\lambda_0\big (b(x,0)-b(x,s)\big )\leq f(x,s)\leq\lambda_0\big (b(x,\delta_0))-b(x,s)\big ).
				\end{equation}
				
				\noindent In particular for a.e. $x\in\Omega$ and for any $s\in [0,\delta_0]$ we have: 
				
				\begin{equation}\label{fcevaimportant}
				|f(x,s)|\leq \lambda_0\big (b(x,\delta_0))-b(x,0)\big ).
				\end{equation} 
				
				\item[\textnormal{\textbf{(2)}}] $f(\cdot,0),f(\cdot,\delta_0)\in L^2(\Omega)$.
				
				\item[\textnormal{\textbf{(3)}}] For any $U\in\mathcal{U}$ we have that $f(\cdot,U(\cdot))\in L^2(\Omega)$. In particular if $U_n\to U$ in $L^2(\Omega)$ and $U_n,U\in\mathcal{U}\subset L^2(\Omega)$ for each $n\geq 1$, then $f(\cdot,U_n)\to f(\cdot,U)$ in $L^2(\Omega)$.
				
				\item[\textnormal{\textbf{(4)}}] If $b(\cdot,0),\ b(\cdot,\delta_0)\in L^{\infty}(\Omega)$ then $f(\cdot,0),f(\cdot,\delta_0)\in L^{\infty}(\Omega)$. Furthermore, for any $U\in\mathcal{U}$ one also gets that $f(\cdot,U(\cdot))\in L^{\infty}(\Omega)$.
				
				\item[\textnormal{\textbf{(5)}}] $f\in L^{\infty}\bigl(0,\delta_0;L^2(\Omega)\bigr)$ and in particular $f\in L^2\bigl((0,\delta_0)\times\Omega\bigr)$.
			\end{enumerate}
			
		\end{proposition}
		
		\begin{proof} \textbf{(1)} 	\noindent From \textnormal{\textbf{(H13)}} we have for a.e. $x\in\Omega$: $\lambda_0 b(x,0)\leq f(x,0)+\lambda_0 b(x,0)\leq f(x,s)+\lambda_0 b(x,s)\leq f(x,\delta_0)+\lambda_0 b(x,\delta_0)\leq \lambda_0 b(x,\delta_0)$. Combining this with \textbf{(H12)} yields:
			
			\begin{equation}
				 -\lambda_0\big (b(x,\delta_0))-b(x,0)\big )\leq\lambda_0\big (b(x,0)-b(x,s)\big )\leq f(x,s)\leq\lambda_0\big (b(x,\delta_0))-b(x,s)\big )\leq  \lambda_0\big (b(x,\delta_0))-b(x,0)\big ),
			\end{equation}

\noindent as needed.

		\noindent\textbf{(2)} From \textbf{(1)} we get that $|f(x,0)|,|f(x,\delta_0)|\leq \lambda_0\big (b(x,\delta_0))-b(x,0)\big )\in L^2(\Omega)$. So $f(\cdot,0),\ f(\cdot,\delta_0)\in L^2(\Omega)$.
			
		\noindent\textbf{(3)} From \textbf{(1)} we have for a.e. $x\in\Omega$:
		
		\begin{equation}\label{oecuatie2}
			|f(x,U(x))|\leq\lambda_0\big (b(x,\delta_0))-b(x,0)\big ).
		\end{equation}
		
		\noindent The right-hand side is in $L^2(\Omega)$ (see \textbf{(H11)}) and from Remark \ref{measurabilityfigueiredo} we already know that $f(\cdot,U(\cdot))$ is measurable. Thus we first get that $f(\cdot,U(\cdot))\in L^2(\Omega)$. For the second part of \textbf{(3)} we refer the reader to \cite[Theorem 19.1, page 155]{vain}.
		
		\noindent\textbf{(4)} From \textbf{(1)} we get that for each $s\in [0,\delta_0]$: $|f(x,s)|\leq \lambda_0\big (b(x,\delta_0))-b(x,0)\big )\in L^{\infty}(\Omega)$. So  $f(\cdot,s)\in L^{\infty}(\Omega)$ and in particular $f(\cdot,0),\ f(\cdot,\delta_0)\in L^{\infty}(\Omega)$. From \eqref{oecuatie2} we also get that $f(\cdot,U(\cdot))\in L^{\infty}(\Omega)$.
		
		\noindent\textbf{(5)} Using Remark \ref{measurabilityfigueiredo} we first have that $f$ is jointly measurable, i.e. $f:\Omega\times [0,\delta_0]\to\mathbb{R}$ is measurable. Next, from \eqref{fcevaimportant}, we have for a.e. $x\in\Omega$ that:
		
		\begin{equation}
			|f(x,s)|\leq \lambda_0\bigl(b(x,\delta_0)-b(x,0)\bigr),\ \forall\ s\in [0,\delta_0].
		\end{equation}
		
		\noindent Thus, from \textbf{(H11)}: $\Vert f(\cdot,s)\Vert_{L^2(\Omega)}\leq \lambda_0\Vert b(\cdot,\delta_0)-b(\cdot,0)\Vert_{L^2(\Omega)}$ for any $s\in [0,\delta_0]$. This shows that:
		
		\begin{equation}
			\begin{cases} \Vert f\Vert_{L^{\infty}(0,\delta_0;L^2(\Omega))}&=\underset{s\in [0,\delta_0]}{\operatorname{ess\ sup}}\ \Vert f(\cdot,s)\Vert_{L^2(\Omega)}\leq \lambda_0\Vert b(\cdot,\delta_0)-b(\cdot,0)\Vert_{L^2(\Omega)}<\infty\\[3mm] 
				\Vert f\Vert_{L^{2}((0,\delta_0)\times\Omega)}&\leq \lambda_0\Vert b(\cdot,\delta_0)-b(\cdot,0) \Vert_{L^2((0,\delta_0)\times\Omega)}= \lambda_0\sqrt{\delta_0}\Vert b(\cdot,\delta_0)-b(\cdot,0)\Vert_{L^2(\Omega)},\end{cases}
		\end{equation}
		
		\noindent i.e. $f\in L^{\infty}\bigl(0,\delta_0;L^2(\Omega)\bigr)$ and $f\in  L^2\bigl((0,\delta_0)\times\Omega\bigr)$.
			
		\end{proof}

			\begin{proposition}\label{prop36} The function $\overline{f}:\Omega\times\mathbb{R}\to\mathbb{R}$ has the following properties:
			
\begin{enumerate}
	\item[\textnormal{\textbf{(1)}}] $\overline{f}$ is a Carath\'{e}odory function.
	
	\item[\textnormal{\textbf{(2)}}] For any $\lambda\geq \lambda_0$ the function $\mathbb{R}\ni s\longmapsto\overline{f}(x,s)+\lambda\cdot\overline{b}(x,s)$ is strictly increasing for a.e. $x\in\Omega$.
	
	\item[\textnormal{\textbf{(3)}}] For any $V\in L^2(\Omega)$ we have that $\overline{f}(\cdot,V(\cdot))\in L^2(\Omega)$, i.e. the Nemytskii operator associated to $\overline{f}$ maps $L^2(\Omega)$ into itself. In particular if $V_n\to V$ in $L^2(\Omega)$ then $\overline{f}(\cdot, V_n)\to \overline{f}(\cdot,V)$ in $L^2(\Omega)$. Also, the following inequality holds for a.e. $x\in \Omega$ and any $s\in\mathbb{R}$:
	
	\begin{equation}
		|\overline{f}(x,s)|\leq \lambda_0 |s|+f(x,0)-f(x,\delta_0)+\lambda_0\big (b(x,\delta_0)-b(x,0) \big ).
	\end{equation}
	
	\item[\textnormal{\textbf{(4)}}] If $b(\cdot,\delta_0)\in L^{\infty}(\Omega)$ then for any $V\in L^{\infty}(\Omega)$ we get that $\overline{f}(\cdot,V(\cdot)\big )\in L^{\infty}(\Omega)$.
	
	\item[\textnormal{\textbf{(5)}}]  If \textnormal{\textbf{(EH)}} holds and $b(\cdot,0),\ b(\cdot,\delta_0)\in L^{\infty}(\Omega)$ then for any $\tilde{\lambda}_0\geq \underset{x\in\Omega}{\operatorname{ess\ sup}} -\dfrac{\alpha-1}{\delta_0}f(x,\delta_0)$ the function $[0,\infty)\ni s\mapsto -\overline{F}(x,\sqrt[\alpha]{s})$ is convex and the function $(0,\infty)\ni s\mapsto\dfrac{\overline{f}(x,\sqrt[\alpha]{s})}{\sqrt[\alpha]{s^{\alpha-1}}}$ is decreasing for a.e. $x\in\Omega$.
	
	\item[\textnormal{\textbf{(6)}}]  If \textnormal{\textbf{(EH$_f$)}} holds and $b(\cdot,0),\ b(\cdot,\delta_0)\in L^{\infty}(\Omega)$ then for any $\tilde{\lambda}_0>0$ with $\tilde{\lambda}_0\geq \underset{x\in\Omega}{\operatorname{ess\ sup}} -\dfrac{\alpha-1}{\delta_0}f(x,\delta_0)$ then the function $[0,\infty)\ni s\mapsto -\overline{F}(x,\sqrt[\alpha]{s})$ is \textbf{strictly convex} and the function $(0,\infty)\ni s\mapsto\dfrac{\overline{f}(x,\sqrt[\alpha]{s})}{\sqrt[\alpha]{s^{\alpha-1}}}$ is \textbf{strictly decreasing} for a.e. $x\in\Omega$.\footnote{We can choose $\tilde{\lambda}_0\geq \underset{x\in\Omega}{\operatorname{ess\ sup}} -\dfrac{\alpha-1}{\delta_0}f(x,\delta_0)$ if we increase the value of $\lambda_0$ so that $\lambda_0>\underset{x\in\Omega}{\operatorname{ess\ sup}} -\dfrac{\alpha-1}{\delta_0}f(x,\delta_0)$. We can do this if needed, from Remark \ref{remlambdaincrease}.}
\end{enumerate}

\end{proposition}

		\begin{proof} \textbf{(1)} If $s<0$ then $\Omega\ni x\mapsto\overline{f}(x,s)=f(x,0)-\dfrac{\lambda}{2}s$ which is measurable, knowing \textbf{(H9)}. If $s\in [0,\delta_0]$, $\Omega\ni x\mapsto\overline{f}(x,s)=f(x,s)$ is also measurable from \textbf{(H9)}. For $s>\delta_0$ then $\Omega\ni x\mapsto\overline{f}(x,s)=f(x,\delta_0)-\tilde{\lambda}_0(s-\delta_0)$ is clearly measurable from the same hypothesis \textbf{(H9)}. For a.e. $x\in\Omega$ it is easy to see that $\mathbb{R}\ni s\mapsto \overline{f}(x,s)$ is continuous, using again \textbf{(H9)}. Thus $\overline{f}$ is a Carath\'{e}odory function.

			\noindent \textbf{(2)} For almost all $x\in\Omega$ the function:
			
			\begin{equation}
				\mathbb{R}\ni s\mapsto\overline{f}(x,s)+\lambda\cdot\overline{b}(x,s)=\begin{cases} f(x,0)+\lambda b(x,0)+\left (\lambda-\dfrac{\lambda_0}{2}\right )s, & s<0\\[3mm] f(x,s)+\lambda_0 b(x,s)+(\lambda-\lambda_0) b(x,s), & s\in [0,\delta_0]\\[3mm] f(x,\delta_0)+\lambda b(x,\delta_0)+(\lambda-\tilde{\lambda}_0)(s-\delta_0), & s>\delta_0 \end{cases}
			\end{equation}
			
			\noindent is continuous and strictly increasing on each branch, from \textbf{(H12)} and \textbf{(H13)}. Here we have used that $\tilde{\lambda}_0<\lambda_0$. This proves the assertion.
			
			\noindent\textbf{(3)} For $V\in L^2(\Omega)$ we have that:
			
			\begin{align}
				|\overline{f}(x,V(x))|&=\begin{cases} \left|f(x,0)-\dfrac{\lambda_0}{2}V(x)\right |, & V(x)<0\\[3mm] |f(x,V(x))|, & V(x)\in [0,\delta_0] \\[3mm] |f(x,\delta_0)-\tilde{\lambda}_0 (V(x)-\delta_0)|, & V(x)>\delta_0\end{cases}\nonumber\\ 
			\eqref{fcevaimportant}\ \ \ \ \ \ \ \	&\leq \begin{cases} f(x,0)+\dfrac{\lambda_0}{2}|V(x)|, & V(x)<0\\[3mm] \lambda_0\big (b(x,\delta_0)-b(x,0) \big ), & V(x)\in [0,\delta_0] \nonumber\\[3mm] -f(x,\delta_0)+\tilde{\lambda}_0|V(x)|-\delta_0\tilde{\lambda}_0, & V(x)>\delta_0\end{cases}\\
				&\leq f(x,0)-f(x,\delta_0)+\lambda_0|V(x)|+\lambda_0\big (b(x,\delta_0)-b(x,0) \big )\in L^2(\Omega)\label{ineqf}.
			\end{align}
			
			\noindent So $\overline{f}(\cdot,V(\cdot))\in L^2(\Omega)$. Choosing $V\equiv s\in\mathbb{R}$ gives us the desired inequality.
			
				\noindent\textbf{(4)} From Proposition \ref{prop35} \textbf{(3)} we have that $f(\cdot,0),f(\cdot,\delta_0)\in L^{\infty}(\Omega)$. We also know that $b(\cdot,0),b(\cdot,\delta_0)\in L^{\infty}(\Omega)$ and $V\in L^{\infty}(\Omega)$. Therefore \eqref{ineqf} implies that $\overline{f}(\cdot,V(\cdot))\in L^{\infty}(\Omega)$.
			
				\noindent\textbf{(5)} Since for a.e. $x\in\Omega$, $\overline{f}(x,\cdot)$ is continuous ($f$ being a Carath\'{e}odory function) we deduce that its antiderivative $\overline{F}(x,\cdot)$ is continuously differentiable on $\mathbb{R}$.
			Thus, from the chain rule, for any $s>0$, we get that: $\dfrac{\partial }{\partial s}\overline{F}(x,\sqrt[\alpha]{s})=\overline{f}(x,\sqrt[\alpha]{s})\dfrac{1}{\alpha\sqrt[\alpha]{s^{\alpha-1}}}=\dfrac{1}{\alpha}\begin{cases}\dfrac{f(x,\sqrt[\alpha]{s})}{\sqrt[\alpha]{s^{\alpha-1}}}, & s\in (0,\delta_0^\alpha]\\[3mm] \dfrac{f(x,\delta_0)-\tilde{\lambda}_0(\sqrt[\alpha]{s}-\delta_0)}{\sqrt[\alpha]{s^{\alpha-1}}}, & s>\delta_0^\alpha \end{cases}$. From \textbf{(EH)} we get that this function is decreasing on $(0,\delta_0^\alpha]$. Being continuous at $s=\delta_0^\alpha$ we only need to check if the second branch defines a decreasing function too. Note that for any $s>\delta_0^\alpha$ we have:
			
			\begin{equation}
				\dfrac{\partial}{\partial s}\left [ \dfrac{f(x,\delta_0)-\tilde{\lambda}_0(\sqrt[\alpha]{s}-\delta_0)}{\sqrt[\alpha]{s^{\alpha-1}}}\right ]=\dfrac{s^{\frac{1}{\alpha}-2}}{\alpha}\big [\underbrace{(1-\alpha)}_{<0}\big (f(x,\delta_0)+\tilde{\lambda}_0\delta_0 \big ) -\tilde{\lambda}_0\underbrace{(2-\alpha)}_{\geq 0}\sqrt[\alpha]{s}\big ]\leq 0,
			\end{equation}
			
			\noindent because $-\tilde{\lambda}_0(2-\alpha)\sqrt[\alpha]{s}\leq-\tilde{\lambda}_0(2-\alpha)\delta_0\leq (\alpha-1)\big (f(x,\delta_0)+\tilde{\lambda}_0\delta_0)\Longleftrightarrow \tilde{\lambda}_0\geq -\dfrac{\alpha-1}{\delta_0}f(x,\delta_0)$, which is true from the statement.
			
			\noindent  We conclude that $-\dfrac{\partial }{\partial s}\overline{F}(x,\sqrt[\alpha]{s})\geq 0$ for any $s>0$. Therefore $-\overline{F}(x,\sqrt[\alpha]{s})$ is convex with respect to $s\in [0,\infty)$.
			
			\noindent\textbf{(6)} In a similar manner as above we obtain $-\dfrac{\partial }{\partial s}\overline{F}(x,\sqrt[\alpha]{s})>0$ for any $s>0$. It is essential that $\alpha<2$ and $\tilde{\lambda}_0>0$. Therefore $-\overline{F}(x,\sqrt[\alpha]{s})$ is strictly convex with respect to $s\in [0,\infty)$.
			
		\end{proof}
		
		\begin{proposition}\label{prop37} The function $\overline{F}:\Omega\times\mathbb{R}\to\mathbb{R}$ has the following properties:
			
			\begin{enumerate}
				\item[\textnormal{\textbf{(1)}}] $\overline{F}$ is a Carath\'{e}odory function. In particular its restriction $F$ is a Carath\'{e}odory function.

				\item[\textnormal{\textbf{(2)}}] For a.e. $x\in\Omega$ we have that $\overline{F}(x,\cdot)\in C^1(\mathbb{R})$ and $\dfrac{\partial\overline{F}}{\partial s}(x,s)=\overline{f}(x,s)$ for any $s\in\mathbb{R}$. In particular, for a.e. $x\in\Omega$ we have that $F(x,\cdot)\in C^1([0,\delta_0])$ and $\dfrac{\partial F}{\partial s}(x,s)=f(x,s)$ for any $s\in [0,\delta_0]$.
				
				\item[\textnormal{\textbf{(3)}}] For any $V\in L^2(\Omega)$ we have that $\overline{F}\big (\cdot,V(\cdot)\big )\in L^1(\Omega)$.
				
				\item[\textnormal{\textbf{(4)}}] For a.e. $x\in\Omega$ and any $s\in\mathbb{R}$ the following inequality holds:
				
				\begin{equation}
					\overline{F}(x,s)\leq \dfrac{f(x,0)^2}{\lambda_0}+\dfrac{f(x,\delta_0)^2}{2\tilde{\lambda}_0}+\delta_0\lambda_0 (b(x,\delta_0)-b(x,0)).
				\end{equation}
			\end{enumerate}
			
		\end{proposition}
	
\begin{proof} \textbf{(1)} and \textbf{(2)} Since from Proposition \ref{prop35} \textbf{(5)} we know that $f\in L^2\bigl((0,\delta_0)\times\Omega\bigr)$, we can use \textit{Fubini's theorem} given as \cite[Proposition 5.2.2, page 148]{Cohn} in order to obtain that for any fixed $s\in [0,\delta_0]$ the function $\Omega\ni x\mapsto\overline{F}(x,s)=\displaystyle\int_{0}^s\overline{f}(x,\tau)\ d\tau=\displaystyle\int_{0}^s f(x,\tau)\ d\tau$ is measurable (in fact from $L^1(\Omega)$). In particular this shows that $F(\cdot,\delta_0)$ is measurable. For $s<0$ the function $\Omega\ni x\mapsto f(x,0)s-\dfrac{\lambda_0}{4}s^2$ is clearly measurable. Similarly for $s>\delta_0$ the function $\Omega\ni x\mapsto F(x,\delta_0)+f(x,\delta_0)(s-\delta_0)-\dfrac{\tilde{\lambda}_0}{2}(s-\delta_0)^2$ is also measurable, since from Proposition \ref{prop35} \textbf{(2)}, $f(\cdot,\delta_0)$ is measurable. Thus for any $s\in\mathbb{R}$ we proved that $\overline{F}(\cdot,s)$ is measurable. Lastly, from the \textit{Fundamental theorem of calculus}, because $\mathbb{R}\ni s\mapsto\overline{f}(x,s)$ is continuous for a.e. $x\in\Omega$,  we get that $\overline{F}(x,\cdot)\in C^1(\mathbb{R})\subset C(\mathbb{R})$ and $\dfrac{\partial\overline{F}}{\partial s}(x,s)=\overline{f}(x,s)$ for any $s\in\mathbb{R}$.

\bigskip

\noindent\textbf{(3)} Using Proposition \ref{prop32} \textbf{(2)} and Proposition \ref{prop36} \textbf{(2)}, \textbf{(3)}, we have for a.e. $x\in\Omega$ that:

\begin{align*}
	&\overline{F}(x,V(x))+\lambda_0\overline{B}(x,V(x))=\\
	&=\int_{0}^{V(x)}\overline{f}(x,\tau)+\lambda_0\overline{b}(x,\tau)\ d\tau\ \begin{cases} \leq V(x)\underbrace{\big (\overline{f}(x,V(x))+\lambda_0\overline{b}(x,V(x))\big )}_{\in L^2(\Omega)}\in L^1(\Omega) \\[3mm] \geq V(x)\underbrace{\big (f(x,0)+\lambda_0 b(x,0)\big )}_{\in L^2(\Omega)}\in L^1(\Omega)\end{cases},
\end{align*}

\noindent from \textit{Cauchy inequality}. Being between two $L^1(\Omega)$ functions we deduce that $\overline{F}(\cdot,V(\cdot))+\lambda_0\overline{B}(\cdot,V(\cdot))\in L^1(\Omega)$. But, from Proposition \ref{prop33} \textbf{(3)} we already have that $\overline{B}(\cdot,V(\cdot))\in L^1(\Omega)$. By difference we get that $\overline{F}(\cdot,V(\cdot))\in L^1(\Omega)$, as needed.

\noindent\textbf{(4)} For a.e. $x\in\Omega$ and for all $s\in\mathbb{R}$ we have that:

\begin{align*}
	\overline{F}(x,s)&=\begin{cases} f(x,0)s-\dfrac{\lambda_0}{4} s^2, & s<0 \\[3mm] \displaystyle\int_{0}^s f(x,\tau)\ d\tau, &  s\in [0,\delta_0] \\[3mm] \displaystyle\int_{0}^{\delta_0} f(x,\tau)\ d\tau+f(x,\delta_0)(s-\delta_0)-\dfrac{\tilde{\lambda}_0}{2}(s-\delta_0)^2, & s>\delta_0 \end{cases}\\[3mm]
	&\leq \begin{cases} \dfrac{f(x,0)^2}{\lambda_0}-\dfrac{\lambda_0}{4}\left (s-\dfrac{2f(x,0)}{\lambda_0} \right )^2, & s<0 \\[3mm] \displaystyle\int_{0}^s \lambda_0 (b(x,\delta_0)-b(x,0))\ d\tau, &  s\in [0,\delta_0] \\[3mm] \displaystyle\int_{0}^{\delta_0} \lambda_0 (b(x,\delta_0)-b(x,0))\ d\tau+f(x,\delta_0)(s-\delta_0)-\dfrac{\tilde{\lambda}_0}{2}(s-\delta_0)^2, & s>\delta_0 \end{cases}\\[3mm]
	&\leq  \begin{cases} \dfrac{f(x,0)^2}{\lambda_0}, & s<0 \\[3mm]  \delta_0\lambda_0 (b(x,\delta_0)-b(x,0)), &  s\in [0,\delta_0] \\[3mm]  \delta_0\lambda_0 (b(x,\delta_0)-b(x,0))+\dfrac{f(x,\delta_0)^2}{2\tilde{\lambda}_0}-\dfrac{\tilde{\lambda}_0}{2}\left (s-\delta_0-\dfrac{f(x,\delta_0)}{\tilde{\lambda}_0}\right )^2, & s>\delta_0 \end{cases}\\[3mm]
	&\leq \dfrac{f(x,0)^2}{\lambda_0}+\dfrac{f(x,\delta_0)^2}{2\tilde{\lambda}_0}+\delta_0\lambda_0 (b(x,\delta_0)-b(x,0)).
\end{align*}

\end{proof}

			\begin{proposition}\label{prop38} The functional $\overline{\mathcal{F}}:W^{1,p(x)}(\Omega)\to\mathbb{R}$ has the following properties:
			
			\begin{enumerate}
			\item[\textnormal{\textbf{(1)}}] $\overline{\mathcal{F}}$ is well-defined.
			
			\item[\textnormal{\textbf{(2)}}] $\overline{\mathcal{F}}(V)=\mathcal{F}(V)$ for every $V\in X$.
			
			\item[\textnormal{\textbf{(3)}}] $\displaystyle\inf_{V\in W^{1,p(x)}(\Omega)} -\overline{\mathcal{F}}(V)>-\infty$, i.e. $-\overline{\mathcal{F}}$ is bounded from below, or $\overline{\mathcal{F}}$ is bounded from above.
			
			\item[\textnormal{\textbf{(4)}}] $\displaystyle\overline{\mathcal{F}}\in C^1\big (W^{1,p(x)}(\Omega)\big )$ and for any $V\in W^{1,p(x)}(\Omega)$:
			
			\begin{equation}
				\langle\overline{\mathcal{F}}'(V),\phi\rangle=\int_{\Omega} \overline{f}\big (x,V(x)\big )\phi(x)\ dx,\ \forall\  \phi\in W^{1,p(x)}(\Omega).
			\end{equation}
			\end{enumerate}
			
		\end{proposition}
		
		\begin{proof} \textbf{(1)} For any $V\in W^{1,p(x)}(\Omega)$, from \textbf{(H2)} we get that $V\in L^2(\Omega)$. From Proposition \ref{prop37} \textbf{(3)} we have that $\overline{F}(\cdot,V(\cdot))\in L^1(\Omega)$ and therefore $\overline{\mathcal{F}}(V)=\displaystyle\int_{\Omega}\overline{F}(x,V(x))\ dx\in\mathbb{R}$.
			
		\noindent\textbf{(2)} $\overline{\mathcal{F}}(V)=\displaystyle\int_{\Omega}\overline{F}(x,V(x))\ dx=\displaystyle\int_{\Omega} F(x,V(x))\ dx=\mathcal{F}(V)$, because $V(x)\in [0,\delta_0]$ for a.e. $x\in\Omega$.

		\noindent\textbf{(3)} For any $V\in W^{1,p(x)}(\Omega)$, using Proposition \ref{prop37} \textbf{(4)}, we get that:
		
		\begin{equation}
			\overline{\mathcal{F}}(V)=\int_{\Omega} \overline{F}(x,V(x))\ dx\leq \int_{\Omega}\dfrac{f(x,0)^2}{\lambda_0}+\dfrac{f(x,\delta_0)^2}{2\tilde{\lambda}_0}+\delta_0\lambda_0 (b(x,\delta_0)-b(x,0)) \ dx<\infty.
		\end{equation}

		\noindent\textbf{(4)} \noindent Next, we'll show that $\overline{\mathcal{F}}$ is G\^ateaux-differentiable on $W^{1,p(x)}(\Omega)$. Set a direction $\phi\in W^{1,p(x)}(\Omega)$ (so $\phi\in L^{2}(\Omega)$) and consider any sequence $(\varepsilon_n)_{n\geq 1}\subset\mathbb{R}^*$ that is convergent to $0$. Without loss of generality we can assume that $|\varepsilon_n|\leq 1,\ \forall\ n\geq 1$. 
		
		\noindent Consider for each $n\geq 1$ the function: $g_n:\Omega\to\mathbb{R},\ g_n(x)=\dfrac{\overline{F}(x,V(x)+\varepsilon_n\phi(x))-\overline{F}(x,V(x))}{\varepsilon_n}$. From Proposition \ref{prop37} \textbf{(2)} and the classical 1-dimensional \textit{mean value theorem} we have that $g_n(x)=\dfrac{\partial \overline{F}}{\partial s}(x,V(x)+s_n(x)\phi(x))\phi(x)=\overline{f}(x,V(x)+s_n(x)\phi(x))\phi(x)$ for some real number $s_n(x)$ with $|s_n(x)|\leq\varepsilon_n$. Using the continuity of $\overline{f}(x,\cdot)$ for a.e. $x\in\Omega$ (see Proposition \ref{prop36} \textbf{(1)}) we get that:
		
		\begin{equation}
			\lim\limits_{n\to\infty} g_n(x)=\lim\limits_{n\to\infty}\overline{f}(x,V(x)+s_n(x)\phi(x))\phi(x)\stackrel{|s_n(x)|\leq \varepsilon_n}{=}\overline{f}(x,V(x))\phi(x),\ \text{for a.e.}\ x\in\Omega.
		\end{equation}
		
		\noindent This proves that $g_n$ has a pointwise limit. The next step will be to show that $|g_n|$ is bounded by a function from $L^1(\Omega)$ for any $n\geq 1$. Indeed, from Proposition \ref{prop36} \textbf{(3)} we get that:
		
		\begin{align*}
			|g_n(x)|&=|\overline{f}(x,V(x)+s_n(x)\phi(x))|\cdot |\phi(x)|\leq \\
			&\leq \big [\lambda_0\big (b(x,\delta_0)-b(x,0)\big )+\lambda_0|V(x)+s_n(x)\phi(x)|+f(x,0)-f(x,\delta_0) \big ]\cdot |\phi(x)|\\
			&\leq \big [\lambda_0\big (b(x,\delta_0)-b(x,0)\big )+\lambda_0|V(x)|+\lambda_0|\phi(x)|+f(x,0)-f(x,\delta_0) \big ]\cdot |\phi(x)|\\
			\textbf{(H11)}\ \ \ \	&\leq \underbrace{\big [\lambda_0\big (b(x,\delta_0)-b(x,0)\big )+\lambda_0|V(x)|+\lambda_0|\phi(x)|+f(x,0)-f(x,\delta_0) \big ]}_{\in L^2(\Omega)}\cdot \underbrace{|\phi(x)|}_{\in L^2(\Omega)}\in L^1(\Omega),\\
		\end{align*}
		
		\noindent because $V,\phi\in W^{1,p(x)}(\Omega)\hookrightarrow L^2(\Omega)$. We have also used here Proposition \ref{prop35} \textbf{(2)}.
		
		\noindent From the \textit{Lebesgue dominated convergence theorem} we get that:
		
		\begin{align*}
			\lim\limits_{n\to \infty}\dfrac{\overline{\mathcal{F}}(V+\varepsilon_n\phi)-\overline{\mathcal{F}}(V)}{\varepsilon_n}&=\lim\limits_{n\to\infty} \int_{\Omega}\dfrac{\overline{F}(x,V(x)+\varepsilon_n\phi(x))-\overline{F}(x,V(x))}{\varepsilon_n}\ dx\\
			&=\lim\limits_{n\to\infty}\int_{\Omega} g_n(x)\ dx=\int_{\Omega} \overline{f}(x,V(x))\phi(x)\ dx.
		\end{align*}
		
		\noindent It is easy to remark that $\partial\overline{\mathcal{F}}(V):W^{1,p(x)}(\Omega)\to\mathbb{R},\ \partial_{\phi}\overline{\mathcal{F}}(V)=\displaystyle\int_{\Omega} \overline{f}(x,V(x))\phi(x)\ dx$ is a linear operator. We show now that it is in fact a bounded linear operator. Indeed:
		
		\begin{align*}
			\left |\partial_{\phi}\overline{\mathcal{F}}(V) \right |&=\left | \int_{\Omega} \overline{f}(x,V(x))\phi(x)\ dx \right |\leq \int_{\Omega} |\overline{f}(x,V(x))|\cdot |\phi(x)|\ dx\\
			\text{Proposition \ref{prop36}\ \bf{(3)}}\ \ \ \ \	&\leq \int_{\Omega} \left ( \lambda_0\big (b(x,\delta_0)-b(x,0)\big )+\lambda_0|V(x)|+f(x,0)-f(x,\delta_0) \right )|\phi(x)|\ dx\\
			\text{(Cauchy ineq.)}	\ \ \ \ \ &\leq \big\Vert\lambda_0\big (b(x,\delta_0)-b(x,0)\big )+\lambda_0|V(x)|+f(x,0)-f(x,\delta_0)\big \Vert_{L^2(\Omega)}\Vert \phi\Vert_{L^2(\Omega)}\\
			\big (W^{1,p(x)}(\Omega)\hookrightarrow L^{2}(\Omega)\big )\ \ \ \ \ \ \ &\lesssim \big\Vert \lambda_0\big (b(x,\delta_0)-b(x,0)\big )+\lambda_0|V(x)|+f(x,0)-f(x,\delta_0)\big \Vert_{L^2(\Omega)}\Vert \phi\Vert_{W^{1,p(x)}(\Omega)}.
		\end{align*}
		
		\noindent Henceforth: $\displaystyle\sup_{\phi\in W^{1,p(x)}(\Omega)\setminus\{0\}}\dfrac{\big |\partial_{\phi}\overline{\mathcal{F}}(V) \big |}{\Vert \phi\Vert_{W^{1,p(x)}(\Omega)}}<\infty$, which shows that $\partial\overline{\mathcal{F}}(V)\in\mathcal{L}(W^{1,p(x)}(\Omega);\mathbb{R})=W^{1,p(x)}(\Omega)^*$, for any $V\in W^{1,p(x)}(\Omega)$. This observation allows us to write that: $\partial\overline{\mathcal{F}}:W^{1,p(x)}(\Omega)\to W^{1,p(x)}(\Omega)^*$. The following step is to prove that $\partial\overline{\mathcal{F}}$ is a continuous mapping.

		\noindent Consider a sequence $(V_n)_{n\geq 1}\subset W^{1,p(x)}(\Omega)$ and some $V\in W^{1,p(x)}(\Omega)$ such that $V_n\to V$ in $W^{1,p(x)}(\Omega)$. Since $W^{1,p(x)}(\Omega)\hookrightarrow L^{2}(\Omega)$ we also have that $V_n\to V$ in $L^{2}(\Omega)$.

		\begin{align*}\big \Vert\partial\overline{\mathcal{F}}(V_n)-\partial\overline{\mathcal{F}}(V)\big\Vert_{W^{1,p(x)}(\Omega)^*}&=\sup_{\phi\in W^{1,p(x)}(\Omega)\setminus\{0\}}\dfrac{\big |\partial_{\phi}\overline{\mathcal{F}}(V_n)-\partial_{\phi}\overline{\mathcal{F}}(V)\big |}{\Vert\phi\Vert_{W^{1,p(x)}(\Omega)}}\\
			&=\sup_{\phi\in W^{1,p(x)}(\Omega)\setminus\{0\}}\dfrac{\left |\displaystyle\int_{\Omega} \big (\overline{f}(x,V_n(x))-\overline{f}(x,V(x)))\cdot\phi(x) dx\right |}{\Vert\phi\Vert_{W^{1,p(x)}(\Omega)}} \\
			&\leq \sup_{\phi\in W^{1,p(x)}(\Omega)\setminus\{0\}}\dfrac{\displaystyle\int_{\Omega} \big |\overline{f}(x,V_n(x))-\overline{f}(x,V(x))|\cdot|\phi(x)| dx}{\Vert\phi\Vert_{W^{1,p(x)}(\Omega)}}\\
			\text{(Cauchy ineq.)}\ \ \ \ \ \ \ \	&\leq \sup_{\phi\in W^{1,p(x)}(\Omega)\setminus\{0\}}\dfrac{\Vert \overline{f}(\cdot,V_n)-\overline{f}(\cdot,V) \Vert_{L^{2}(\Omega)}\Vert \phi\Vert_{L^{2}(\Omega)}}{\Vert\phi\Vert_{W^{1,p(x)}(\Omega)}}\\
			\big (W^{1,p(x)}(\Omega)\hookrightarrow L^{2}(\Omega)\big )\ \ \ \ \ \ \ &\lesssim \sup_{\phi\in W^{1,p(x)}(\Omega)\setminus\{0\}}\dfrac{\Vert \overline{f}(\cdot,V_n)-\overline{f}(\cdot,V)\Vert_{L^{2}(\Omega)}\Vert \phi\Vert_{W^{1,p(x)}(\Omega)}}{\Vert\phi\Vert_{W^{1,p(x)}(\Omega)}}\\
			&=\Vert \overline{f}(\cdot,V_n)-\overline{f}(\cdot,V)\Vert_{L^{2}(\Omega)}\stackrel{n\to\infty}{\longrightarrow} 0,
		\end{align*}

		\noindent from Proposition \ref{prop36} \textbf{(3)}.
		\noindent Finally, from Theorem \ref{athmgatfre}, we obtain that $\overline{\mathcal{F}}\in C^1\big (W^{1,p(x)}(\Omega)\big )$ and:
		
		\begin{equation*}
			\langle \overline{\mathcal{F}}'(V),\phi \rangle=\displaystyle\int_{\Omega} \overline{f}(x,V(x))\phi(x)\ dx,\ \forall\ V,\phi\in W^{1,p(x)}(\Omega).
		\end{equation*}

		\end{proof}

		\begin{proposition}\label{prop39} The functional $\overline{\mathcal{C}}:W^{1,p(x)}(\Omega)\to\mathbb{R},\ \overline{\mathcal{C}}(V)=\displaystyle\int_{\Omega} \dfrac{|V(x)|^{p(x)}}{p(x)}\ dx$ is well-defined, $\mathcal{C}\in C^1\big (W^{1,p(x)}(\Omega) \big )$, and:
			
			\begin{equation}
				\langle \mathcal{C}'(V),\phi\rangle=\int_{\Omega} |V(x)|^{p(x)-2}V(x)\phi(x)\ dx,\ \forall V,\phi\in W^{1,p(x)}(\Omega).
			\end{equation}
			
		\end{proposition}
		
		\begin{proof} First remark that $0\leq \mathcal{C}(V)\leq \dfrac{1}{p^-}\int_{\Omega} |V(x)|^{p(x)}\ dx<\infty$ from the definition of the Sobolev space $W^{1,p(x)}(\Omega)\ni V$. So $\mathcal{C}$ is well-defined.
			
		\noindent For a fixed $V\in W^{1,p(x)}(\Omega)$ we consider the linear operator $\ell_V:W^{1,p(x)}(\Omega)\to\mathbb{R}$ given by:
		
		\begin{equation}
			\ell_V(\phi)=\int_{\Omega} |V(x)|^{p(x)-2}V(x)\phi(x)\ dx,\ \forall\ \phi\in W^{1,p(x)}(\Omega).
		\end{equation}
		
		\noindent The first think to observe is that $\ell_V\in W^{1,p(x)}(\Omega)^*$. Indeed, for all $\phi\in W^{1,p(x)}(\Omega)$:
		
		\begin{equation*}
			\left | \ell_V(\phi) \right |\leq \int_{\Omega} |V|^{p(x)-1}\phi(x)\ dx\leq 2\Vert V^{p(x)-1}\Vert_{L^{p'(x)}(\Omega)}\Vert\phi\Vert_{L^{p(x)}(\Omega)}\leq 2\Vert V^{p(x)-1}\Vert_{L^{p'(x)}(\Omega)}\Vert\phi\Vert_{W^{1,p(x)}(\Omega)}.
		\end{equation*}
		
		\noindent So $\displaystyle\sup_{\phi\in W^{1,p(x)}(\Omega)\setminus\{0\}}\dfrac{|\ell_V(\phi)|}{\Vert\phi\Vert_{W^{1,p(x)}(\Omega)}}\leq 2\Vert V^{p(x)-1}\Vert_{L^{p'(x)}(\Omega)}\ \Rightarrow\ \ell_V\in W^{1,p(x)}(\Omega)^*$.

		\noindent Fix some $V\in W^{1,p(x)}(\Omega)$ and a direction $\phi\in W^{1,p(x)}(\Omega)$. For any sequence of real non-zero numbers $(\varepsilon_n)_{n\geq 1}$ convergent to $0$ we define the sequence of function $(v_n)_{n\geq 1}$ by:
		
		\begin{equation}
			v_n(x)=\dfrac{|V(x)+\varepsilon_n\phi(x)|^{p(x)}-|V(x)|^{p(x)}}{p(x)\varepsilon_n},\ \text{for a.e.}\ x\in\Omega,\ \forall\ n\geq 1.
		\end{equation}
		
		\noindent Using Lemma \ref{alemhos} ($p(x)\geq p^->1$ for all $x\in\overline{\Omega}$) we have that: $\lim\limits_{n\to\infty} v_n(x)=|V(x)|^{p(x)-2}V(x)\phi(x)$ for a.e. $x\in\Omega$.

		\noindent Taking into account that both $V+\varepsilon_n\phi, \phi\in W^{1,p(x)}(\Omega)\subset L^{p(x)}(\Omega)$ we deduce from Theorem \ref{athleb} \textbf{(24)} that $|V+\varepsilon\phi|^{p(x)},|V(x)|^{p(x)}\in L^1(\Omega)$. Thence $v_n\in L^1(\Omega)$ for each integer $n\geq 1$.
		
		\noindent From the \textit{Lebesgue dominated convergence theorem} we obtain that:
		
		\begin{align*}
			\lim\limits_{n\to\infty}\dfrac{\mathcal{C}(V+\varepsilon_n\phi)-\mathcal{C}(V)}{\varepsilon_n}&=\lim\limits_{n\to\infty}\int_{\Omega}\dfrac{|V(x)+\varepsilon_n\phi(x)|^{p(x)}-|V(x)|^{p(x)}}{p(x)\varepsilon_n}\ dx\\
			&=\lim\limits_{n\to\infty}\int_{\Omega} v_n(x)\ dx=\int_{\Omega} |V(x)|^{p(x)-2}V(x)\phi(x)\ dx\\
			&=\ell_V(\phi),\ \forall\ V,\phi\in W^{1,p(x)}(\Omega).
		\end{align*}
		
		\noindent So $\partial\mathcal{C}(V):W^{1,p(x)}(\Omega)\to\mathbb{R},\ \partial\mathcal{C}(V)\phi=\ell_{V}(\phi)$ is a bounded linear operator as we have seen before. Thus we can define the G\^ateaux differential of $\mathcal{C}$ as $\partial\mathcal{C}:W^{1,p(x)}(\Omega)\to W^{1,p(x)}(\Omega)^*$.
		
		\noindent Next, we need to show that $\partial\mathcal{C}$ is a continuous (nonlinear) operator. In that sense consider some $V\in W^{1,p(x)}(\Omega)$ and any sequence $(V_n)_{n\geq 1}\subset W^{1,p(x)}(\Omega)$ with $V_n\longrightarrow V$ in $W^{1,p(x)}(\Omega)$. Also denote: $\Omega_1=\{x\in\overline{\Omega}\ |\ p(x)<2\}$ and $\Omega_2=\{x\in\overline{\Omega}\ |\ p(x)\geq 2\}$. Therefore:
		
		\begin{align*}
			&\Vert \partial \mathcal{C}(V_n)-\partial\mathcal{C}(V)\Vert_{W^{1,p(x)}(\Omega)^*}=\sup_{\phi\in W^{1,p(x)}(\Omega)\setminus\{0\}}\dfrac{\left |\displaystyle\int_{\Omega} \big (|V_n(x)|^{p(x)-2}V_n(x)-|V(x)|^{p(x)-2}V(x)\big )\phi(x) \ dx\right |}{\Vert\phi\Vert_{W^{1,p(x)}(\Omega)}}\\
			&\leq \sup_{\phi\in W^{1,p(x)}(\Omega)\setminus\{0\}}\dfrac{\displaystyle\int_{\Omega} \left | |V_n(x)|^{p(x)-2}V_n(x)-|V(x)|^{p(x)-2}V(x)\right |\cdot |\phi(x)| \ dx}{\Vert\phi\Vert_{W^{1,p(x)}(\Omega)}}\\
			&\leq\sup_{\phi\in W^{1,p(x)}(\Omega)\setminus\{0\}}\dfrac{\displaystyle\int_{\Omega_1} \left | |V_n(x)|^{p(x)-2}V_n(x)-|V(x)|^{p(x)-2}V(x)\right |\cdot |\phi(x)| \ dx}{\Vert\phi\Vert_{W^{1,p(x)}(\Omega)}}+\\
			&+\dfrac{\displaystyle\int_{\Omega_2} \left | |V_n(x)|^{p(x)-2}V_n(x)-|V(x)|^{p(x)-2}V(x)\right |\cdot |\phi(x)| \ dx}{\Vert\phi\Vert_{W^{1,p(x)}(\Omega)}}\\
			&\leq\sup_{\phi\in W^{1,p(x)}(\Omega)\setminus\{0\}}\dfrac{c_1\displaystyle\int_{\Omega_1} |V_n-V|^{p(x)-1}|\phi(x)|\ dx}{\Vert\phi\Vert_{W^{1,p(x)}(\Omega)}}+\dfrac{c_2\displaystyle\int_{\Omega_2}\big (|V_n|+|V| \big )^{p(x)-2} |V_n-V|\cdot|\phi(x)|\ dx}{\Vert\phi\Vert_{W^{1,p(x)}(\Omega)}},\\
			%\text{(H\"{o}lder)} &\leq \sup_{\phi\in W^{1,p(x)}(\Omega)\setminus\{0\}}
		\end{align*}

		\noindent from \text{(Lemma \ref{alplap})}.

		\noindent We have here two separated problems. From $|V_n-V|\in W^{1,p(x)}(\Omega)\subset L^{p(x)}(\Omega)$ we deduce from Theorem \ref{athleb} \textbf{(24)} that $|V_n-V|^{p(x)-1}\in L^{p'(x)}(\Omega)$ for each $n\geq 1$. Since $|\phi|\in W^{1,p(x)}(\Omega)\subset L^{p(x)}(\Omega)$ we get from \textit{H\"{o}lder inequality} that:
		
		\begin{align*}
			\dfrac{1}{\Vert\phi\Vert_{W^{1,p(x)}(\Omega)}}\int_{\Omega_1} |V_n-V|^{p(x)-1}|\phi(x)|\ dx &\leq \dfrac{2}{\Vert\phi\Vert_{W^{1,p(x)}(\Omega)}}\Vert |V_n-V|^{p(x)-1}\Vert_{L^{p'(x)}(\Omega_1)}\Vert\phi\Vert_{L^{p(x)}(\Omega_1)}\\
			&\leq \dfrac{2}{\Vert\phi\Vert_{W^{1,p(x)}(\Omega)}}\Vert |V_n-V|^{p(x)-1}\Vert_{L^{p'(x)}(\Omega)}\Vert\phi\Vert_{L^{p(x)}(\Omega)}\\
			&\leq  2\Vert |V_n-V|^{p(x)-1}\Vert_{L^{p'(x)}(\Omega)}\stackrel{n\to\infty}{\longrightarrow} 0\\
			\text{(Theorem 2.58 in \cite{Cruz})}\ 	&\Longleftrightarrow\rho_{p'(x)}\left (|V_n-V|^{p(x)-1} \right )\stackrel{n\to\infty}{\longrightarrow} 0\\
			&\Longleftrightarrow\int_{\Omega}\left (|V_n-V|^{p(x)-1} \right )^{p'(x)}\ dx\stackrel{n\to\infty}{\longrightarrow} 0\\
			&\Longleftrightarrow\int_{\Omega}|V_n-V|^{p(x)}\ dx\stackrel{n\to\infty}{\longrightarrow} 0\\
			&\Longleftrightarrow\rho_{p(x)}\left (V_n-V \right )\stackrel{n\to\infty}{\longrightarrow} 0\\
			\text{(Theorem 2.58 in \cite{Cruz})}\ 	&\Longleftrightarrow\Vert V_n-V\Vert_{L^{p(x)}(\Omega)}\stackrel{n\to\infty}{\longrightarrow} 0,
		\end{align*}
		
		\noindent which is true, because $V_n\to V$ in $W^{1,p(x)}(\Omega)$ so that: $0\leq \Vert V_n-V\Vert_{L^{p(x)}(\Omega)}\leq\Vert V_n-V\Vert_{W^{1,p(x)}(\Omega)} \stackrel{n\to\infty}{\longrightarrow} 0$. 
		
		\noindent For the other term notice that $z_n:=|V_n|+|V|\in W^{1,p(x)}(\Omega)\subset L^{p(x)}(\Omega)$, and therefore, from Theorem \ref{athleb} \textbf{(24)}, we have that $z_n^{p(x)-2}\in L^{\frac{p(x)}{p(x)-2}}(\Omega)$ for each $n\geq 1$. Another important remark is that $z_n\to 2|V|$ in $L^{p(x)}(\Omega)$ because $\big |z_n-2|V|\big |=\big ||V_n|-|V|\big |\leq |V_n-V|$, which implies that $0\leq\rho_{p(x)}\big (z_n-2|V| \big )\leq \rho_{p(x)}\big (V_n-V \big )\stackrel{n\to\infty}{\longrightarrow} 0$.\footnote{See Theorem \ref{athleb} \textbf{(14)}.}
		
		\noindent From the \textit{Generalized H\"{o}lder inequality}\footnote{It can be found at \cite[Corollary 2.30, page 30]{Cruz}.}, taking into account that $\left (\dfrac{p(x)}{p(x)-2} \right )^{-1}+\dfrac{1}{p(x)}+\dfrac{1}{p(x)}=1$ we obtain that:
		
		\begin{align*}
			\dfrac{\displaystyle\int_{\Omega_2}\big (|V_n|+|V| \big )^{p(x)-2} |V_n-V|\cdot|\phi(x)|\ dx}{\Vert\phi\Vert_{W^{1,p(x)}(\Omega)}} &\lesssim \Vert z_n^{p(x)-2}\Vert_{L^{\frac{p(x)}{p(x)-2}}(\Omega_2)}\Vert V_n-V\Vert_{L^{p(x)}(\Omega_2)}\dfrac{\Vert\phi\Vert_{L^{p(x)}(\Omega_2)}}{\Vert\phi\Vert_{W^{1,p(x)}(\Omega)}}      \\
			&\leq\Vert z_n^{p(x)-2}\Vert_{L^{\frac{p(x)}{p(x)-2}}(\Omega)}\Vert V_n-V\Vert_{L^{p(x)}(\Omega)}\dfrac{\Vert\phi\Vert_{L^{p(x)}(\Omega)}}{\Vert\phi\Vert_{W^{1,p(x)}(\Omega)}} \\
			&\leq \Vert z_n^{p(x)-2}\Vert_{L^{\frac{p(x)}{p(x)-2}}(\Omega)}\Vert V_n-V\Vert_{L^{p(x)}(\Omega)}\stackrel{n\to\infty}{\longrightarrow} 0,
		\end{align*}
		
		\noindent because the sequence of real positive numbers $\left (\Vert z_n^{p(x)-2}\Vert_{L^{\frac{p(x)}{p(x)-2}}(\Omega)} \right )_{n\geq 1}$ is bounded. Let's explain this fact in details: From $z_n\to 2|V|$ in $L^{p(x)}(\Omega)$ we get that there is some $n_0\geq 1$ such that:
		
		\begin{equation}
			\Vert z_n\Vert_{L^{p(x)}(\Omega)}\leq 2\Vert V\Vert_{L^{p(x)}(\Omega)}+\Vert z_n-2|V|\Vert_{L^{p(x)}(\Omega)}\leq 2\Vert V\Vert_{L^{p(x)}(\Omega)}+1:=\alpha,\ \forall\ n\geq n_0.
		\end{equation}
		
		\noindent Thus $\left\Vert\dfrac{z_n}{\alpha}\right \Vert_{L^{p(x)}(\Omega)}\leq 1$ for any $n\geq n_0$. Using now Theorem \ref{athleb} \textbf{(25)} we deduce that $\rho_{p(x)}\left (\dfrac{z_n}{\alpha} \right )\leq \left\Vert\dfrac{z_n}{\alpha}\right \Vert_{L^{p(x)}(\Omega)}\leq 1$ for $n\geq n_0$. Because we have chosen $\alpha\geq 1$ we observe that:
		
		\begin{equation}
			\int_{\Omega} \dfrac{z_n^{p(x)}}{\alpha^{p^+}}\ dx\leq \int_{\Omega}\left (\dfrac{z_n}{\alpha}\right )^{p(x)}\ dx=\rho_{p(x)}\left (\dfrac{z_n}{\alpha} \right )\leq 1,\ \forall\ n\geq n_0.
		\end{equation}
		
		\noindent On the other hand 
		
		\begin{align*}
			\Vert z_n^{p(x)-2}\Vert_{L^{\frac{p(x)}{p(x)-2}}(\Omega)}&=\inf\left\{ \lambda>0\ |\ \rho_{\frac{p(x)}{p(x)-2}}\left(z_n^{p(x)-2}\right)\leq 1\right\}\\
			&=\inf\underbrace{\left\{ \lambda>0\ \bigg |\ \int_{\Omega}\dfrac{z_n^{p(x)}}{\lambda^{\frac{p(x)}{p(x)-2}}}\ dx\leq 1 \right\}}_{:=Z_n}
		\end{align*}
		
		\noindent If we choose $\lambda>1$ sufficiently large so that $\lambda^{\frac{p^+}{p^+-2}}\geq\alpha^{p^+}$, i.e. $\lambda\geq\alpha^{p^+-2}$ then:
		
		\begin{equation}
			\lambda^{\frac{p(x)}{p(x)-2}}=\lambda^{1+\frac{2}{p(x)-2}}\stackrel{\lambda>1}{\geq}\lambda^{1+\frac{2}{p^+-2}}=\lambda^{\frac{p^+}{p^+-2}}\geq\alpha^{p^+},\ \forall\ x\in\overline{\Omega}.
		\end{equation}
		
		\noindent This shows that $\displaystyle\int_{\Omega}\dfrac{z_n^{p(x)}}{\lambda^{\frac{p(x)}{p(x)-2}}}\ dx\leq\int_{\Omega} \dfrac{z_n^{p(x)}}{\alpha^{p^+}}\ dx\leq 1, \ \forall\ n\geq n_0$, which means in particular that $\lambda_1:=\alpha^{p^+-2}\in Z_n,\ \forall\ n\geq n_0$. Therefore:
		
		\begin{equation}
			\Vert z_n^{p(x)-2}\Vert_{L^{\frac{p(x)}{p(x)-2}}(\Omega)}=\inf Z_n\leq \lambda_1, \ \forall\ n\geq n_0,
		\end{equation}
		
		\noindent so that $\left (\Vert z_n^{p(x)-2}\Vert_{L^{\frac{p(x)}{p(x)-2}}(\Omega)} \right )_{n\geq 1}$ is bounded, as needed. This completes the proof of the fact that $\partial\mathcal{C}$ is a continuous operator. Using now Theorem \ref{athmgatfre} we conclude that $\mathcal{C}\in C^1\big (W^{1,p(x)}(\Omega) \big )$ and:
		
		\begin{equation}\label{3eqc3}
			\langle\mathcal{C}'(V),\phi\rangle=\int_{\Omega} |V(x)|^{p(x)-2}V(x)\phi(x)\ dx,\ \forall\ V,\phi\in W^{1,p(x)}(\Omega).
		\end{equation}

		\end{proof}

		\section{Introducing the auxiliary problem}
		
		\noindent For any $\lambda>0$ and any $g\in\mathcal{M}_{\lambda}\subseteq L^2(\Omega)$ consider the following problem:
		
		\begin{equation}\tag{$DE_{\lambda}$}\label{problemdelambda}
			\begin{cases}-\operatorname{div}\mathbf{a}\big (x,\nabla V(x)\big )+\lambda \overline{b}\big (x,V(x)\big )=g(x), & x\in\Omega\\[3mm] \mathbf{a}(x,\nabla V)\cdot\nu=0, & x\in \partial\Omega  \end{cases}
		\end{equation}
		
		\begin{definition}\label{auxdefinition}
			We say that $V\in W^{1,p(x)}(\Omega)$ is a weak solution of \eqref{problemdelambda} if for each test function $\phi\in W^{1,p(x)}(\Omega)$ we have that:
			
			\begin{equation}\label{auxeqdef}
				\int_{\Omega} \mathbf{a}\big (x,\nabla V(x)\big )\cdot\nabla \phi(x)\ dx+\lambda\int_{\Omega} \overline{b}\big (x,V(x)\big )\phi(x)\ dx=\int_{\Omega} g(x)\phi(x)\ dx.
			\end{equation}
		\end{definition}
		
		\begin{remark} Let $V$ be a weak solution of \eqref{problemdelambda}. From Proposition \ref{4prop1} \textnormal{\textbf{(6)}} $\mathbf{a}(\cdot,\nabla V)\cdot\nabla\phi\in L^1(\Omega)$. Since $V\in W^{1,p(x)}(\Omega)\hookrightarrow L^2(\Omega)$ we get from Proposition \ref{prop32} \textnormal{\textbf{(2)}} that $\overline{b}(\cdot,V)\in L^2(\Omega)$. Therefore, since $\phi\in W^{1,p(x)}(\Omega)\hookrightarrow L^2(\Omega)$, we deduce from \textit{Cauchy inequality} that $\overline{b}(\cdot,V)\phi\in L^1(\Omega)$. Similarly, since $g,\phi\in L^2(\Omega)\Rightarrow\ g\phi\in L^1(\Omega)$. So the integrals in \eqref{auxeqdef} are finite.
			
		\end{remark}
		
		\begin{proposition}\label{prop41} If $V\in W^{1,p(x)}(\Omega)$ is a weak solution of the problem \eqref{problemdelambda} then $V\in\mathcal{U}$.
		\end{proposition}
		
		\begin{proof}  For the first part of the proof note that $\forall\ \phi\in W^{1,p(x)}(\Omega)$:
			
			\begin{equation}
				\begin{cases}\displaystyle\int_{\Omega} \mathbf{a}(x,\nabla V(x))\cdot\nabla\phi(x)\ dx+\lambda\int_{\Omega}\overline{b}(x,V(x))\phi(x)\ dx=\int_{\Omega} g(x)\phi(x)\ dx\\[5mm]
					\displaystyle\int_{\Omega}\mathbf{a}(x,\nabla 0)\cdot\nabla\phi(x)\ dx+\lambda\int_{\Omega} \overline{b}(x,0)\cdot\phi(x)\ dx=\int_{\Omega} \lambda b(x,0)\phi(x)\ dx\end{cases}.
			\end{equation}
			
			\noindent Substracting these equations:
			
			\begin{equation}
				\int_{\Omega} \big (\mathbf{a}(x,\nabla V)-\mathbf{a}(x,\nabla 0)\big )\cdot\nabla\phi\ dx+\lambda\int_{\Omega}\big [\overline{b}(x,V)-\overline{b}(x,0)\big ]\phi\ dx=\int_{\Omega} \bigl (g-\lambda\overline{b}(x,0)\bigr)\phi\ dx.
			\end{equation}

			\noindent We choose $\phi=V^{-}=\begin{cases} 0, \ \text{a.e. on}\ \{x\in\Omega\ |\ V(x)\geq 0\}\\[3mm] -V(x),\ \text{a.e. on}\ \{x\in\Omega\ |\ V(x)<0\}:=\Omega_0\end{cases}\in W^{1,p(x)}(\Omega)$ as the test function. Notice that:
			
			\begin{equation}
				\nabla \phi=\nabla V^{-}=\begin{cases}0,\ \text{a.e. on}\ \{x\in\Omega\ |\ V(x)\geq 0\}\\[3mm]-\nabla V, \ \text{a.e. on}\ \{x\in\Omega\ |\ V(x)<0\}=\Omega_0 \end{cases}.
			\end{equation}
			
			\noindent We obtain that:
			
			\begin{equation}
				-\int_{\Omega_0} \underbrace{\mathbf{a}(x,\nabla V)\cdot\nabla V\ dx}_{\geq 0}-\lambda\int_{\Omega_0} \underbrace{\big [\overline{b}(x,V)-\overline{b}(x,0)\big ](V-0)}_{\geq 0}\ dx=\int_{\Omega} \underbrace{\bigl(g-\lambda\overline{b}(x,0)\bigr)}_{\geq 0} V^{-}\ dx.
			\end{equation}

			\noindent So the left hand side is negative while the right hand side is positive, since $g\in\mathcal{M}_{\lambda}$. Therefore $\displaystyle\int_{\Omega_0} \big [\overline{b}(x,V)-\overline{b}(x,0)\big ](V-0)\ dx=0$ which shows, using Proposition \ref{prop32} \textbf{(4)}, that $|\Omega_0|=0$, i.e. $V^{-}\equiv 0$ on $\Omega$ ($\lambda>0$). In conclusion $V=V^{+}-V^-=V^+\geq 0$ a.e. on $\Omega$.

%			First we show that $V\geq 0$ a.e. on $\Omega$. Take $\phi=V^{-}\in W^{1,p(x)}(\Omega)^{+}$ as the test function. Then $\nabla \phi=\nabla V^{-}=\begin{cases}0,\ \text{a.e. on}\ \{x\in\Omega\ |\ V(x)\geq 0\}\\[3mm]-\nabla V, \ \text{a.e. on}\ \{x\in\Omega\ |\ V(x)<0\}:=\Omega_0 \end{cases}$. Using Proposition \ref{prop32} \textbf{(4)} we obtain that:
%			
%			\begin{equation}0\geq -\int_{\Omega_0} \Phi(x,|\nabla V|)|\nabla V|\ dx-\lambda\int_{\Omega_0} \underbrace{\overline{b}(x,V)V}_{\geq 0}\ dx=\int_{\Omega} gV^{-}\ dx\geq 0,
%			\end{equation}
%
%		
%		\noindent As $\lambda>0$ we get that $\displaystyle\int_{\Omega_0} \overline{b}(x,V)V=0$, i.e. $V\equiv 0$ on $\Omega_0$ (false) or $|\Omega_0|=0$ (i.e. $V^-\equiv 0$). So $V=V^+\geq 0$ a.e. on $\Omega$.
			
			\noindent For the second part of the proof first note that $\forall\ \phi\in W^{1,p(x)}(\Omega)$:
			
			\begin{equation}
				\begin{cases}\displaystyle\int_{\Omega}\mathbf{a}(x,\nabla \delta_0)\cdot\nabla\phi(x)\ dx+\lambda\int_{\Omega} \overline{b}(x,\delta_0)\cdot\phi(x)\ dx=\int_{\Omega} \lambda b(x,\delta_0)\phi(x)\ dx\\[5mm]
					\displaystyle\int_{\Omega} \mathbf{a}(x,\nabla V(x))\cdot\nabla\phi(x)\ dx+\lambda\int_{\Omega}\overline{b}(x,V(x))\phi(x)\ dx=\int_{\Omega} g(x)\phi(x)\ dx\end{cases}.
			\end{equation}
			
			\noindent Substracting both equations leads to:
			
			\begin{equation}
				\int_{\Omega} \big (\mathbf{a}(x,\nabla \delta_0)-\mathbf{a}(x,\nabla V)\big )\cdot\nabla\phi\ dx+\lambda\int_{\Omega}\big [\overline{b}(x,\delta_0)-\overline{b}(x,V)\big ]\phi\ dx=\int_{\Omega} (\lambda\overline{b}(x,\delta_0)-g)\phi\ dx.
			\end{equation}

			\noindent We choose $\phi=(\delta_0-V)^{-}=\begin{cases} 0, \ \text{a.e. on}\ \{x\in\Omega\ |\ V(x)\leq \delta_0\}\\[3mm] -(\delta_0-V(x)),\ \text{a.e. on}\ \{x\in\Omega\ |\ V(x)>\delta_0\}:=\Omega_1\end{cases}\in W^{1,p(x)}(\Omega)$ as the test function. Notice that:
			
			\begin{equation}
				\nabla\phi=\nabla (\delta_0-V)^{-}=\begin{cases} 0,\ \text{a.e. on}\ \{x\in\Omega\ |\ V(x)\leq \delta_0\}\\[3mm] -(\nabla \delta_0-\nabla V)=\nabla V,\ \text{a.e. on}\ \Omega_1\end{cases}.
			\end{equation}
			
			\noindent We obtain that:
			
			\begin{equation}
				-\int_{\Omega_1} \underbrace{\mathbf{a}(x,\nabla V)\cdot\nabla V\ dx}_{\geq 0}-\lambda\int_{\Omega_1} \underbrace{\big [\overline{b}(x,\delta_0)-\overline{b}(x,V)\big ](\delta_0-V)}_{\geq 0}\ dx=\int_{\Omega} \underbrace{(\lambda\overline{b}(x,\delta_0)-g)}_{\geq 0} (\delta_0-V)^{-}\ dx.
			\end{equation}

			\noindent So the left hand side is negative while the right hand side is positive, since $g\in\mathcal{M}_{\lambda}$. Therefore $\displaystyle\int_{\Omega_1} \big [\overline{b}(x,\delta_0)-\overline{b}(x,V)\big ](\delta_0-V)\ dx=0$ which shows, using Proposition \ref{prop32} \textbf{(4)}, that $(\delta_0-V)^{-}\equiv 0$ on $\Omega$ ($\lambda>0$). In conclusion $\delta_0-V=(\delta_0-V_{\varepsilon})^{+}\geq 0$ a.e. on $\Omega$.
			
		\end{proof}
		
		\begin{remark}
			Proposition \ref{prop41} shows that if $V$ is a weak solution of \eqref{problemdelambda} then actually $V$ is a weak solution of the problem:
				
			\begin{equation}
					\begin{cases}-\operatorname{div}\mathbf{a}\big (x,\nabla V(x)\big )+\lambda b\big (x,V(x)\big )=g(x), & x\in\Omega\\[3mm] \mathbf{a}(x,\nabla V)\cdot\nu=0, & x\in \partial\Omega\\[3mm] 0\leq V(x)\leq \delta_0, & x\in\Omega \end{cases}.
				\end{equation}
				
		\end{remark}
		
			\begin{definition} The energy functional associated to \eqref{eqedglambda} is: $\mathcal{J}_{\lambda}:W^{1,p(x)}(\Omega)\to\mathbb{R}$ given by:
			
			\begin{equation}
				\mathcal{J}_{\lambda}(V)=\int_{\Omega} A\big (x,\nabla V(x)\big )\ dx+\lambda\int_{\Omega} \overline{B}\big (x,V(x)\big )\ dx-\int_{\Omega} g(x)V(x)\ dx.
			\end{equation}
			
		\end{definition}
		
		\begin{remark} From Proposition \ref{4prop1} \textnormal{\textbf{(6)}} and Proposition \ref{prop34} \textnormal{\textbf{(4)}}, we have that $\mathcal{J}_{\lambda}\in C^1\big (W^{1,p(x)}(\Omega)\big )$ and $\mathcal{J}_{\lambda}':W^{1,p(x)}(\Omega)\to W^{1,p(x)}(\Omega)^*$ is given by:
			
			\begin{equation}
				\langle\mathcal{J}_{\lambda}'(V),\phi\rangle=\int_{\Omega} \mathbf{a}\big (x,\nabla V\big )\cdot\nabla \phi\ dx+\lambda\int_{\Omega} \overline{b}\big (x,V\big )\phi\ dx-\int_{\Omega} g\phi\ dx,\ \forall\ V,\phi\in W^{1,p(x)}(\Omega).
			\end{equation}
			
			\noindent In particular $V\in W^{1,p(x)}(\Omega)$ is a weak solution of \eqref{eqedg} iff $\mathcal{J}'_{\lambda}(V)\equiv 0$.
			
		\end{remark}

		\section{The perturbed problem}
		
		\noindent We associate to \eqref{eqedg} the following perturbed problem:
		
		\begin{equation}\tag{$DE_{\lambda,\epsilon}$}\label{eqedgeps}
			\begin{cases}-\operatorname{div}\mathbf{a}\big (x,\nabla V(x)\big )+\lambda \overline{b}\big (x,V(x)\big )+\epsilon |V(x)|^{p(x)-2}V(x)=g(x), & x\in\Omega\\[3mm] \mathbf{a}(x,\nabla V)\cdot\nu=0, & x\in\partial\Omega\end{cases}
		\end{equation}
		
			\begin{definition} 
			We say that $V\in W^{1,p(x)}(\Omega)$ is a weak solution of \eqref{eqedgeps} if for any test function $\phi\in W^{1,p(x)}(\Omega)$ we have that:
			
			\begin{equation}
				\int_{\Omega} \mathbf{a}\big (x,\nabla V(x)\big )\cdot\nabla \phi\ dx+\lambda\int_{\Omega} \overline{b}\big (x,V(x)\big )\phi\ dx+\epsilon\int_{\Omega} |V(x)|^{p(x)-2}V(x)\phi\ dx=\int_{\Omega} g(x)\phi\ dx.
			\end{equation}
		\end{definition}
		
		\begin{proposition} If $V_{\epsilon}\in W^{1,p(x)}(\Omega)$ is a weak solution of the problem \eqref{eqedgeps} then $V_{\epsilon}\in\mathcal{U}$.
		\end{proposition}
		
		\begin{proof} First we show that $V_{\epsilon}\geq 0$ a.e. on $\Omega$. Take $\phi=V_{\epsilon}^{-}\in W^{1,p(x)}(\Omega)^{+}$ as the test function. Then $\nabla \phi=\nabla V_{\epsilon}^{-}=\begin{cases}0,\ \text{a.e. on}\ \{x\in\Omega\ |\ V_{\epsilon}(x)\geq 0\}\\[3mm]-\nabla V_{\epsilon}, \ \text{a.e. on}\ \{x\in\Omega\ |\ V_{\epsilon}(x)<0\}:=\Omega_0 \end{cases}$. We obtain that:
			
			\begin{align*}0&\geq -\int_{\Omega_0} \underbrace{\Phi(x,|\nabla V_{\epsilon}|)|\nabla V_{\epsilon}|}_{\geq 0}\ dx-\lambda\int_{\Omega_0} \underbrace{\big[\overline{b}(x,V_{\epsilon})-\overline{b}(x,0) \big ](V_{\epsilon}-0)}_{\geq 0} dx-\epsilon\int_{\Omega_0} |V_{\epsilon}|^{p(x)}\ dx\\
				&=\int_{\Omega} \underbrace{\big (g-\lambda\overline{b}(x,0)\big )}_{\geq 0}V_{\epsilon}^{-}\ dx\geq 0.
			\end{align*}
			
			\noindent As $\lambda>0$ we get that $\displaystyle\int_{\Omega_0} \big[\overline{b}(x,V_{\epsilon})-\overline{b}(x,0) \big ](V_{\epsilon}-0) dx=0$, i.e. $V_{\epsilon}^{-}\equiv 0$. So $V_{\epsilon}=V_{\epsilon}^+\geq 0$ a.e. on $\Omega$.
			
			\noindent For the second part of the proof first note that $\forall\ \phi\in W^{1,p(x)}(\Omega)$:
			
			\begin{equation}
				\begin{cases}\displaystyle\int_{\Omega}\mathbf{a}(x,\nabla \delta_0)\cdot\nabla\phi\ dx+\lambda\int_{\Omega} \overline{b}(x,\delta_0)\phi\ dx+\epsilon\int_{\Omega} \delta_0^{p(x)-1}\phi\ dx=\int_{\Omega} \big (\lambda b(x,\delta_0)+\epsilon\delta_0^{p(x)-1}\big )\phi\ dx\\[5mm]
					\displaystyle\int_{\Omega} \mathbf{a}(x,\nabla V_{\epsilon}(x))\cdot\nabla\phi\ dx+\lambda\int_{\Omega} \overline{b}(x,V_{\epsilon}(x))\phi\ dx+\epsilon\int_{\Omega} |V_{\epsilon}(x)|^{p(x)-2}V_{\epsilon}(x)\phi\ dx=\int_{\Omega} g(x)V_{\epsilon}(x)\ dx\end{cases}.
			\end{equation}
			
			\noindent Substracting both equations, and using the fact that $V_{\epsilon}\geq 0$ a.e. on $\Omega$ we have that: $|V_{\epsilon}|^{p(x)-2}V_{\epsilon}=V_{\epsilon}^{p(x)-1}$, and this leads to:
			
			\begin{align*}
				\int_{\Omega} \big (\mathbf{a}(x,\nabla \delta_0)-\mathbf{a}(x,\nabla V_{\epsilon})\big )\cdot\nabla\phi\ dx&+\lambda\int_{\Omega}\big [\overline{b}(x,\delta_0)-\overline{b}(x,V_{\epsilon}))]\phi\ dx+\epsilon\int_{\Omega} \big (\delta_0^{p(x)-1}-V_{\epsilon}^{p(x)-1}\big )\phi\ dx\\
				=\int_{\Omega} \big (\lambda b(x,\delta_0)+\epsilon\delta_0^{p(x)-1}-g\big )\phi\ dx.
			\end{align*}

			\noindent We choose $\phi=(\delta_0-V_{\epsilon})^{-}=\begin{cases} 0, \ \text{a.e. on}\ \{x\in\Omega\ |\ V_{\epsilon}(x)\leq \delta_0\}\\[3mm] -(\delta_0-V_{\epsilon}(x)),\ \text{a.e. on}\ \{x\in\Omega\ |\ V_{\epsilon}(x)>\delta_0\}:=\Omega_1\end{cases}\in W^{1,p(x)}(\Omega)$ as the test function. Notice that:
			
			\begin{equation}
				\nabla\phi=\nabla (\delta_0-V_{\epsilon})^{-}=\begin{cases} 0,\ \text{a.e. on}\ \{x\in\Omega\ |\ V_{\epsilon}(x)\leq \delta_0\}\\[3mm] -(\nabla \delta_0-\nabla V_{\epsilon})=\nabla V_{\epsilon},\ \text{a.e. on}\ \Omega_1\end{cases}.
			\end{equation}
			
			\noindent We obtain that:

				\begin{align*}
				0\geq&-\int_{\Omega} \big (\mathbf{a}(x,\nabla \delta_0)-\mathbf{a}(x,\nabla V_{\epsilon})\big )\cdot\big (\nabla\delta_0-\nabla V_{\epsilon}\big )\ dx-\lambda\int_{\Omega}\big [\overline{b}(x,\delta_0)-\overline{b}(x,V_{\epsilon}))](\delta_0-V_{\epsilon})\ dx\\
				&-\epsilon\int_{\Omega} \big (\delta_0^{p(x)-1}-V_{\epsilon}^{p(x)-1}\big )(\delta_0-V_{\epsilon})\ dx=\int_{\Omega} \big (\lambda b(x,\delta_0)+\epsilon\delta_0^{p(x)-1}-g\big )(\delta_0-V_{\epsilon})^{-}\geq 0\ dx.
			\end{align*}

			\noindent It is straightforward to check that for each $x\in\overline{\Omega}$ the function $\theta\mapsto \theta^{p(x)-1}$ is a strictly increasing function on $[0,\infty)$, as $p(x)\geq p^->1$ for all $x\in\overline{\Omega}$. So the left hand side is negative while the right hand side is positive as $g\in\mathcal{M}_{\lambda}$. Therefore $\lambda\displaystyle\int_{\Omega}\big [\overline{b}(x,\delta_0)-\overline{b}(x,V_{\epsilon}))](\delta_0-V_{\epsilon})\ dx=0$ which shows that $(\delta_0-V_{\epsilon})^{-}\equiv 0$ on $\Omega$ ($\lambda>0$). In conclusion $\delta_0-V_{\epsilon}=(\delta_0-V_{\epsilon})^{+}\geq 0$ a.e. on $\Omega$.
			
		\end{proof}
		
		\begin{definition} The energy functional associated to \eqref{eqedgeps} is: $\mathcal{J}_{\lambda,\epsilon}:W^{1,p(x)}(\Omega)\to\mathbb{R}$ given by:
			
			\begin{equation}
				\mathcal{J}_{\lambda,\epsilon}(V)=\int_{\Omega} A\big (x,\nabla V(x)\big )\ dx+\lambda\int_{\Omega} \overline{B}\big (x,V(x)\big )\ dx+\epsilon\int_{\Omega}\dfrac{|V(x)|^{p(x)}}{p(x)}\ dx-\int_{\Omega} g(x)V(x)\ dx.
			\end{equation}
			
		\end{definition}
		
		\begin{remark} From Proposition \ref{4prop1} \textnormal{\textbf{(6)}}, Proposition \ref{prop35} \textnormal{\textbf{(4)}} and Proposition \ref{prop39} we have that $\mathcal{J}_{\lambda,\epsilon}\in C^1\big (W^{1,p(x)}(\Omega)\big )$ and its differential $\mathcal{J}_{\lambda,\epsilon}:W^{1,p(x)}(\Omega)\to W^{1,p(x)}(\Omega)^*$ is given by:
			
			\begin{equation}
				\langle\mathcal{J}_{\lambda,\epsilon}'(V),\phi\rangle=\int_{\Omega} \mathbf{a}\big (x,\nabla V\big )\cdot\nabla \phi\ dx+\lambda\int_{\Omega} \overline{b}\big (x,V\big )\phi\ dx+\epsilon\int_{\Omega} |V|^{p(x)-2}V\phi\ dx-\int_{\Omega} g\phi\ dx,
			\end{equation}
			
			\noindent for every $V,\phi\in W^{1,p(x)}(\Omega)$. In particular $V\in W^{1,p(x)}(\Omega)$ is a weak solution of \eqref{eqedgeps} iff $\mathcal{J}'_{\lambda,\epsilon}(V)\equiv 0$, i.e. $V$ is a critical point of $\mathcal{J}_{\lambda,\epsilon}$.
			
		\end{remark}
		
		\begin{remark} For $\epsilon=0$ we shall write $\mathcal{J}_{\lambda,0}=\mathcal{J}_{\lambda}$.	
		\end{remark}
		
		\begin{theorem}\label{thmepsilon}
			For each $\epsilon>0$, problem \eqref{eqedgeps} admits a unique weak solution that will be denoted by $V_{\epsilon}\in\mathcal{U}$.
		\end{theorem}
		
		\begin{proof} \textbf{(Uniqueness)} Suppose that $V_1,V_2\in W^{1,p(x)}(\Omega)$ are two weak solutions of \eqref{eqedgeps}. Choosing $\phi=V_1-V_2\in W^{1,p(x)}(\Omega)$ as the test function we easily get that:
			
			\begin{align*}
				&\int_{\Omega} \big (\mathbf{a}(x,\nabla V_1)-\mathbf{a}(x,\nabla V_2) \big)\cdot (\nabla V_1-\nabla V_2)\ dx+\lambda\int_{\Omega} \big [\overline{b}(x,V_1)-\overline{b}(x,V_2)\big ](V_1-V_2)\ dx+\\
				&+\epsilon\int_{\Omega} \big (|V_1|^{p(x)-2}V_1-|V_2|^{p(x)-2}V_2\big )(V_1-V_2)\ dx=0.
			\end{align*}
			
			\noindent Using the fact that $\mathbf{a}(x,\cdot)$ is monotone (see Proposition \ref{4prop1}), $\lambda>0$ and the function $\mathbb{R}\ni\theta\mapsto |\theta|^{p(x)-2}\theta=|\theta|^{p(x)-1}\operatorname{sgn}(\theta)$ is strictly increasing on $\mathbb{R}$ for all $x\in\overline{\Omega}$ as $p(x)\geq p^{-}>1,\ \forall x\in\overline{\Omega}$, we obtain that every integral in the left hand side is positive. We infer that $\displaystyle\int_{\Omega} \big [\overline{b}(x,V_1)-\overline{b}(x,V_2)\big ](V_1-V_2)\ dx=0$. So, from Proposition \ref{prop32} \textbf{(4)}, $V_1\equiv V_2$, as needed.
			
			\bigskip
			
			\noindent\textbf{(Existence)} We will show the existence in 5 steps:
			
			\bigskip

			\noindent $\blacktriangleright$ \textbf{Fact I: $\mathcal{J}_{\lambda,\epsilon}$ is well-defined and bounded from below.}
			
			\noindent Using the minimum of a quadratic function and the fact that $\lambda b(x,0)\leq g(x)\leq\lambda b(x,\delta_0)$ for a.e. $x\in \Omega$ we find that for every $s\in\mathbb{R}$ and a.e. $x\in\Omega$:
			
			\begin{align*}
				\lambda\overline{B}(x,s)-g(x)s&=\begin{cases}\lambda b(x,0)s+\dfrac{\lambda s^2}{2}-g(x)s, & s<0 \\[3mm] \displaystyle\int_{0}^s\lambda b(x,\tau)\ d\tau-g(x)s=\displaystyle\int_{0}^s\lambda b(x,\tau)-g(x)\ d\tau, & s\in [0,\delta_0]\\[3mm]\lambda\displaystyle\int_{0}^{\delta_0}b(x,\tau)\ d\tau+\lambda b(x,\delta_0)(s-\delta_0)-g(x)s+\dfrac{\lambda(s-\delta_0)^2}{2}, & s>\delta_0\end{cases}\\[3mm]
				&\geq \begin{cases}-\dfrac{(\lambda b(x,0)-g(x))^2}{2\lambda}, & s<0\\[3mm] \displaystyle\int_{0}^s\lambda b(x,0)-\lambda b(x,\delta_0)\ d\tau, & s\in [0,\delta_0]\\[3mm] \lambda\displaystyle\int_{0}^{\delta_0}b(x,\tau)\ d\tau-\delta_0 g(x)-\dfrac{(\lambda b(x,\delta_0)-g(x))^2}{2\lambda}, & s>\delta_0 \end{cases}\\[3mm]
				&\geq \begin{cases}-\dfrac{\lambda}{2} (b(x,\delta_0)-b(x,0))^2, & s<0\\[3mm] -\lambda s\big (b(x,\delta_0)-b(x,0)), & s\in [0,\delta_0]\\[3mm] \lambda\delta_0 b(x,0)-\delta_0 g(x)-\dfrac{\lambda}{2} (b(x,\delta_0)-b(x,0))^2, & s>\delta_0 \end{cases}\\[3mm]
				&\geq \begin{cases}-\dfrac{\lambda}{2} (b(x,\delta_0)-b(x,0))^2, & s<0\\[3mm] -\lambda \delta_0\big (b(x,\delta_0)-b(x,0)), & s\in [0,\delta_0]\\[3mm] -\lambda \delta_0\big (b(x,\delta_0)-b(x,0))-\dfrac{\lambda}{2} (b(x,\delta_0)-b(x,0))^2, & s>\delta_0 \end{cases}\\[3mm]
				&\geq -\lambda \delta_0\big (b(x,\delta_0)-b(x,0))-\dfrac{\lambda}{2} (b(x,\delta_0)-b(x,0))^2.
			\end{align*}
			
			\noindent Therefore, the following inequality holds:
			
			\begin{align}\label{eqlowerbound}
				\nonumber \mathcal{J}_{\lambda,\epsilon}(V)&=\int_{\Omega} \underbrace{A\big (x,\nabla V(x)\big )}_{\geq 0}\ dx+\epsilon\int_{\Omega}\underbrace{\dfrac{|V(x)|^{p(x)}}{p(x)}}_{\geq 0}\ dx+\lambda\int_{\Omega} \overline{B}\big (x,V(x)\big )\ dx-\int_{\Omega} g(x)V(x)\ dx\\ \nonumber
				&\geq \int_{\Omega}\lambda\overline{B}\big (x,V(x)\big )-g(x)V(x)\ dx\\ 
				&\geq -\lambda\delta_0\int_{\Omega} b(x,\delta_0)-b(x,0)\ dx-\dfrac{\lambda}{2}\int_{\Omega} \big (b(x,\delta_0)-b(x,0)\big )^2\ dx>-\infty.
			\end{align}

			\noindent$\blacktriangleright$ \textbf{Fact II:} $\mathcal{J}_{\lambda,\epsilon}\in C^1\big (W^{1,p(x)}(\Omega)\big )$ 
			
			\noindent $\blacktriangleright$ \textbf{Fact III: $\mathcal{J}_{\lambda,\epsilon}$ is a strictly convex functional.}
			
			\noindent Here use the fact that for a.e. $x\in\Omega$ the following functions are strictly convex: $A(x,\cdot)$ (see Proposition \ref{4prop1} \textbf{(4)}), $\overline{B}(x,\cdot)$ (see Proposition \ref{prop33} \textbf{(5)}) and $[0,\infty)\ni s\mapsto \dfrac{s^{p(x)}}{p(x)}$ because its first derivative $[0,\infty)\ni s\mapsto s^{p(x)-1}$ is a strictly increasing function as $p(x)\geq p^{-}>1$ for all $x\in\overline{\Omega}$. This shows that $\mathcal{J}_{\lambda,\epsilon}$ is strictly convex.

			\bigskip
			
			\noindent $\blacktriangleright$ \textbf{Fact IV: $\mathcal{J}_{\lambda,\epsilon}$ is a weakly lower semicontinuous functional.}
			
			\noindent Let now some $V\in W^{1,p(x)}(\Omega)$ and a sequence $(V_n)_{n\geq 1}\subset W^{1,p(x)}(\Omega)$ such that $V_n\rightharpoonup V$ in $W^{1,p(x)}(\Omega)$. Since $\mathcal{J}_{\lambda,\epsilon}\in C^1\big (W^{1,p(x)}(\Omega)\big)$ we have that $\mathcal{J}'_{\lambda,\epsilon}(V)\in W^{1,p(x)}(\Omega)^*$. Therefore: $\lim\limits_{n\to\infty} \langle \mathcal{J}'_{\lambda,\epsilon}(V),V_n-V\rangle=0$.
			
			\noindent Using now the convexity of $\mathcal{J}_{\lambda,\epsilon}$ and the fact that $\mathcal{J}_{\lambda,\epsilon}\in C^1\big (W^{1,p(x)}(\Omega)\big)$ we have from \cite[Proposition 42.6, page 247]{zeidler} that for each $n\geq 1$:
			
			\begin{equation*}
				\mathcal{J}_{\lambda,\epsilon}(V_n)\geq \mathcal{J}_{\lambda,\epsilon}(V)+\langle \mathcal{J}'_{\lambda,\epsilon}(V),V_n-V\rangle.
			\end{equation*}
			
			\noindent So, $\liminf\limits_{n\to\infty} \mathcal{J}_{\lambda,\epsilon}(V_n)\geq \mathcal{J}_{\lambda,\epsilon}(V)+\liminf\limits_{n\to\infty}\langle\mathcal{J}'_{\lambda,\epsilon}(V),V_n-V\rangle=\mathcal{J}_{\lambda,\epsilon}(V)$. This shows that $\mathcal{J}_{\lambda,\epsilon}$ is weakly lower semicontinuous.
			
			\bigskip
			
			\noindent $\blacktriangleright$ \textbf{Fact V: $\mathcal{J}_{\lambda,\epsilon}$ is a coercive functional for any $\epsilon>0$.}
			
			\noindent Let $(V_n)_{n\geq 1}\subset W^{1,p(x)}(\Omega)$ with $\Vert V_n\Vert_{W^{1,p(x)}(\Omega)}\longrightarrow\infty$. From Theorem \ref{athsob} \textbf{(4)} we have that $\varrho_{p(x)}(V_n)=\rho_{p(x)}(V_n)+\rho_{p(x)}(|\nabla V_n|)\longrightarrow\infty$. This means that $\rho_{p(x)}(V_n)\to\infty$ or $\rho_{p(x)}(|\nabla V_n|)=\displaystyle\int_{\Omega}|\nabla V_n|^{p(x)}\ dx\to\infty$. In the second case, from \textbf{(H7)}, we get that $\displaystyle\int_{\Omega} A(x,\nabla V_n)\ dx\to\infty$. Therefore:
			
			\begin{align}
				\mathcal{J}_{\lambda,\epsilon}(V_n)&=\int_{\Omega} A(x,\nabla V_n)\ dx+\lambda\int_{\Omega} \overline{B}(x,V_n)\ dx+\epsilon\int_{\Omega}\dfrac{|V_n|^{p(x)}}{p(x)}\ dx-\int_{\Omega} gV_n\ dx\nonumber\\
				&\geq \int_{\Omega} A(x,\nabla V_n)\ dx+\dfrac{\epsilon}{p^+}\int_{\Omega} |V_n|^{p(x)}\ dx+\int_{\Omega} \lambda \overline{B}(x,V_n)-gV_n\ dx\nonumber\\
			\eqref{eqlowerbound}\ \ \ \ \	&\geq \int_{\Omega} A(x,\nabla V_n)\ dx+\dfrac{\epsilon}{p^+}\rho_{p(x)}(V_n)-\lambda\delta_0\int_{\Omega} b(x,\delta_0)-b(x,0)\ dx-\dfrac{\lambda}{2}\int_{\Omega} \big (b(x,\delta_0)-b(x,0)\big )^2\ dx\\
				&\stackrel{n\to\infty}{\longrightarrow}\infty\nonumber
			\end{align}

			\noindent Combining the fact that $W^{1,p(x)}(\Omega)$ is a reflexive Banach space, \textbf{Fact III} and \textbf{Fact IV}, we deduce based on Theorem \ref{athmstru}, that there is some $V_{\epsilon}\in W^{1,p(x)}(\Omega)$ such that: $\mathcal{J}_{\lambda,\epsilon}(V_{\epsilon})=\displaystyle\inf_{U\in W^{1,p(x)}(\Omega)} \mathcal{J}_{\lambda,\epsilon}(U)$.
			
			\noindent Since $\mathcal{J}_{\lambda,\epsilon}\in C^1\big (W^{1,p(x)}(\Omega) \big )$ and $V_{\epsilon}$ is a (global) minimum of $\mathcal{J}_{\lambda,\epsilon}$, we deduce that $\mathcal{J}_{\lambda,\epsilon}'(V_{\epsilon})=0$, i.e. $V_{\epsilon}$ is a weak solution of the problem \eqref{eqedgeps}.\footnote{See Corollary 2.5 at page 53 from \cite{coleman}.}
		\end{proof}

		\section{Existence and uniqueness for the auxiliary problem}
		
		\begin{lemma}\label{3lemgama} Let any $V\in W^{1,p(x)}(\Omega)$ and any $(V_n)_{n\geq 1}\subset W^{1,p(x)}(\Omega)$ such that $V_n\longrightarrow V$ in $W^{1,p(x)}(\Omega)$. Then for every sequence of real numbers $(\epsilon_n)_{n\geq 1}$ that converges to $0$ we have that:
			
			\begin{equation}
				\lim\limits_{n\to\infty} \mathcal{J}_{\lambda,\epsilon_n}(V_n)=\mathcal{J}_{\lambda}(V).
			\end{equation}
			
		\end{lemma}
		
		\begin{proof} We have for each $n\geq 1$ that:
			
			\begin{align*} \big |\mathcal{J}_{\lambda,\epsilon_n}(V_n)-\mathcal{J}_{\lambda}(V)\big |&=\left |\mathcal{J}_{\lambda}(V_n)-\mathcal{J}_{\lambda}(V)+\epsilon_n\int_{\Omega}\dfrac{|V_n|^{p(x)}}{p(x)}\ dx\right |\\
				&\leq \big |\mathcal{J}_{\lambda}(V_n)-\mathcal{J}_{\lambda}(V) \big |+|\epsilon_n|\int_{\Omega}\dfrac{|V_n|^{p(x)}}{p(x)}\ dx\\
				&\leq \big |\mathcal{J}_{\lambda}(V_n)-\mathcal{J}_{\lambda}(V) \big |+\dfrac{|\epsilon_n|}{p^{-}}\int_{\Omega}|V_n|^{p(x)}\ dx\\
				&=\big |\mathcal{J}_{\lambda}(V_n)-\mathcal{J}_{\lambda}(V) \big |+\dfrac{|\epsilon_n|}{p^{-}}\rho_{p(x)}(V_n)\\
				&\longrightarrow 0\ \text{as}\ n\to\infty.
			\end{align*}
			
			\noindent Above we have used the continuity of $\mathcal{J}_{\lambda}$ on $W^{1,p(x)}(\Omega)$ and the fact that $\big (\rho_{p(x)}(V_n)\big )_{n\geq 1}$ is a bounded sequence. Let's explain this fact in more details: from Theorem \ref{athleb} \textbf{(11)} we have that:
			
			\begin{equation}
				0\leq \rho_{p(x)}(V_n)=\rho_{p(x)}(V_n-V+V)\leq 2^{p^+-1}\rho_{p(x)}(V_n-V)+2^{p^+-1}\rho_{p(x)}(V).
			\end{equation}
			
			\noindent Since $V_n\to V$ in $W^{1,p(x)}(\Omega)\hookrightarrow L^{p(x)}(\Omega)$ we also have that $V_n\to V$ in $L^{p(x)}(\Omega)$. Therefore, based on Theorem \ref{athleb} \textbf{(7)} we have that $\lim\limits_{n\to\infty}\rho_{p(x)}(V_n-V)=0$. Hence $\big (\rho_{p(x)}(V_n-V)\big )_{n\geq 1}$ is a bounded sequence (being convergent). This shows that $\big (\rho_{p(x)}(V_n)\big )_{n\geq 1}$ is also bounded and the proof is complete.

			\begin{remark} In particular, Lemma \ref{3lemgama} shows that the family of functionals $\big (\mathcal{J}_{\lambda,\epsilon} \big )_{\epsilon>0}$ is $\Gamma$-convergent to $\mathcal{J}_{\lambda}$.
			\end{remark}

		\end{proof}
		
		\begin{theorem}\label{3thmrot} Problem \eqref{eqedglambda} admits a unique weak solution for any $g\in\mathcal{M}_{\lambda}$ that will be denoted by $V\in\mathcal{U}$.
		\end{theorem}
		
		\begin{proof} The proof here is a particular case of the proof of the Theorem \ref{thmepsilon}, because $\mathcal{J}_{\lambda}=\mathcal{J}_{\lambda,0}$. The only step in that proof when it was used that $\epsilon>0$ is in \textbf{Fact V} (coercivity). So the uniqueness and the first 4 facts from the existence part are proved.
			
		\noindent We want to show that all the requirements of the Theorem \ref{athmstru} are satisfied for the functional $\mathcal{J}_{\lambda}:X:=W^{1,p(x)}(\Omega)\cap\mathcal{U}\to\mathbb{R}$.
		
		\begin{enumerate}
			\item[$\bullet$] $W^{1,p(x)}(\Omega)$ is a \textbf{reflexive Banach space}\footnote{See Remark \ref{aremsob} from the Appendix.} and $X$ is a \textbf{weakly closed} subset of $W^{1,p(x)}(\Omega)$. 
			
			\noindent It is easy to see that $X$ is a convex subset of $W^{1,p(x)}(\Omega)$. We now show that it is strongly closed. Indeed, let $V\in W^{1,p(x)}(\Omega)$ and $(V_n)_{n\geq 1}\subset X$ such that $V_n\longrightarrow V$ in $W^{1,p(x)}(\Omega)$. It follows that $V_n\longrightarrow V$ in $L^{p(x)}(\Omega)\hookrightarrow L^{p^-}(\Omega)$. Therefore there is a subsequence $[0,\delta_0]\ni V_{n_k}(x)\stackrel{k\to\infty}{\longrightarrow} V(x)$ for a.e. $x\in \Omega$.\footnote{See the Corollary from \cite[Page 234]{Jones}.} So $V(x)\in [0,\delta_0]$ for a.e. $x\in\Omega$ which means that $V\in W^{1,p(x)}(\Omega)\cap\mathcal{U}=X$.
			
			\noindent Using now \textit{Mazur's theorem}\footnote{See Theorem 3.3.18 and Corollary 3.3.22 in \cite{papa1}.} we conclude that $X$ is weakly closed.
			
			%		\noindent Let $V\in W^{1,p(x)}(\Omega)$ and $(V_n)_{n\geq 1}\subset M$ such that $V_n\rightharpoonup V$ in $W^{1,p(x)}(\Omega)$. First we infer that $(V_n)_{n\geq 1}$ is a bounded sequence in $W^{1,p(x)}(\Omega)$. Since the embedding $W^{1,p(x)}\stackrel{\text{c}}{\hookrightarrow} L^{p(x)}(\Omega)$ is compact, we deduce that the bounded sequence $(V_{n})_{n\geq 1}$ has a (strongly) convergent subsequence in $L^{p(x)}(\Omega)$.  
			
			\item[$\bullet$] $\mathcal{J}_{\lambda}=\mathcal{J}_{\lambda,0}$ is \textbf{weakly lower semicontinuous}.
			
			\item[$\bullet$] $\mathcal{J}_{\lambda}:X\to\mathbb{R}$ is \textbf{coercive}, meaning that for $V\in X$: $\lim\limits_{\Vert V\Vert_{W^{1,p(x)}(\Omega)}\to \infty} \mathcal{J}_{\lambda}(V)=\infty$.
			
			\noindent Let $(V_n)_{n\geq 1}\subset X=W^{1,p(x)}(\Omega)\cap\mathcal{U}$ with $\Vert V_n\Vert_{W^{1,p(x)}(\Omega)}\longrightarrow\infty$. From Theorem \ref{athsob} \textbf{(4)} we have that $\varrho_{p(x)}(V_n)=\rho_{p(x)}(V_n)+\rho_{p(x)}(|\nabla V_n|)\longrightarrow\infty$. The key step here is to observe that for each $n\geq 1$:
			
			\begin{align}
				\rho_{p(x)}(V_n)=\int_{\Omega} |V_n(x)|^{p(x)}\ dx\leq \int_{\Omega} \delta_0^{p(x)}\ dx\leq\max\{\delta_0^{p^+},\delta_0^{p^-}\}|\Omega|\\ 
				\Longrightarrow\ \rho_{p(x)}(|\nabla V_n|)=\int_{\Omega} |\nabla V_n|^{p(x)}\ dx\longrightarrow\infty.
			\end{align}
			
			\noindent Therefore, from \textbf{(H7)}:

			\begin{align}
				\mathcal{J}_{\lambda}(V_n)&=\int_{\Omega} A(x,\nabla V_n)\ dx+\lambda\int_{\Omega} \overline{B}(x,V_n)\ dx-\int_{\Omega} gV_n\ dx\nonumber\\
				&\geq\int_{\Omega} A(x,\nabla V_n)\ dx+\int_{\Omega} \lambda\overline{B}(x,V_n)-gV_n\ dx\nonumber\\
			\eqref{eqlowerbound}\ \ \ \ \ 	&\geq\int_{\Omega} A(x,\nabla V_n)\ dx-\lambda\delta_0\int_{\Omega} b(x,\delta_0)-b(x,0)\ dx-\dfrac{\lambda}{2}\int_{\Omega} \big (b(x,\delta_0)-b(x,0)\big )^2\ dx\stackrel{n\to\infty}{\longrightarrow}\infty.
			\end{align}

		\end{enumerate}
		
		\noindent Having now all the requirements fulfilled we deduce that there is some $V\in X= W^{1,p(x)}(\Omega)\cap\mathcal{U}$ such that $\mathcal{J}_{\lambda}(V)=\displaystyle\min_{U\in X} \mathcal{J}_{\lambda}(U)$. In what follows, we prove a very important lemma:

			\begin{lemma}\label{lemlocglob} We have that: $m:=\displaystyle\inf_{U\in X} \mathcal{J}_{\lambda}(U)=\mathcal{J}_{\lambda}(V)=\inf_{U\in W^{1,p(x)}(\Omega)} \mathcal{J}_{\lambda}(U):=d$.
		\end{lemma}
		
		\begin{proof}

			\noindent Note that, since $\mathcal{J}_{\lambda}=\mathcal{J}_{\lambda,0}$ is bounded from below, we have that $m,d\in\mathbb{R}$. It is obvious that since $X\subset W^{1,p(x)}(\Omega)$ we have that $m=\displaystyle\inf_{U\in X} \mathcal{J}_{\lambda}(U)\geq\inf_{U\in W^{1,p(x)}(\Omega)} \mathcal{J}_{\lambda}(U)=d$. We need to show that the reverse inequality is also true: $d\geq m$.

			\noindent For each $\epsilon>0$ we take $V_{\epsilon}\in X=W^{1,p(x)}(\Omega)\cap\mathcal{U}$ the unique weak solution of the problem \eqref{eqedgeps}, i.e. the unique (global) minimum of the functional $\mathcal{J}_{\lambda,\epsilon}$. First, since $V_{\epsilon}\in X$ we have that:
			
			\begin{align*}
				\min_{U\in W^{1,p(x)}(\Omega)} \mathcal{J}_{\lambda,\epsilon}(U)&=\mathcal{J}_{\lambda,\epsilon}(V_\epsilon)=\mathcal{J}_{\lambda}(V_{\epsilon})+\epsilon\int_{\Omega}\dfrac{|V_{\epsilon}|^{p(x)}}{p(x)}\ dx\geq \mathcal{J}_{\lambda}(V_{\epsilon})+0\\
				(V_{\epsilon}\in X)\ \ \ \ \ \ \ 		&\geq \mathcal{J}_{\lambda}(V)=m.
			\end{align*}

			\noindent Secondly we can remark that:
			
			\begin{align*} \mathcal{J}_{\lambda,\epsilon}(V_{\epsilon})=\min_{U\in W^{1,p(x)}(\Omega)} \mathcal{J}_{\lambda,\epsilon}(U)&\leq \mathcal{J}_{\lambda,\epsilon}(V)=\mathcal{J}_{\lambda}(V)+\epsilon\int_{\Omega}\dfrac{|V|^{p(x)}}{p(x)}\ dx\\
				(V\in \mathcal{U})\ \ \ \ \ \ \ \	&\leq \mathcal{J}_{\lambda}(V)+\dfrac{\epsilon}{p^{-}}\int_{\Omega} 1^{p(x)}\ dx\\
				&=m+\dfrac{\epsilon |\Omega|}{p^{-}}
			\end{align*}
			
			\noindent Thence:
			
			\begin{equation}
				m\leq \min_{U\in W^{1,p(x)}(\Omega)} \mathcal{J}_{\lambda,\epsilon}(U)\leq m+\dfrac{\epsilon |\Omega|}{p^{-}},\ \forall\ \epsilon>0.
			\end{equation}
			
			\noindent This shows that:
			
			\begin{equation}\label{3eqmd}
				\lim\limits_{\epsilon\to 0^+} \min_{U\in W^{1,p(x)}(\Omega)} \mathcal{J}_{\lambda,\epsilon}(U)=m=\min_{U\in X} \mathcal{J}_{\lambda}(U),\ \text{or}\ \lim\limits_{\epsilon\to 0^+}\mathcal{J}_{\lambda,\epsilon}(V_{\epsilon})=\mathcal{J}_{\lambda}(V).
			\end{equation}
			
			\noindent Since $d=\displaystyle\inf_{U\in W^{1,p(x)}(\Omega)} \mathcal{J}_{\lambda}(U)$ we can take a minimizing sequence $(U_n)_{n\geq 1}\subset W^{1,p(x)}(\Omega)$ with:
			
			\begin{equation}
				d\leq \mathcal{J}_{\lambda}(U_n)\leq d+\dfrac{1}{n},\ \forall\ n\geq 1.
			\end{equation}
			
			\noindent The key step comes here. For any $\epsilon>0$ and any $n\geq 1$ we have that:
			
			\begin{align*}
				d\leq\mathcal{J}_{\lambda}(U_n)&\leq \mathcal{J}_{\lambda,\epsilon}(U_n)=\mathcal{J}_{\lambda}(U_n)+\epsilon\int_{\Omega} \dfrac{|U_n(x)|^{p(x)}}{p(x)}\ dx\\
				&\leq d+\dfrac{1}{n}+\dfrac{\epsilon}{p^{-}}\int_{\Omega} |U_n(x)|^{p(x)}\ dx\\
				&=d+\dfrac{1}{n}+\dfrac{\epsilon}{p^{-}}\rho_{p(x)}(U_n).
			\end{align*}
			
			\noindent Because $U_n\in W^{1,p(x)}(\Omega)\subset L^{p(x)}(\Omega)$ we have that $0\leq\rho_{p(x)}(U_n)<\infty$. So for $\epsilon_n=\dfrac{p^{-}}{n\big (\rho_{p(x)}(U_n)+1\big )}\leq\dfrac{p^{-}}{n}$ we have that:
			
			\begin{equation}
				d\leq\mathcal{J}_{\lambda}(U_n)\leq \mathcal{J}_{\lambda,\epsilon_n}(U_n)\leq d+\dfrac{2}{n},\ \forall\ n\geq 1.
			\end{equation}
			
			\noindent From the \textit{squeezing principle} we obtain that $\lim\limits_{n\to\infty} \mathcal{J}_{\lambda,\epsilon_n}(U_n)=d$. Finally from \eqref{3eqmd}, since $\epsilon_n\to 0^{+}$, and Lemma \ref{3lemgama} we get that:
			
			\begin{equation}
				d=\lim\limits_{n\to\infty} \mathcal{J}_{\lambda,\epsilon_n}(U_n)\geq\lim\limits_{n\to\infty} \mathcal{J}_{\lambda,\epsilon_n}(V_{\epsilon_n})=\mathcal{J}_{\lambda}(V)=m.
			\end{equation}

		\end{proof}
		
			\noindent Since $\mathcal{J}_{\lambda}\in C^1\big (W^{1,p(x)}(\Omega) \big )$ and $V$ is a (global) minimum of $\mathcal{J}_{\lambda}$ as $\mathcal{J}_{\lambda}(V)=m=d=\displaystyle\inf_{U\in W^{1,p(x)}(\Omega)} \mathcal{J}_{\lambda}(U)$, we deduce that $\mathcal{J}_{\lambda}'(V)=0$, i.e. $V$ is a weak solution of the problem \eqref{problemdelambda}.\footnote{See Corollary 2.5 at page 53 from \cite{coleman}.}
	\end{proof}

		\section{Steady-states}
		
		\noindent This section is dedicated to the problem \eqref{eqedg}:
		
		\begin{equation}\tag{$DE$}
			\begin{cases}-\operatorname{div}\mathbf{a}\big (x,\nabla U(x)\big )=f\big (x,U(x)\big ), & x\in\Omega\\[3mm] \mathbf{a}(x,\nabla U)\cdot\nu=0, & x\in \partial\Omega\\[3mm] 0\leq U(x)\leq \delta_0, & x\in\Omega\end{cases}
		\end{equation}
		
		\begin{definition}  We say that $U\in W^{1,p(x)}(\Omega)$ is a \textbf{weak solution} for \eqref{eqedg} if $\delta_0\geq U\geq 0$ a.e. on $\Omega$ and for any test function $\phi\in W^{1,p(x)}(\Omega)$ we have that:
			\begin{equation}\label{weaksoldef}
				\int_{\Omega} \mathbf{a}(x,\nabla U(x))\cdot \nabla\phi(x)\ dx=\int_{\Omega} f(x,U(x))\phi(x)\ dx.
			\end{equation}
			
			 We say that $U\in W^{1,p(x)}(\Omega)$ is a \textbf{weak subsolution} for \eqref{eqedg} if $\delta_0\geq U\geq 0$ a.e. on $\Omega$ and for any test function $\phi\in W^{1,p(x)}(\Omega)^+$ we have that:
			\begin{equation}
				\int_{\Omega} \mathbf{a}(x,\nabla U(x))\cdot \nabla\phi(x)\ dx\leq\int_{\Omega} f(x,U(x))\phi(x)\ dx.
			\end{equation}
			
			 We say that $U\in W^{1,p(x)}(\Omega)$ is a \textbf{weak supersolution} for \eqref{eqedg} if $\delta_0\geq U\geq 0$ a.e. on $\Omega$ and for any test function $\phi\in W^{1,p(x)}(\Omega)^+$ we have that:
			\begin{equation}
				\int_{\Omega} \mathbf{a}(x,\nabla U(x))\cdot \nabla\phi(x)\ dx\geq\int_{\Omega} f(x,U(x))\phi(x)\ dx.
			\end{equation}
		\end{definition}
		
		\begin{remark} It is easy to notice from \textnormal{\textbf{(H10)}} that $0$ is a \textbf{weak subsolution} of $\eqref{eqedg}$ and $\delta_0$ is a \textbf{weak supersolution} of $\eqref{eqedg}$. Indeed:
			
			\begin{equation*}
				\begin{cases}\displaystyle\int_{\Omega} \mathbf{a}(x,\nabla 0)\cdot \nabla\phi(x)\ dx=0\leq \int_{\Omega} f(x,0)\phi(x)\ dx,\ \forall\ \phi\in W^{1,p(x)}(\Omega)^+ \\[5mm]
					\displaystyle\int_{\Omega} \mathbf{a}(x,\nabla \delta_0)\cdot \nabla\phi(x)\ dx=0\geq \int_{\Omega} f(x,\delta_0)\phi(x)\ dx,\ \forall\ \phi\in W^{1,p(x)}(\Omega)^+
				\end{cases}.
			\end{equation*}
		\end{remark}
		
		\begin{remark} Because $U\in W^{1,p(x)}(\Omega)\cap\mathcal{U}\subseteq L^2(\Omega)\cap\mathcal{U}$ we have, from Proposition \ref{prop36} \textbf{(3)}, that $f(\cdot,U(\cdot))\in L^2(\Omega)$. Knowing that $\phi\in W^{1,p(x)}(\Omega)\hookrightarrow L^2(\Omega)$, from \textnormal{\textbf{(H2)}}, we conclude that $f(\cdot,U(\cdot))\phi\in L^1(\Omega)$ and \eqref{weaksoldef} makes sense.
			
		\end{remark}
		
		\begin{proposition}\label{prop71} If $U$ is a weak solution of \eqref{eqedg} then $\displaystyle\int_{\Omega} f\big (x,U(x)\big )\ dx=0$.
		\end{proposition}
		
		\begin{proof} Taking $\phi=1\in W^{1,p(x)}(\Omega)$ as a test function in \eqref{weaksoldef} then:
			
			\begin{equation}
				\int_{\Omega} f\big (x,U(x)\big )\ dx=\int_{\Omega}  f\big (x,U(x)\big )\cdot 1\ dx=\int_{\Omega} \mathbf{a}(x,\nabla U(x))\cdot \nabla 1\ dx=0.
			\end{equation}
			
		\end{proof}
		
		\bigskip
		\bigskip
		
		\section{The unconstrained problem}
		
		\noindent Consider the following unconstrained problem:
		\begin{equation}\tag{$\overline{DE}$}\label{eqedgover}
			\begin{cases}-\operatorname{div}\mathbf{a}\big (x,\nabla U(x)\big )=\overline{f}\big (x,U(x)\big ), & x\in\Omega\\[3mm] \mathbf{a}(x,\nabla U)\cdot\nu=0, & x\in \partial\Omega\end{cases}
		\end{equation}
		
		\begin{proposition}\label{propuover}
			If $U\in W^{1,p(x)}(\Omega)$ is a weak solution for \eqref{eqedgover} then $U\in\mathcal{U}$.

		\end{proposition}
		
		\begin{proof} First we show that $U\geq 0$ a.e. on $\Omega$. Take $\phi=U^{-}\in W^{1,p(x)}(\Omega)^{+}$ as the test function. Then $\nabla \phi=\nabla U^{-}=\begin{cases}0,\ \text{a.e. on}\ \{x\in\Omega\ |\ U(x)\geq 0\}\\[3mm]-\nabla U, \ \text{a.e. on}\ \{x\in\Omega\ |\ U(x)<0\}:=\Omega_0 \end{cases}$. Then:
			
			\begin{equation}0\geq -\int_{\Omega_0} \Phi(x,|\nabla U|)|\nabla U|\ dx=\int_{\Omega_0}\overline{f}(x,U)\cdot (-U)\ dx=\int_{\Omega_0} \underbrace{-f(x,0)U}_{\geq 0}+\dfrac{\lambda_0}{2}U^2\geq 0,
			\end{equation}

			\noindent As $\lambda_0>0$ we get that $\displaystyle\int_{\Omega_0} U^2=0$, i.e. $U\equiv 0$ on $\Omega_0$ (false) or $|\Omega_0|=0$ (i.e. $U^-\equiv 0$). So $U=U^+\geq 0$ a.e. on $\Omega$.
			
			\noindent Now take $\phi=(\delta_0-U)^{-}=\begin{cases} 0, \ \text{a.e. on}\ \{x\in\Omega\ |\ U(x)\leq \delta_0\}\\[3mm] -(\delta_0-U(x)),\ \text{a.e. on}\ \{x\in\Omega\ |\ U(x)>\delta_0\}:=\Omega_1\end{cases}\in W^{1,p(x)}(\Omega)$ as the test function. Notice that:

			\begin{equation}
				\nabla\phi=\nabla (\delta_0-U)^{-}=\begin{cases} 0,\ \text{a.e. on}\ \{x\in\Omega\ |\ U(x)\leq \delta_0\}\\[3mm] -(\nabla \delta_0-\nabla U)=\nabla U,\ \text{a.e. on}\ \Omega_1\end{cases}.
			\end{equation}
			
			\noindent We obtain that:
			
			\begin{equation}
				0\geq -\int_{\Omega_1} \underbrace{\mathbf{a}(x,\nabla U)\cdot\nabla U\ dx}_{\geq 0}=\int_{\Omega_1} \overline{f}(x,U) (\delta_0-U)\ dx=\int_{\Omega_1} \underbrace{\overline{f}(x,\delta_0)(\delta_0-U)}_{\geq 0}+\tilde{\lambda}_0(\delta_0-U)^2\ dx\geq 0.
			\end{equation}

			\noindent Therefore $\displaystyle\int_{\Omega_1}(\delta_0-U)^2\ dx=0$ which shows that $(\delta_0-U)^{-}\equiv 0$ on $\Omega$ ($\lambda>0$). In conclusion $\delta_0-U=(\delta_0-U)^{+}\geq 0$ a.e. on $\Omega$.

		\end{proof}
		
		\begin{remark} Any weak solution of \eqref{eqedgover} is a weak solution of \eqref{eqedg} and vice-versa, meaning that the two problems are equivalent.
		\end{remark}
		
		\begin{definition} We consider the energy functional $\mathcal{J}:W^{1,p(x)}(\Omega)\to\mathbb{R}$ associated to problem \eqref{eqedgover}, that is defined by:
			
			\begin{equation}
				\mathcal{J}(U)=\int_{\Omega} A(x,\nabla U(x))\ dx-\int_{\Omega} \overline{F}(x,U(x))\ dx,\ \forall\ U\in W^{1,p(x)}(\Omega).
			\end{equation}	
		\end{definition}
		
		\begin{remark} From Proposition \ref{prop38} \textnormal{\textbf{(4)}}, we have that $\mathcal{J}\in C^1\big (W^{1,p(x)}(\Omega) \big )$ and:
			
			\begin{equation}
				\langle \mathcal{J}'(U),\phi\rangle=\int_{\Omega} \mathbf{a}(x,\nabla U(x))\cdot\nabla\phi(x)\ dx-\int_{\Omega} \overline{f}(x,U(x))\phi(x)\ dx,\ \forall\ \phi\in W^{1,p(x)}(\Omega).
			\end{equation}
			
			\noindent In particular $U\in W^{1,p(x)}(\Omega)$ is a \textbf{weak solution} of \eqref{eqedgover} iff $\mathcal{J}'(U)\equiv 0$. 
			
		\end{remark}

		\begin{definition}
		Define the following (nonlinear) operator:
		
		\begin{equation}\label{kboperator}
			\mathcal{K}:\mathcal{U}\to \mathcal{U},\ \mathcal{K}(U)=V,\ \text{where}\ \begin{cases}-\operatorname{div}\mathbf{a}(x,\nabla V)+\lambda_0 \overline{b}(x,V(x))=g_{\lambda_0}(x,U(x)), & x\in\Omega\\[3mm] \mathbf{a}(x,\nabla V)\cdot\nu=0, & x\in\partial\Omega \end{cases}.
		\end{equation}
		
		\end{definition}
		
		\begin{remark}
			The operator $\mathcal{K}$ is well-defined, because $g_{\lambda_0}(\cdot, U)=f(\cdot,U)+\lambda_0 b(\cdot,U)\in\mathcal{M}_{\lambda_0}$ as $U\in\mathcal{U}$. Moreover from Proposition \ref{prop41} and Theorem \ref{3thmrot} we have that problem \eqref{kboperator} has a unique solution $V\in W^{1,p(x)}(\Omega)\cap\mathcal{U}$.
		\end{remark}
		
		\begin{proposition}\label{propK} The (nonlinear) operator $\mathcal{K}:\mathcal{U}\to\mathcal{U}$ has the following properties:
		\begin{enumerate}
			\item[\textbf{(1)}] $\mathcal{K}(U)\in W^{1,p(x)}(\Omega)\cap\mathcal{U}$ for any $U\in\mathcal{U}$;
			
			\item[\textbf{(2)}] $\mathcal{K}$ is strictly monotone, i.e. for any $U_1,U_2\in\mathcal{U}$ with $U_1\leq U_2$ a.e. on $\Omega$ and $U_1\not\equiv U_2$ we have that $\mathcal{K}(U_1)\leq \mathcal{K}(U_2)$ a.e. on $\Omega$ and $\mathcal{K}(U_1)\not\equiv\mathcal{K}(U_2)$;
			
			\item[\textbf{(3)}] If $\varepsilon\in [0,\delta_0]$ and $f(x,\varepsilon)\geq 0$ for a.e. $x\in\Omega$ then $\mathcal{K}(\varepsilon)\geq \varepsilon$ a.e. on $\Omega$. Similarly, if $\delta\in [0,\delta_0]$ and $f(x,\delta)\leq 0$ for a.e. $x\in\Omega$ then $\mathcal{K}(\delta)\leq\delta$ a.e. on $\Omega$.
			
			\item[\textbf{(4)}] $\mathcal{K}$ is continuous with respect to the norm of $L^2(\Omega)$, i.e. for any sequence $(U_n)_{n\geq 1}\subset\mathcal{U}$ and $U\in\mathcal{U}$ with $U_n\to U$ in $L^2(\Omega)$ we have that: $\mathcal{K}(U_n)\to \mathcal{K}(U)$ in $L^2(\Omega)$. This result still holds if we replace $L^2(\Omega)$ by $L^r(\Omega)$, for any $r\in [1,\infty)$.
			
		\end{enumerate}
		
		\end{proposition}
		
		\begin{proof}\textbf{(1)} This is obvious since $\mathcal{K}(U)$ is a weak solution of \eqref{kboperator}. Hence $\mathcal{K}(U)\in W^{1,p(x)}(\Omega)\cap\mathcal{U}=X$.
			
			\bigskip
			
			\noindent\textbf{(2)} Denote $V_1=\mathcal{K}(U_1)$ and $V_2=\mathcal{K}(U_2)$. Therefore for any $\phi\in W^{1,p(x)}(\Omega)$:
			
			\begin{equation}
				\begin{cases}\displaystyle\int_{\Omega}\mathbf{a}(x,\nabla V_1)\cdot \nabla\phi\ dx+\lambda_0 \int_{\Omega} b(x,V_1)\phi\ dx=\int_{\Omega} (f(x,U_1)+\lambda_0 b(x,U_1))\phi(x)\ dx\\[3mm] \displaystyle\int_{\Omega}\mathbf{a}(x,\nabla V_2)\cdot \nabla\phi\ dx+\lambda_0 \int_{\Omega} b(x,V_2)\phi\ dx=\int_{\Omega} (f(x,U_2)+\lambda_0 b(x,U_2))\phi(x)\ dx \end{cases}.
			\end{equation}
			
			\noindent Substract the two relations, set $\phi=(V_2-V_1)^-\in W^{1,p(x)}(\Omega)$ and $\Omega_0=\{x\in\Omega\ |\ V_2(x)<V_1(x)\}$. Then we get:
			
			\begin{align*}
				0&\geq -\displaystyle\int_{\Omega_0}\underbrace{\big (\mathbf{a}(x,\nabla V_2)-\mathbf{a}(x,\nabla V_1)\big )\cdot \big (\nabla V_2-\nabla V_1\big )}_{\geq 0}\ dx-\lambda_0 \int_{\Omega_0} \underbrace{\big (b(x,V_2)-b(x,V_1)\big )(V_2-V_1)}_{\geq 0}\ dx\\
				&=\int_{\Omega} \underbrace{\big (g_{\lambda_0}(x,U_2)-g_{\lambda_0}(x,U_1)\big )}_{\geq 0}(V_2-V_1)^{-}\ dx\geq 0.
			\end{align*}
			
			\noindent Thus $\displaystyle\int_{\Omega_0}\big (b(x,V_2)-b(x,V_1)\big )(V_2-V_1)\ dx=0$, which means that $V_1\equiv V_2$ a.e. on $\Omega_0$, i.e. $|\Omega_0|=0$ which is that $V_2\geq V_1$ a.e. on $\Omega$.
			
			\noindent If $V_2\equiv V_1$ then $\displaystyle\int_{\Omega} \big [g_{\lambda_0}(x,U_2(x))-g_{\lambda_0}(x,U_1(x))\big ]\phi(x)\ dx=0$ for all $\phi\in W^{1,p(x)}(\Omega)$. From the \textit{Fundamental theorem of the calculus of variations} we get that $g_{\lambda_0}(\cdot,U_1)=g_{\lambda_0}(\cdot,U_1)$ a.e. on $\Omega$, which is the same as $U_1(x)=U_2(x)$ ($g(x,\cdot)$ is strictly increasing on $[0,\delta_0]$) for a.e. $x\in\Omega$. So $U_1\equiv U_2$. Thus $V_1\not\equiv V_2$ and we are done.	
			
			\noindent\textbf{(3)} Let $V=\mathcal{K}(\varepsilon)$. Then for any $\phi\in W^{1,p(x)}(\Omega)$:
			
			\begin{equation}
				\begin{cases}\displaystyle\int_{\Omega}\mathbf{a}(x,\nabla V)\cdot \nabla\phi\ dx+\lambda_0 \int_{\Omega} b(x,V)\phi\ dx=\int_{\Omega} (f(x,\varepsilon)+\lambda_0 b(x,\varepsilon))\phi(x)\ dx\\[3mm] \displaystyle\int_{\Omega}\mathbf{a}(x,\nabla \varepsilon)\cdot \nabla\phi\ dx+\lambda_0 \int_{\Omega} b(x,\varepsilon)\phi\ dx=\int_{\Omega} \lambda_0 b(x,\varepsilon)\phi(x)\ dx \end{cases}.
			\end{equation}
			
			\noindent Substracting the two relations and setting $\phi=(V-\varepsilon)^-\in W^{1,p(x)}(\Omega)$ and $\Omega_0=\{x\in\Omega\ |\ V(x)<\varepsilon\}$ yields:
			
			\begin{align*}
				0&\geq -\displaystyle\int_{\Omega_0}\underbrace{\big (\mathbf{a}(x,\nabla V)-\mathbf{a}(x,\nabla \varepsilon)\big )\cdot \big (\nabla V-\nabla \varepsilon\big )}_{\geq 0}\ dx-\lambda_0 \int_{\Omega_0} \underbrace{\big (b(x,V)-b(x,\varepsilon)\big )(V-\varepsilon)}_{\geq 0}\ dx\\
				&=\int_{\Omega} \underbrace{f(x,\varepsilon)}_{\geq 0}(V-\varepsilon)^{-}\ dx\geq 0.
			\end{align*}
			
			\noindent Thus $\displaystyle\int_{\Omega_0}\big (b(x,V)-b(x,\varepsilon)\big )(V-\varepsilon)\ dx=0$, which means that $V\equiv \varepsilon$ a.e. on $\Omega_0$ (false), or $|\Omega_0|=0$ which gives that $V\geq \varepsilon$ a.e. on $\Omega$, as needed.
			
			\noindent In the same manner we will show that $V:=\mathcal{K}(\delta)\leq \delta$ a.e. on $\Omega$. We have for each $\phi\in W^{1,p(x)}(\Omega)$:
			
			\begin{equation}
				\begin{cases}\displaystyle\int_{\Omega}\mathbf{a}(x,\nabla \delta)\cdot \nabla\phi\ dx+\lambda_0 \int_{\Omega} b(x,\delta)\phi\ dx=\int_{\Omega} \lambda_0 b(x,\delta)\phi(x)\ dx\\[3mm] \displaystyle\int_{\Omega}\mathbf{a}(x,\nabla V)\cdot \nabla\phi\ dx+\lambda_0 \int_{\Omega} b(x,V)\phi\ dx=\int_{\Omega} \big (f(x,\delta)+\lambda_0 b(x,\delta)\big )\phi(x)\ dx \end{cases}.
			\end{equation}
			
			\noindent Substracting the two relations and setting $\phi=(\delta-V)^-\in W^{1,p(x)}(\Omega)$ and $\Omega_1=\{x\in\Omega\ |\ V(x)>\delta\}$ yields:
			
			\begin{align*}
				0&\geq -\displaystyle\int_{\Omega_1}\underbrace{\big (\mathbf{a}(x,\nabla \delta)-\mathbf{a}(x,\nabla V)\big )\cdot \big (\nabla \delta-\nabla V\big )}_{\geq 0}\ dx-\lambda_0 \int_{\Omega_1} \underbrace{\big (b(x,\delta)-b(x,V)\big )(\delta-V)}_{\geq 0}\ dx\\
				&=\int_{\Omega} \underbrace{-f(x,\delta)}_{\geq 0}(\delta-V)^{-}\ dx\geq 0.
			\end{align*}
			
			\noindent Thus $\displaystyle\int_{\Omega_1}\big (b(x,\delta)-b(x,V)\big )(\delta-V)\ dx=0$, which means that $V\equiv \delta$ a.e. on $\Omega_1$ (false), or $|\Omega_1|=0$ which gives that $V\leq \delta$ a.e. on $\Omega$, as needed.

			\noindent\textbf{(4)} Set $V_n=\mathcal{K}(U_n)$ for each integer $n\geq 1$ and $V=\mathcal{K}(U)$.
			
			\begin{align*}
				\lambda_0\int_{\Omega} \big (b(x,V_n)-b(x,V) \big )(V_n-V)\ dx&\leq \int_{\Omega}\big (\mathbf{a}(x,\nabla V_n)-\mathbf{a}(x,\nabla V) \big )\cdot \big (\nabla V_n-\nabla V\big ) \ dx+\\
				&+\lambda_0\int_{\Omega} \big (b(x,V_n)-b(x,V) \big )(V_n-V)\ dx\\
				&=\int_{\Omega} \big (g_{\lambda_0}(x,U_n)-g_{\lambda_0}(x,U))(V_n-V)\ dx\\
				\text{(Cauchy ineq.)}\ \ \ \ \ \ \ \ &\leq \Vert g_{\lambda_0}(\cdot,U_n)-g_{\lambda_0}(x\cdot,U)\Vert_{L^2(\Omega)}\Vert V_n-V\Vert_{L^2(\Omega)}\\
				&\leq \Vert g_{\lambda_0}(\cdot,U_n)-g_{\lambda_0}(x\cdot,U)\Vert_{L^2(\Omega)}\Vert \delta_0\Vert_{L^2(\Omega)}\\
				&=\delta_0\sqrt{|\Omega|}\cdot\Vert g_{\lambda_0}(\cdot,U_n)-g_{\lambda_0}(x\cdot,U)\Vert_{L^2(\Omega)}\\
				&\leq \delta_0\sqrt{|\Omega|}\cdot\left(\Vert f(\cdot,U_n)-f(\cdot,U)\Vert_{L^2(\Omega)}+\lambda_0\Vert b(\cdot,U_n)-b(\cdot,U)\Vert_{L^2(\Omega)}\right) \\
				&\stackrel{n\to\infty}{\longrightarrow} 0,
			\end{align*}
			
			\noindent from Proposition \ref{propprel1} \textbf{(1)} and Proposition \ref{prop35} \textbf{(3)}. So we have proved that $\big (b(\cdot,V_n)-b(\cdot,V)\big )(V_n-V)\to 0$ in $L^1(\Omega)$. Take any subsequence $(V_{n_k})_{k\geq 1}$. Therefore $\big (b(\cdot,V_{n_k})-b(\cdot,V)\big )(V_{n_k}-V)\to 0$ in $L^1(\Omega)$. Using a Corollary of the Riez-Fischer Theorem (see \cite[page 234]{Jones}) we get that there is a further subsequence that is pointwise convergent to zero, i.e. $\big (b(\cdot,V_{n_{k_\ell}})-b(\cdot,V)\big )(V_{n_{k_\ell}}-V)\to 0$ pointwise a.e. on $\Omega$.

			\noindent  Next, we show that $V_{n_{k_\ell}}\to V$ pointwise a.e. on $\Omega$. Suppose the contrary. Then there is a set $\Omega_0\subset\Omega$ with $|\Omega_0|>0$ such that $|V_{n_{k_\ell}}(x)-V(x)|\nrightarrow 0$ for each $x\in\Omega_0$. Therefore there is for each $x\in\Omega_0$ an $\varepsilon_x>0$ and a further subsequence (still denoted by $V_{n_{k_\ell}}(x)$) such that $|V_{n_{k_\ell}}(x)-V(x)|>\varepsilon_x$ for any $\ell\geq 1$. There are two cases:
			
			\begin{itemize}
				\item If $V_{n_{k_\ell}}(x)>V(x)+\varepsilon_x$ then from the strict monotonicity of $b(x,\cdot)$ we get 
				
				$$\big (b(x,V_{n_{k_\ell}}(x))-b(x,V(x))\big )(V_{n_{k_\ell}}(x)-V(x))>\underbrace{\big (b(x,V(x)+\varepsilon_x)-b(x,V(x))\big )\varepsilon_x}_{>0}\nrightarrow 0,\ \text{as}\ \ell\to\infty.$$
				
				\item If $V_{n_{k_\ell}}(x)<V(x)-\varepsilon$ then from the monotonicity of $b(x,\cdot)$ we get 
				
				$$\big (b(x,V_{n_{k_\ell}}(x))-b(x,V(x))\big )(V_{n_{k_\ell}}(x)-V(x))>\underbrace{\big (b(x,V(x))-b(x,V(x)-\varepsilon_x)\big )\varepsilon_x}_{>0}\nrightarrow 0,\ \text{as}\ \ell\to\infty.$$
			\end{itemize}
			
			\noindent Thus we obtained the desired contradiction: $(b(x,V_{n_{k_\ell}}(x))-b(x,V(x))\big )(V_{n_{k_\ell}}(x)-V(x))\nrightarrow 0$ for any $x\in\Omega_0$. Up to this point we have that $V_{n_{k_\ell}}\to V$ pointwise. But sice $|V_{n_{k_\ell}}|\leq\delta_0\in L^r(\Omega)$, for any $r\in [1,\infty)$, we obtain from \textbf{Lebesgue dominated convergence theorem} that $V_{n_{k_\ell}}\to V$ in $L^r(\Omega)$, i.e. $\Vert V_{n_{k_\ell}}-V\Vert_{L^r(\Omega)}\to 0$. So for the sequence $\big (\Vert V_n-V\Vert_{L^r(\Omega)}\big )_{n\geq 1}$ we have proved that any subsequence of it has a further subsequence that converges to $0$.

			\noindent We give without proof the following elementary lemma:
			
			\begin{lemma}
				A sequence of real numbers converges to a limit $\ell\in\overline{\mathbb{R}}$ if and only if every subsequence has a further subsequence which converges to $\ell$.
			\end{lemma}
			
			\noindent Applying this lemma to the sequence of real numbers $\Vert V_n-V\Vert_{L^r(\Omega)}$ gives us that $V_n\to V$ in $L^r(\Omega)$ as desired.
			
		\end{proof}

		\begin{theorem}\label{theorem81} Problem \eqref{eqedg} admits at least one weak solution. Also there is a maximal weak solution $\overline{U}$ and a minimal weak solution $\underline{U}$, meaning that any weak solution $U$ of \eqref{eqedg} satisfies $\underline{U}\leq U\leq \overline{U}$ almost everywhere on $\Omega$. If $\underline{U}\equiv\overline{U}$ then \eqref{eqedg} has a unique solution.
		\end{theorem}
		
		\begin{proof}\noindent\textbf{(Existence)} We introduce the following recurrence:
			
			\begin{equation}
				\begin{cases}U_{n+1}=\mathcal{K}(U_n),\ n\geq 0\\[3mm] U_0\in\mathcal{U}\end{cases}.
			\end{equation}
			
			\noindent $\bullet$ If we take $U_0\equiv \delta_0\in\mathcal{U}$ then $U_1=\mathcal{K}(U_0)\in\mathcal{U}\ \Rightarrow\ U_1\leq \delta_0=U_0$ a.e. on $\Omega$. Using the monotonicity of $\mathcal{K}$ we get immediately by induction that $U_2=\mathcal{K}(U_1)\leq \mathcal{K}(U_0)=U_1$, ... , $U_{n+1}=\mathcal{K}(U_n)\leq \mathcal{K}(U_{n-1})=U_n$ for any $n\geq 1$. Also, from $\delta_0\equiv U_0\geq 0$ a.e. on $\Omega$ we obtain that $U_1=\mathcal{K}(U_0)\geq \mathcal{K}(0)\geq 0$ because $\mathcal{K}(0)\in\mathcal{U}$. Repeating the same process $U_2=\mathcal{K}(U_1)\geq \mathcal{K}(0)\geq 0$ and in general we get that $U_n\geq 0$ for any $n\geq 1$. Hence:
			
			\begin{equation}
				\delta_0\equiv U_0\geq U_1\geq U_2\geq\dots\geq U_n\geq U_{n+1}\dots\geq 0
			\end{equation}
			
			\noindent So, for a.e. $x\in\Omega$, the sequence of real numbers $\big (U_n(x)\big )_{n\geq 1}\subset [0,\delta_0]$ is decreasing and bounded. We conclude that it is convergent to some limit denoted by $\overline{U}(x)\in [0,\delta_0]$. We can state that $U_n\to\overline{U}$ pointwise a.e. on $\Omega$. Now, since $U_n=|U_n|\leq \delta_0\in L^2(\Omega)$ ($\Omega$ is bounded) for any $n\geq 1$ we obtain from \textit{Lebesgue dominated convergence theorem} that $U_n\to \overline{U}$ in $L^2(\Omega)$.

			\noindent Using Proposition \ref{propK} \textbf{(4)} we get that $\mathcal{K}(U_n)\to\mathcal{K}(\overline{U})$ in $L^2(\Omega)$. But $\mathcal{K}(U_n)=U_{n+1}\to \overline{U}$ in $L^2(\Omega)$. From the uniqueness of the limit we conclude that $\mathcal{K}(\overline{U})=\overline{U}$. So $\overline{U}\in W^{1,p(x)}(\Omega)\cap \mathcal{U}=X$ and $\overline{U}$ is a weak solution of the problem \eqref{eqedg}.
			
			\begin{remark} If $\delta_0\not\equiv\mathcal{K}(\delta_0)$ then $\delta_0\gneq U_1\gneq U_2\gneq\dots\gneq U_{n}\gneq U_{n+1}\gneq\dots\gneq \overline{U}$. This follows from the strict monotonicity of $\mathcal{K}$.
			\end{remark}
			
			\noindent $\bullet$ If we set $U_0\equiv 0$ then in the same manner as above we will get that:
			
			\begin{equation}
				\delta_0\geq\dots U_{n+1}\geq U_n\geq\dots\geq U_2\geq U_1\geq U_0\equiv 0.
			\end{equation}
			
			\noindent Thus for a.e. $x\in\Omega$ the sequence $\big (U_n(x) \big )_{n\geq 1}\subset [0,\delta_0]$ is increasing and bounded, hence convergent to some $\underline{U}(x)\in [0,\delta_0]$. So $U_n\to\underline{U}$ pointwise a.e. on $\Omega$. From $|U_n|=U_n\leq \delta_0\in L^2(\Omega)$ we deduce from the \textit{Lebesgue dominated convergence theorem} that $U_n\to\underline{U}$ in $L^2(\Omega)$. 
			
			\noindent Now from the continuity of $\mathcal{K}$ with respect to the $L^2(\Omega)$-norm we conclude that $U_{n+1}=\mathcal{K}(U_n)\to\mathcal{K}(\underline{U})$ in $L^2(\Omega)$. But since $U_{n+1}\to \underline{U}$ in $L^2(\Omega)$ we deduce that $\underline{U}=\mathcal{K}(\underline{U})$. Thus $\underline{U}\in W^{1,p(x)}(\Omega)\cap \mathcal{U}=X$ is also a weak solution of the problem \eqref{eqedg}.
			
			\begin{remark} If $0\not\equiv\mathcal{K}(0)$ then $0\lneq U_1\lneq U_2\lneq\dots\lneq U_{n}\lneq U_{n+1}\lneq\dots\lneq \underline{U}$. 
			\end{remark}
			
			\bigskip
			
			\noindent Now if $U\in X$ is any weak solution of the problem \eqref{eqedg} then: $\underline{U}(x)\leq U(x)\leq\overline{U}(x)$ for a.e. $x\in\Omega$. Indeed, from $U\in\mathcal{U}$, we have that: $0\leq U\leq \delta_0$ a.e. on $\Omega$. Applying $\mathcal{K}$ we obtain: $\mathcal{K}(0)\leq\mathcal{K}(U)=U\leq \mathcal{K}(\delta_0)$. Repeating this process $n\geq 1$ times we get:
				
				\begin{equation}
					\mathcal{K}^n(0)\leq U\leq\mathcal{K}^n(\delta_0).
				\end{equation}
				
				\noindent Making $n\to\infty$ we get that $\underline{U}\leq U\leq\overline{U}$ a.e. on $\Omega$. Of course, if $\underline{U}=\overline{U}$ then problem \eqref{eqedg} has exactly one weak solution.	
		\end{proof}
		
			\begin{remark} If $f(x,0)\equiv 0$ then $0\equiv\mathcal{K}(0)$ and therefore $\underline{U}\equiv 0$. This will be called the \textbf{trivial solution} of the problem \eqref{eqedg}.
		\end{remark}

		\begin{remark} If $U,\widetilde{U}\in X$ are two weak solutions of \eqref{eqedg} then the following relation holds:\footnote{Starting with this relation I have proved the uniqueness of the steady-state for the Neumann laplacian in the paper \cite{max1}.}
			
			\begin{align}
				\int_{\Omega} \big (\Psi(x,|\nabla\widetilde{U}|)-\Psi(x,|\nabla U|)\big )\nabla\tilde{U}\cdot\nabla U\ dx&=\int_{\Omega} g(x,\widetilde{U})U-g(x,U)\tilde{U}\ dx\nonumber \\
				&=\int_{\Omega} f(x,\widetilde{U})U-f(x,U)\tilde{U}\ dx.
			\end{align}
			
		\end{remark}

		\section{Solutions in a given subinterval}
		
		\noindent If there are some $0\leq\varepsilon\leq\delta\leq\delta_0$ such that $\begin{cases}f(x,\varepsilon)\geq 0, \ \text{for a.e.}\ x\in\Omega\\[3mm] f(x,\delta)\leq 0,\ \text{for a.e.}\ x\in\Omega \end{cases}$ we may consider the following problem:
		
		\begin{equation}\tag{$E_{[\varepsilon,\delta]}$}\label{eqedgepsdelt}
			\begin{cases}-\operatorname{div}\mathbf{a}\big (x,\nabla U(x)\big )=f\big (x,U(x)\big ), & x\in\Omega\\[3mm] \mathbf{a}(x,\nabla U)\cdot\nu=0, & x\in \partial\Omega\\[3mm] \varepsilon\leq U(x)\leq \delta, & x\in\Omega\end{cases}
		\end{equation}
		
		\begin{definition} We say that $U\in W^{1,p(x)}(\Omega)$ is a \textbf{weak solution} for \eqref{eqedgepsdelt} if $U\in\mathcal{U}_{[\varepsilon,\delta]}$ and for any test function $\phi\in W^{1,p(x)}(\Omega)$ we have that:
			\begin{equation}
				\int_{\Omega} \mathbf{a}(x,\nabla U(x))\cdot \nabla\phi(x)\ dx=\int_{\Omega} f(x,U(x))\phi(x)\ dx.
			\end{equation}
		\end{definition}
		
		\begin{theorem}\label{theorem91}
			Problem \eqref{eqedgepsdelt} admits at least one weak solution. Furthermore there is a minimal and a maximal weak solution of \eqref{eqedgepsdelt} denoted by $\underline{U}_{[\varepsilon,\delta]}$ and $\overline{U}_{[\varepsilon,\delta]}$, i.e. any other weak solution $U$ of \eqref{eqedgepsdelt} satisfies $\underline{U}_{[\varepsilon,\delta]}\leq U\leq \overline{U}_{[\varepsilon,\delta]}$ a.e. on $\Omega$.
		\end{theorem}
		
		\begin{proof} From Proposition \ref{propK} \textbf{(3)} we can work with the restriction of the operator $\mathcal{K}$ on $\mathcal{U}_{[\varepsilon,\delta]}$, that is: $\mathcal{K}:\mathcal{U}_{[\varepsilon,\delta]}\to\mathcal{U}_{[\varepsilon,\delta]}$.
			
		\noindent We introduce the following recurrence:
		
		\begin{equation}
			\begin{cases}U_{n+1}=\mathcal{K}(U_n),\ n\geq 0\\[3mm] U_0\in\mathcal{U}_{[\varepsilon,\delta]}\end{cases}.
		\end{equation}
		
		\noindent $\bullet$ If we take $U_0\equiv \delta\in\mathcal{U}_{[\varepsilon,\delta]}$ then $U_1=\mathcal{K}(U_0)\in\mathcal{U}_{[\varepsilon,\delta]}\ \Rightarrow\ U_1\leq \delta=U_0$ a.e. on $\Omega$. Using the monotonicity of $\mathcal{K}$ we get immediately by induction that $U_2=\mathcal{K}(U_1)\leq \mathcal{K}(U_0)=U_1$, ... , $U_{n+1}=\mathcal{K}(U_n)\leq \mathcal{K}(U_{n-1})=U_n$ for any $n\geq 1$. Also, from $\delta\equiv U_0\geq \varepsilon$ a.e. on $\Omega$ we obtain that $U_1=\mathcal{K}(U_0)\geq \mathcal{K}(\varepsilon)\geq \varepsilon$ because $\mathcal{K}(\varepsilon)\in\mathcal{U}_{[\varepsilon,\delta]}$. Repeating the same process $U_2=\mathcal{K}(U_1)\geq \mathcal{K}(\varepsilon)\geq \varepsilon$ and in general we get that $U_n\geq \varepsilon$ for any $n\geq 1$. Hence:
		
		\begin{equation}
			\delta\equiv U_0\geq U_1\geq U_2\geq\dots\geq U_n\geq U_{n+1}\dots\geq \varepsilon.
		\end{equation}
		
		\noindent So, for a.e. $x\in\Omega$, the sequence of real numbers $\big (U_n(x)\big )_{n\geq 1}\subset [\varepsilon,\delta]$ is decreasing and bounded. We conclude that it is convergent to some limit denoted by $\overline{U}_{[\varepsilon,\delta]}(x)\in [\varepsilon,\delta]$. We can state that $U_n\to\overline{U}_{[\varepsilon,\delta]}$ pointwise a.e. on $\Omega$. Now, since $U_n=|U_n|\leq \delta\in L^2(\Omega)$ ($\Omega$ is bounded) for any $n\geq 1$ we obtain from \textit{Lebesgue dominated convergence theorem} that $U_n\to \overline{U}_{[\varepsilon,\delta]}$ in $L^2(\Omega)$.

		\noindent Using Proposition \ref{propK} \textbf{(4)} we get that $\mathcal{K}(U_n)\to\mathcal{K}(\overline{U}_{[\varepsilon,\delta]})$ in $L^2(\Omega)$. But $\mathcal{K}(U_n)=U_{n+1}\to \overline{U}_{[\varepsilon,\delta]}$ in $L^2(\Omega)$. From the uniqueness of the limit we conclude that $\mathcal{K}(\overline{U}_{[\varepsilon,\delta]})=\overline{U}_{[\varepsilon,\delta]}$. So $\overline{U}_{[\varepsilon,\delta]}\in W^{1,p(x)}(\Omega)\cap \mathcal{U}_{[\varepsilon,\delta]}$ and $\overline{U}_{[\varepsilon,\delta]}$ is a weak solution of the problem \eqref{eqedgepsdelt}.
		
		\begin{remark} If $\delta\not\equiv\mathcal{K}(\delta_0)$ then $\delta\gneq U_1\gneq U_2\gneq\dots\gneq U_{n}\gneq U_{n+1}\gneq\dots\gneq \overline{U}_{[\varepsilon,\delta]}$. This follows from the strict monotonicity of $\mathcal{K}$.
		\end{remark}
		
		\noindent $\bullet$ If we set $U_0\equiv \varepsilon$ then in the same manner as above we will get that:
		
		\begin{equation}
			\delta\geq\dots U_{n+1}\geq U_n\geq\dots\geq U_2\geq U_1\geq U_0\equiv \varepsilon.
		\end{equation}
		
		\noindent Thus for a.e. $x\in\Omega$ the sequence $\big (U_n(x) \big )_{n\geq 1}\subset [\varepsilon,\delta]$ is increasing and bounded, hence convergent to some $\underline{U}_{[\varepsilon,\delta]}(x)\in [\varepsilon,\delta]$. So $U_n\to\underline{U}_{[\varepsilon,\delta]}$ pointwise a.e. on $\Omega$. From $|U_n|=U_n\leq \delta\in L^2(\Omega)$ we deduce from the \textit{Lebesgue dominated convergence theorem} that $U_n\to\underline{U}_{[\varepsilon,\delta]}$ in $L^2(\Omega)$. 
		
		\noindent Now from the continuity of $\mathcal{K}$ with respect to the $L^2(\Omega)$-norm we conclude that $U_{n+1}=\mathcal{K}(U_n)\to\mathcal{K}(\underline{U}_{[\varepsilon,\delta]})$ in $L^2(\Omega)$. But since $U_{n+1}\to \underline{U}_{[\varepsilon,\delta]}$ in $L^2(\Omega)$ we deduce that $\underline{U}_{[\varepsilon,\delta]}=\mathcal{K}(\underline{U}_{[\varepsilon,\delta]})$. Thus $\underline{U}_{[\varepsilon,\delta]}\in W^{1,p(x)}(\Omega)\cap \mathcal{U}_{[\varepsilon,\delta]}$ is also a weak solution of the problem \eqref{eqedgepsdelt}.
		
		\begin{remark} If $\varepsilon\not\equiv\mathcal{K}(\varepsilon)$ then $\varepsilon\lneq U_1\lneq U_2\lneq\dots\lneq U_{n}\lneq U_{n+1}\lneq\dots\lneq \underline{U}_{[\varepsilon,\delta]}$. 
		\end{remark}
		
		\bigskip
		
		\noindent Now if $U\in W^{1,p(x)}\cap\mathcal{U}_{[\varepsilon,\delta]}$ is any weak solution of the problem \eqref{eqedgepsdelt} then: $\underline{U}_{[\varepsilon,\delta]}(x)\leq U(x)\leq\overline{U}_{[\varepsilon,\delta]}(x)$ for a.e. $x\in\Omega$. Indeed, from $U\in\mathcal{U}_{[\varepsilon,\delta]}$, we have that: $\varepsilon\leq U\leq \delta$ a.e. on $\Omega$. Applying $\mathcal{K}$ we obtain: $\mathcal{K}(\varepsilon)\leq\mathcal{K}(U)=U\leq \mathcal{K}(\delta)$. Repeating this process $n\geq 1$ times we get:
		
		\begin{equation}
			\mathcal{K}^n(\varepsilon)\leq U\leq\mathcal{K}^n(\delta).
		\end{equation}
		
		\noindent Making $n\to\infty$ we get that $\underline{U}_{[\varepsilon,\delta]}\leq U\leq\overline{U}_{[\varepsilon,\delta]}$ a.e. on $\Omega$. Of course, if $\underline{U}_{[\varepsilon,\delta]}=\overline{U}_{[\varepsilon,\delta]}$ then problem \eqref{eqedgepsdelt} has exactly one weak solution.	
		\end{proof}

		\section{Uniqueness of the steady state}
		
		\subsection*{General results}
		
		\begin{theorem} If for a.e. $x\in\Omega$ the function $[\varepsilon,\delta]\ni s\mapsto f(x,s)$ is strictly decreasing then problem \eqref{eqedgepsdelt} admits exactly one weak solution.
		\end{theorem}
		
		\begin{proof} Let $U_1,U_2\in\mathcal{U}_{[\varepsilon,\delta]}\cap W^{1,p(x)}(\Omega)$ be two weak solutions of \eqref{eqedgepsdelt}. Therefore for any $\phi\in W^{1,p(x)}(\Omega)$ we have that:
			
		\begin{equation}
			\begin{cases}\displaystyle\int_{\Omega} \mathbf{a}(x,\nabla U_1)\cdot\nabla\phi\ dx=\int_{\Omega} f(x,U_1)\phi\ dx \\[3mm]\displaystyle\int_{\Omega} \mathbf{a}(x,\nabla U_2)\cdot\nabla\phi\ dx=\int_{\Omega} f(x,U_2)\phi\ dx \end{cases}.
		\end{equation}
		
		\noindent Substracting these relations leads us to: 
		
		\begin{equation}
		\displaystyle\int_{\Omega} \big [\mathbf{a}(x,\nabla U_1)-\mathbf{a}(x,\nabla U_2)\big ]\cdot\nabla\phi\ dx=\int_{\Omega} \big [f(x,U_1)-f(x,U_2)\big ]\phi\ dx.
		\end{equation}
		
		\noindent Choosing the test function $\phi=U_1-U_2\in W^{1,p(x)}(\Omega)$ we get that:
		
		\begin{equation}
			0\leq \int_{\Omega} \big [\mathbf{a}\big (x,\nabla U_1\big )-\mathbf{a}\big (x,\nabla U_2\big )\big ]\cdot \big (\nabla U_1-\nabla U_2\big )\ dx=\int_{\Omega} \big [f(x,U_1)-f(x,U_2)\big ]\big (U_1-U_2\big )\ dx\leq 0.
		\end{equation}
		
		\noindent Thus $U_1\equiv U_2$ and the uniqueness follows.
			
		\end{proof}
		
		\noindent The following important uniqueness result is taken from \cite[Theorem 6.2]{max3}.
		
\begin{theorem}[\textbf{Uniqueness}]\label{uniq} The following results hold:
	
	\begin{enumerate}
		\item[\textbf{(1)}] If hypothesis \textnormal{\textbf{(EH)}} holds and $U_1,U_2$ are two distinct weak solutions ($U_1\not\equiv U_2$) of \eqref{eqedg} that are \textbf{strongly positive}, i.e. $\underset{\Omega}{\operatorname{ess\ inf}}\ U_1$, $ \underset{\Omega}{\operatorname{ess\ inf}}\ U_2>0$, then there is a constant $\lambda>0$ such that:
		
		\medskip
		
		\begin{enumerate}
			\item[$\bullet$] $U_2(x)=\lambda U_1(x)$, for a.e. $x\in\Omega$,
			
			\medskip
			
			\item[$\bullet$] $\Phi\big (x,\lambda|\nabla U_1(x)|\big )=\lambda^{\alpha-1}\Phi\big (x,|\nabla U_1(x)|\big ),\ \text{for a.e.}\ x\in\Omega$,
			
			\medskip
			
			\item[$\bullet$] $f\big (x,\lambda U_1(x)\big )=\lambda^{\alpha-1} f\big (x,U_1(x)\big )$, for a.e. $x\in\Omega$.
		\end{enumerate} 
		
		\bigskip
		
		\item[\textbf{(2)}] If hypothesis \textnormal{\textbf{(EH$_\Phi$)}} holds then problem \eqref{eqedg} has at most one weak solution $U$ that is \textbf{strongly positive}.
		
		\bigskip
		
		\item[\textbf{(3)}] If hypothesis \textnormal{\textbf{(EH$_f$)}} holds then problem \eqref{eqedg} has at most one weak solution $U$ that is \textbf{strongly positive}.
	\end{enumerate}
\end{theorem}

		\subsection*{The source $f$ is independent of the spatial variable}
		
		\noindent In this case $f:[0,\delta_0]\to\mathbb{R}$ is just a continuous functions and the problem \eqref{eqedg} becomes:
		
		\begin{equation}\label{eqedgomo}
			\begin{cases}-\operatorname{div}\mathbf{a}\big (x,\nabla U(x)\big )=f\big (U(x)\big ), & x\in\Omega\\[3mm] \mathbf{a}(x,\nabla U)\cdot\nu=0, & x\in \partial\Omega\\[3mm] 0\leq U(x)\leq \delta_0, & x\in\Omega\end{cases}
		\end{equation}
		
		\noindent In this case from Proposition \ref{prop71} we have that for any weak solution $U$ of \eqref{eqedgomo}:
		
		\begin{equation}\label{ecuatiedebaza}
			\int_{\Omega} f\big (U(x)\big )\ dx=0.
		\end{equation}
		
		\begin{proposition}
			If $f(s)\neq 0$ on $[0,\delta_0]$ then \eqref{eqedgomo} has no weak solution.
		\end{proposition}
		
		\begin{proof} From the continuity of $f$ we deduce that $f$ has constant sign on $[0,\delta_0]$ (either strictly positive or strictly negative). If $U$ is a weak solution of \eqref{eqedgomo} then from \eqref{ecuatiedebaza} we have that $\displaystyle\int_{\Omega} f\big (U(x)\big )\ dx=0$, which is impossible to hold in both cases.
			
		\end{proof}
		
		\begin{proposition}
			If for some constant $c\in [0,\delta_0]$ we have that $f(c)=0$, then $U\equiv c$ is a weak solution of \eqref{eqedgomo}. In other words any zero of $f$ is a weak solution \eqref{eqedgomo}.
		\end{proposition}
		
		\noindent Consider that $c_1,c_2,\hdots,c_n$ are the only zeros of $f$ and $0\leq c_1< c_2<\hdots<c_{n}\leq\delta_0$. Being continuous we have that on each of the intervals $(c_1,c_2),\ (c_2,c_3),\hdots, (c_{n-1},c_n)$ the function $f$ has constant sign. In this context we have the following result:

		\begin{theorem}
			If there is some $i\in\overline{1,n}$ such that $f(s)\geq 0$ on $[0,c_i]$ and $f(s)\leq 0$ on $[c_i,\delta_0]$ then $\overline{U}\equiv c_n$ and $\underline{U}=c_1$. Moreover there are no \textbf{non-constant} solutions of \eqref{eqedgomo} with $U\in [0,c_i]$ or $U\in [c_i,\delta_0]$. In particular if $n=1$, $U\equiv c_1$ is the \textbf{unique} solution of \eqref{eqedgomo}.
		\end{theorem}
		
		\begin{proof} The following sign table of $f$ is valid:
			
				\begin{equation}
				\begin{array}{l|lllllllllllllll}
					s	& 0 & \quad & c_1 & \quad & c_2 & \hdots  & c_{i-1} & \quad & c_i & \quad & c_{i+1} &\hdots & c_n & \quad & \delta_0 \\ \hline
					f(s)	& + & + & 0 &  + & 0 & \hdots & 0 & + & 0 & - &0 & \hdots & 0 & - & -
				\end{array}
			\end{equation}
			
			\noindent From $f(c_i)=0\geq 0$ it follows that $\mathcal{K}(c_i)=c_i\geq c_i$. Therefore, from the proof of the Theorem \ref{theorem91} for $\varepsilon=c_i$ we get $\overline{U}_{[c_i,\delta_0]}=\overline{U}_{[0,\delta_0]}=\overline{U}$. This shows that $\overline{U}\geq c_i$ a.e. on $\Omega$. From the sign table of $f$ we conclude that $f(\overline{U})\leq 0$ a.e. on $\Omega$. Since $\overline{U}$ is a weak solution of \eqref{eqedgomo} we have that $\displaystyle\int_{\Omega} f\big (\overline{U}\big )\ dx=0$. Combining the last two facts yield $f(\overline{U})\equiv 0$ which means that $\overline{U}(x)\in\{c_i,c_{i+1},\hdots, c_n\}$ for a.e. $x\in\Omega$ (constant on portions). Because $\overline{U}\in W^{1,p(x)}(\Omega)$, we get that $\nabla \overline{U}\equiv 0$. The key fact is that $\Omega$ is connected, which shows that $\overline{U}$ must be a constant (and just one) a.e. on $\Omega$. Thus $\overline{U}\equiv c_j$ for some $j\in\{i,i+1,...,n\}$. Being the maximal solution we get that $\overline{U}\equiv c_n$.

			\noindent Similarly, from the proof of the Theorem \ref{theorem91} but for $\delta=c_i$ we get $\underline{U}_{[0,c_i]}=\underline{U}_{[0,\delta_0]}=\underline{U}$. This shows that $\underline{U}\leq c_i$ a.e. on $\Omega$. From the sign table of $f$ we conclude that $f(\underline{U})\geq 0$ a.e. on $\Omega$. Since $\underline{U}$ is a weak solution of \eqref{eqedgomo} we have that $\displaystyle\int_{\Omega} f\big (\underline{U}\big )\ dx=0$. Combining the last two facts yield $f(\underline{U})\equiv 0$ which means that $\underline{U}(x)\in\{c_1,c_{2},\hdots, c_i\}$ for a.e. $x\in\Omega$. Because $\underline{U}\in W^{1,p(x)}(\Omega)$, we get that $\nabla \underline{U}\equiv 0$. Using again that $\Omega$ is connected shows that $\underline{U}$ must be constant a.e. on $\Omega$. Thus $\underline{U}\equiv c_j$ for some $j\in\{1,2,...,i\}$. Being the minimal solution we get that $\underline{U}\equiv c_1$. Therefore if $n=1$ the $\underline{U}\equiv\overline{U}\equiv c_1$, which means that $U\equiv c_1$ is the unique solution of \eqref{eqedgomo}.
			
			\noindent Suppose that $U\in [0,c_i]$ is a weak solution of \eqref{eqedgomo}. Since $f(U(x))\geq 0$ a.e. on $\Omega$ and $\displaystyle\int_{\Omega} f\big (\underline{U}\big )\ dx=0$, it follows as above that $U\equiv c_j$ for some $j\in\{1,2,...,i\}$. In the same fashion if $U\in [c_i,\delta_0]$ is a weak solution of \eqref{eqedgomo} then $U\equiv c_j$ for some $j\in\{i,i+1,...,n\}$. The conclusion follows.
		\end{proof}
		
			\begin{remark} Problem \eqref{eqedg} can have more than one nonconstant solution having also constant solutions. An example would be to take $f:\Omega\times [0,\delta_0]\to\mathbb{R}$,
			
			\begin{equation}
				f(x,s)=s\big (c_1(x)-s\big )\big (c_2(x)-s\big )\big (c_3(x)-s\big),
			\end{equation}
			
			\noindent where $c_1,c_2,c_3$ are measurable nonconstant functions for which there are some constants $\alpha_1,\alpha_2,\alpha_3$ with the property that:
			
			\begin{equation}
				0<\alpha_1<c_1(x)<\alpha_2<c_2(x)<\alpha_3<c_3(x)<\delta_0,\ \text{for a.e.}\ x\in\Omega.
			\end{equation}
			
			\noindent The following sign table of $f$ is valid for a.e. $x\in\Omega$:
			
			\begin{equation}
				\begin{array}{l|lllllllll}
					s	& 0 & \quad & \alpha_1 & c_1(x) & \alpha_2 & c_2(x)  & \alpha_3 & c_3(x) &\delta_0 \\ \hline
					f(x,s)	& 0 &  & + & 0  & - &  & + & & -
				\end{array}
			\end{equation}
			
			\noindent Therefore $U_0\equiv 0$ is a constant solution of \eqref{eqedg}, and from Theorem \ref{theorem91} we get that there is a solution $U_1\in [\alpha_1,\alpha_2]$ and a solution $U_2\in [\alpha_3,\delta_0]$. The last two solutions $U_1$ and $U_2$ cannot be constant, because otherwise this will imply that $f(\cdot,U_1)\equiv 0$ or $f(\cdot, U_2)\equiv 0$. This means that $U_1\equiv c_1$ or $U_2\equiv c_3$ (both $c_1$ and $c_3$ are non-constant functions) -- hence contradiction.
			
		\end{remark}

		\subsection*{The symetry of the source generates extra solutions}
		
		\begin{remark}
			Let the relation $f(x,\delta_0-s)=-f(x,s)$ hold for any $s\in [0,\delta_0]$ and for a.e. $x\in\Omega$. Knowing that $U$ is a weak solution of \eqref{eqedg} we also have that $\delta_0-U$ is a weak solution of \eqref{eqedg}. In particular $U\equiv\dfrac{\delta_0}{2}$ is a weak solution of \eqref{eqedg}.
		\end{remark}

		\section{Example: A multiple-phase problem}

		\noindent Consider the following general multiple phase problem for the signomial source $f$ given in Proposition \ref{proposignom}:

		\begin{equation}\tag{$MP$}\label{eqmphase2}
			\begin{cases}-\Delta_{p_1(x),w_1}U-\Delta_{p_2(x),w_2}U-\hdots-\Delta_{p_{\ell}(x),w_{\ell}}U=f\big (x,U(x)\big ), & x\in\Omega\\[3mm] \mathbf{a}(x,\nabla U)\cdot\nu=0, & x\in \partial\Omega\\[3mm] 0\leq U(x)\leq \delta_0, & x\in\Omega\end{cases}
		\end{equation}
		
		\noindent where $-\Delta_{p_{k}(x),w_k}U:=-\operatorname{div}\big (w_k(x)|\nabla U|^{p_k(x)-2}\nabla U \big )$ for each $k\in\overline{1,\ell}$ and:
		
		\begin{enumerate}
			\item[$\bullet$] $p_1,p_2,\dots, p_{\ell}:\overline{\Omega}\to (1,\infty),\ \ell\geq 1$ are continuous variable exponents such that $p:\overline{\Omega}\to (1,\infty),\ p(x)=\max\{p_1(x),p_2(x),\hdots,p_{\ell}(x)\}$ satisfies the inequality $p^-:=\displaystyle\min_{x\in\overline{\Omega}} p(x)>\dfrac{2N}{N+2}$, or
			
			\item[$\bullet$] $p_1,p_2,\dots, p_{\ell}:\overline{\Omega}\to (1,\infty),\ \ell\geq 1$ are log-H\"{o}lder continuous variable exponents such that $p:\overline{\Omega}\to (1,\infty),\ p(x)=\max\{p_1(x),p_2(x),\hdots,p_{\ell}(x)\}$ satisfies the inequality $p^-:=\displaystyle\min_{x\in\overline{\Omega}} p(x)\geq\dfrac{2N}{N+2}$.\footnote{It it important to point out that if $p_1,p_2,\hdots, p_{\ell}$ are all continuous\textbackslash log-H\"{o}lder continuous\textbackslash Lipschitz continuous, the so is $p$. To see why it is so, just note that the following elementary inequality holds: $|p(x)-p(y)|=|\max\{p_1(x),p_2(x),\hdots,p_{\ell}(x)\}-\max\{p_1(y),p_2(y),\hdots, p_{\ell}(y)\}|\leq \max\{|p_1(x)-p_1(y)|,|p_2(x)-p_2(y)|,\hdots, |p_{\ell}(x)-p_{\ell}(y)|\}$ for each $x,y\in\overline{\Omega}$.}

			\item[$\bullet$] The weights $w_1,w_2,\hdots,w_{\ell}\in L^{\infty}(\Omega)$ have the property that there is some constant $\omega>0$ such that $\displaystyle\min_{j\in\overline{1,\ell}}\underset{x\in\Omega}{\operatorname{ess\ inf}}\ w_j(x)\geq\omega$.
			
		\end{enumerate}
		
		\noindent We have in this particular case that:
		
		\begin{enumerate}
			\item[$\bullet$] $\Psi:\overline{\Omega}\times (0,\infty)\to (0,\infty),\ \Psi(x,s)=w_1(x)s^{p_1(x)-2}+w_2(x)s^{p_2(x)-2}+\dots+w_{\ell}(x)s^{p_\ell(x)-2}$. So $\lim\limits_{s\to 0^+} \Psi(x,s)s=0$ for any $x\in\overline{\Omega}$.
			
			\item[$\bullet$] $	\mathbf{a}:\overline{\Omega}\times\mathbb{R}^N\to\mathbb{R}^N$,
				
			\begin{equation}
			\mathbf{a}(x,\xi)=\begin{cases} w_1(x)|\xi|^{p_1(x)-2}\xi+w_2(x)|\xi|^{p_2(x)-2}\xi+\hdots+w_\ell(x)|\xi|^{p_\ell(x)-2}\xi, & \xi\neq \mathbf{0}\\ \mathbf{0}, & \xi=\mathbf{0}\end{cases}.
			\end{equation}
			
			\item[$\bullet$] $\Phi:\overline{\Omega}\times\mathbb{R}\to [0,\infty),\ \Phi(x,s)=w_1(x)|s|^{p_1(x)-1}+w_2(x)|s|^{p_2(x)-1}+\dots+w_{\ell}(x)|s|^{p_{\ell}(x)-1}$ is strictly increasing on $[0,\infty)$ for any $x\in\Omega$.
			
			\item[$\bullet$] Hypothesis \textbf{(H6)} holds, Indeed, taking any $v\in L^{p(x)}(\Omega)$, for every $j\in\overline{1,\ell}$ we also have that $v\in L^{p_{j}(x)}(\Omega)$, since $p(x)\geq p_j(x)$ for all $x\in\overline{\Omega}$. So $|v|^{p_j(x)-1}\in L^{\frac{p_j(x)}{p_j(x)-1}}(\Omega)=L^{p'_j(x)}(\Omega)\subset L^{p'(x)}(\Omega)$,  because $p'(x)\leq p'_j(x)$ for all $x\in\overline{\Omega}$. Summing up we get that $\Phi\big (\cdot,v(\cdot)\big )\in L^{p'(x)}(\Omega)$.
			
			\item[$\bullet$] $A:\overline{\Omega}\times\mathbb{R}^N\to [0,\infty)$ and  $A(x,\xi)=w_1(x)\dfrac{|\xi|^{p_1(x)}}{p_1(x)}+w_2(x)\dfrac{|\xi|^{p_2(x)}}{p_2(x)}+\dots+w_{\ell}(x)\dfrac{|\xi|^{p_\ell(x)}}{p_\ell(x)}$.
			
			\item[$\bullet$] $\mathcal{A}:W^{1,p(x)}(\Omega)\to\mathbb{R},\ \mathcal{A}(u)=\displaystyle\int_{\Omega}w_1(x)\dfrac{|\nabla u|^{p_1(x)}}{p_1(x)}+w_2(x)\dfrac{|\nabla u|^{p_2(x)}}{p_2(x)}+\dots+w_{\ell}(x)\dfrac{|\nabla u|^{p_\ell(x)}}{p_\ell(x)}\ dx$

			\noindent Let's denote for each $j\in\overline{1,\ell}$: $\Omega_j=\{x\in\overline{\Omega}\ |\ p(x)=p_j(x)\}\subset\overline{\Omega}$. It is obvious that all of them are measurable and they cover $\overline{\Omega}$ meaning that: $\Omega_1\cup\Omega_2\cup\dots\cup\Omega_{\ell}=\overline{\Omega}$. To test hypothesis \textbf{(H7)} just note that for each $u\in W^{1,p(x)}(\Omega)$ we have:
			
			\begin{align*}
				\int_{\Omega}A(x,\nabla u)\ dx&=\int_{\Omega}w_1(x)\dfrac{|\nabla u|^{p_1(x)}}{p_1(x)}+w_2(x)\dfrac{|\nabla u|^{p_2(x)}}{p_2(x)}+\dots+w_{\ell}(x)\dfrac{|\nabla u|^{p_\ell(x)}}{p_\ell(x)}\ dx\\
				&\geq \omega\int_{\Omega}\dfrac{|\nabla u|^{p_1(x)}}{p_1(x)}+\dfrac{|\nabla u|^{p_2(x)}}{p_2(x)}+\dots+\dfrac{|\nabla u|^{p_\ell(x)}}{p_\ell(x)}\ dx\\
				&=\omega\left ( \int_{\Omega}\dfrac{|\nabla u|^{p_1(x)}}{p_1(x)}\ dx+\int_{\Omega}\dfrac{|\nabla u|^{p_2(x)}}{p_2(x)}\ dx+\dots+\int_{\Omega}\dfrac{|\nabla u|^{p_\ell(x)}}{p_\ell(x)}\ dx\right ) \\
				&\geq \omega\left (\int_{\Omega_1}\dfrac{|\nabla u|^{p_1(x)}}{p_1(x)}\ dx+\int_{\Omega_2}\dfrac{|\nabla u|^{p_2(x)}}{p_2(x)}\ dx+\dots+\int_{\Omega_{\ell}}\dfrac{|\nabla u|^{p_\ell(x)}}{p_\ell(x)}\ dx\right ) \\
				&=\omega\left (\int_{\Omega_1}\dfrac{|\nabla u|^{p(x)}}{p(x)}\ dx+\int_{\Omega_2}\dfrac{|\nabla u|^{p(x)}}{p(x)}\ dx+\dots+\int_{\Omega_{\ell}}\dfrac{|\nabla u|^{p(x)}}{p(x)}\ dx\right ) \\
				&\geq \omega\int_{\Omega_1\cup\Omega_2\dots\cup\Omega_{\ell}}\dfrac{|\nabla u|^{p(x)}}{p(x)}\ dx\\
				&\geq \omega\int_{\Omega}\dfrac{|\nabla u|^{p(x)}}{p(x)}\ dx\\
				&\geq \dfrac{\omega}{p^+}\int_{\Omega} |\nabla u|^{p(x)}\ dx.
			\end{align*}
		
		\item[$\bullet$] Finally, we will check that \textbf{(EH$_\Phi$)} holds for $\alpha\in (1,p_{-}
		)$, where $p_{-}:=\displaystyle\min_{x\in\overline{\Omega}}\min_{k\in\overline{1,\ell}} p_k(x)>1$. Indeed, for a.e. $x\in\Omega$, the function:
		
		\begin{align*}
			[0,\infty)\ni s\mapsto \dfrac{\Phi(x,s)}{s^{\alpha-1}}&=\dfrac{w_1(x)s^{p_1(x)-1}+w_2(x)s^{p_2(x)-1}+\hdots+w_\ell(x)s^{p_\ell(x)-1}}{s^{\alpha-1}}\\
			&=w_1(x)s^{p_1(x)-\alpha}+w_2(x)s^{p_2(x)-\alpha}+\hdots+w_\ell(x)s^{p_\ell(x)-\alpha},
		\end{align*} 
		
		\noindent is \textbf{strictly increasing} being a sum of strictly increasing functions. It is essential that $p_{k}(x)\geq p_{-}>\alpha$ for every $k\in\overline{1,\ell}$. 
		
	\end{enumerate}
		
		\bigskip
		
		\noindent Therefore applying Theorem \ref{uniq} \textbf{(3)} in this context we obtain the following corollary:
		
		\begin{corollary}\label{cormultiple}
			Problem \eqref{eqmphase2} has at most one strongly positive solution.
		\end{corollary}

		\section{Appendix}
		
		\begin{theorem}\label{athleb} Set $\rho_{p(x)}(u)=\displaystyle\int_{\Omega} |u(x)|^{p(x)}\ dx$ for $u\in L^{p(x)}(\Omega)$. Then:\footnote{See Theorem 1.3 from \cite{Fan2}.}
			
			\begin{enumerate}
				\item[\textnormal{\bf{(1)}}] $\Vert u\Vert_{L^{p(x)}(\Omega)}=a\ \Longleftrightarrow\ \rho_{p(x)}\left (\dfrac{u}{a} \right)=1$.
				\item[\textnormal{\bf{(2)}}] $\Vert u\Vert_{L^{p(x)}(\Omega)}=1\ \Longleftrightarrow\ \rho_{p(x)}\left (u \right)=1$.
				\item[\textnormal{\bf{(3)}}] $\Vert u\Vert_{L^{p(x)}(\Omega)}>1\ \Longleftrightarrow\ \rho_{p(x)}\left (u \right)>1$.
				\item[\textnormal{\bf{(4)}}] $\Vert u\Vert_{L^{p(x)}(\Omega)}<1\ \Longleftrightarrow\ \rho_{p(x)}\left (u \right)<1$.
				\item[\textnormal{\bf{(5)}}] $\Vert u\Vert_{L^{p(x)}(\Omega)}>1\ \Longrightarrow\ \Vert u\Vert_{L^{p(x)}(\Omega)}^{p_{-}}\leq \rho_{p(x)}(u)\leq \Vert u\Vert_{L^{p(x)}(\Omega)}^{p_{+}}$.
				\item[\textnormal{\bf{(6)}}] $\Vert u\Vert_{L^{p(x)}(\Omega)}<1\ \Longrightarrow\ \Vert u\Vert_{L^{p(x)}(\Omega)}^{p_{+}}\leq \rho_{p(x)}(u)\leq \Vert u\Vert_{L^{p(x)}(\Omega)}^{p_{-}}$.
				\item[\textnormal{\bf{(7)}}] For $(u_n)_{n\geq 1}\subset L^{p(x)}(\Omega)$ and $u\in L^{p(x)}(\Omega)$ we have: $\Vert u_n-u\Vert_{L^{p(x)}(\Omega)}\longrightarrow 0\ \Longleftrightarrow\ \rho_{p(x)}(u_n-u)\longrightarrow 0$.\footnote{The proof can be found in \cite[Proposition 2.56 and Corollary 2.58, page 44]{Cruz}.}
				\item[\textnormal{\bf{(8)}}] For $(u_n)_{n\geq 1}\subset L^{p(x)}(\Omega)$ we have: $\Vert u_n\Vert_{L^{p(x)}(\Omega)}\longrightarrow \infty\ \Longleftrightarrow\ \rho_{p(x)}(u_n)\longrightarrow \infty$.
				\item[\textnormal{\bf{(9)}}] For any $u\in L^{p(x)}(\Omega)$ we have that: $\lim\limits_{n\to\infty}\Vert u\chi_{D_n}\Vert_{L^{p(x)}(\Omega)}=0$, where $(D_n)_{n\geq 1}$ is a sequence of measurable sets included in $\Omega$ with the property that $\lim\limits_{n\to\infty} |D_n|=0$.\footnote{See Theorem 1.13 from \cite{Fan2}.}
				\item[\textnormal{\bf(10)}] Let $u\in L^{p(x)}(\Omega)$ and $v:\Omega\to\mathbb{R}$ a measurable function such that $|u|\geq |v|$ a.e. on $\Omega$. Then $v\in L^{p(x)}(\Omega)$ and $\Vert u\Vert_{L^{p(x)}(\Omega)}\geq \Vert v\Vert_{L^{p(x)}(\Omega)}$.
				\item[\textnormal{\bf(11)}] $\rho_{p(x)}(u+v)\leq 2^{p^+-1}\big (\rho_{p(x)}(u)+\rho_{p(x)}(v) \big )$.
				\item[\textnormal{\bf(12)}] If $\lambda\in [0,1]$ we have $\rho_{p(x)}(\lambda u)\leq \lambda \rho(u)$ for each $u\in L^{p(x)}(\Omega)$. 
				\item[\textnormal{\bf(13)}] If $\lambda\geq 1$ we have $\rho_{p(x)}(\lambda u)\geq \lambda \rho(u)$ for each $u\in L^{p(x)}(\Omega)$. 
				\item[\textnormal{\bf(14)}] If $u,v:\Omega\to\mathbb{R}$ are two measurable functions and $|u|\geq |v|$ a.e. on $\Omega$, then $\rho_{p(x)}(u)\geq \rho_{p(x)}(v)$ and the inequality is strict if $|u|\not\equiv|v|$.
				\item[\textnormal{\bf(15)}] For any $u\in L^{p(x)}(\Omega)$, the function $[1,\infty)\ni\lambda\longmapsto \rho_{p(x)}\left (\dfrac{u}{\lambda} \right )$ is continuous and decreasing. Furthermore: $\lim\limits_{\lambda\to\infty} \rho_{p(x)}\left (\dfrac{u}{\lambda}\right )=0$ for every $u\in L^{p(x)}(\Omega)$.\footnote{See \cite[Proposition 2.7, page 17]{Cruz}.}
				\item[\textnormal{\bf(16)}] $\rho_{p(x)}:L^{p(x)}(\Omega)\to\mathbb{R}$ is a convex and continuous functional.
				\item[\textnormal{\bf{(17)}}] If $\Omega\subseteq\Omega_1\cup\Omega_2$ then $\Vert u\Vert_{L^{p(x)}(\Omega)}\leq \Vert u\Vert_{L^{p(x)}(\Omega_1)}+\Vert u\Vert_{L^{p(x)}(\Omega_2)}$ for any $u\in L^{p(x)}(\Omega_1\cup\Omega_2)$.
				\item[\textnormal{\bf{(18)}}] If $D\subset\Omega$ are two measurable sets then $\Vert u\chi_{D}\Vert_{L^{p(x)}(\Omega)}=\Vert u\Vert_{L^{p(x)}(D)}$ for any $u\in L^{p(x)}(\Omega)$.
				\item[\textnormal{\bf{(19)}}] $\rho_{p(x)}:L^{p(x)}(\Omega)\to\mathbb{R}$ is a lower semicontinuous and a weakly lower semicontinuous functional.\footnote{See \cite[Theorem 2.1.17 and Remark 2.1.18]{Hasto}.}
				\item[\textnormal{\bf{(20)}}] If $\lambda\in [1,\infty)$ then $\lambda^{p^-}\rho_{p(x)}(u)\leq\rho_{p(x)}(\lambda u)\leq \lambda^{p^+}\rho_{p(x)}(u)$ for any measurable function $u:\Omega\to\mathbb{R}$.
				\item[\textnormal{\bf{(21)}}] If $\lambda\in (0,1)$ then $\lambda^{p^-}\rho_{p(x)}(u)\geq\rho_{p(x)}(\lambda u)\geq \lambda^{p^+}\rho_{p(x)}(u)$ for any measurable function $u:\Omega\to\mathbb{R}$.
				\item[\textnormal{\bf{(22)}}] If $\Vert u\Vert_{L^{p(x)}(\Omega)}>1$ then: $\rho_{p(x)}(u)^{\frac{1}{p^+}}\leq \Vert u\Vert_{L^{p(x)}(\Omega)}\leq\rho_{p(x)}(u)^{\frac{1}{p^-}}$.
				\item[\textnormal{\bf{(23)}}] If $0<\Vert u\Vert_{L^{p(x)}(\Omega)}<1$ then: $\rho_{p(x)}(u)^{\frac{1}{p^+}}\geq \Vert u\Vert_{L^{p(x)}(\Omega)}\geq\rho_{p(x)}(u)^{\frac{1}{p^-}}$.
				\item[\textnormal{\bf{(24)}}] If $p:\Omega\to [1,\infty)$ is a measurable bounded exponents and $q:\Omega\to (0,\infty)$ is also measurable such that $q(x)\leq p(x)$ a.e. on $\Omega$ then for any $u\in L^{p(x)}(\Omega)$ we have that $|u|^{q(x)}\in L^{\frac{p(x)}{q(x)}}(\Omega)$.
				\item[\textnormal{\bf{(25)}}] If $\Vert u\Vert_{L^{p(x)}(\Omega)}\leq 1$ then $\rho_{p(x)}(u)\leq \Vert u\Vert_{L^{p(x)}(\Omega)}$ and if $\Vert u\Vert_{L^{p(x)}(\Omega)}\geq 1$ then $\rho_{p(x)}(u)\geq \Vert u\Vert_{L^{p(x)}(\Omega)}$.\footnote{This is Corollary 2.22 from \cite[page 24]{Cruz}.}
				
			\end{enumerate}	
		\end{theorem}
		
		\begin{remark}\label{aremsob} For $1<p^{-}\leq p^{+}<\infty$ we have that $\big (W^{1,p(x)}(\Omega),\Vert \cdot\Vert_{W^{1,p(x)}(\Omega)}\big )$ and $W^{1,p(x)}_0(\Omega)$ are uniform convex\footnote{See \cite[Theorem 7]{Mendez}.} (hence reflexive) and separable\footnote{See \cite[Theorem 2.1]{Fan2} or \cite[Theorem 3.1]{kovacik}.} Banach spaces. Moreover, if $p$ is log-H\"{o}lder continuous and $\Omega$ is a Lipschitz bounded domain then $C^{\infty}(\Omega)$ is dense in $W^{1,p(x)}(\Omega)$ and $W^{1,p(x)}_0(\Omega)=W^{1,p(x)}(\Omega)\cap W^{1,1}_0(\Omega)$.\footnote{See \cite[Theorem 2.6]{Fan2}.}
		\end{remark}
		
		\begin{theorem}\label{athsob} Set $\varrho_{p(x)}(u)=\displaystyle\int_{\Omega} |u(x)|^{p(x)}+|\nabla u(x)|^{p(x)}\ dx$ for $u\in W^{1,p(x)}(\Omega)$. Then:
			
			\begin{enumerate}
				\item[\textnormal{\bf{(1)}}] $\Vert u\Vert_{W^{1,p(x)}(\Omega)}>1\ \Longrightarrow\ \Vert u\Vert_{W^{1,p(x)}(\Omega)}^{p_{-}}\leq \varrho_{p(x)}(u)\leq \Vert u\Vert_{W^{1,p(x)}(\Omega)}^{p_{+}}$.
				\item[\textnormal{\bf{(2)}}] $\Vert u\Vert_{W^{1,p(x)}(\Omega)}<1\ \Longrightarrow\ \Vert u\Vert_{W^{1,p(x)}(\Omega)}^{p_{+}}\leq \varrho_{p(x)}(u)\leq \Vert u\Vert_{W^{1,p(x)}(\Omega)}^{p_{-}}$.
				\item[\textnormal{\bf{(3)}}] $\Vert u_n-u\Vert_{W^{1,p(x)}(\Omega)}\longrightarrow 0\ \Longleftrightarrow\ \varrho_{p(x)}(u_n-u)\longrightarrow 0$.
				\item[\textnormal{\bf{(4)}}] $\Vert u_n\Vert_{W^{1,p(x)}(\Omega)}\longrightarrow \infty\ \Longleftrightarrow\ \varrho_{p(x)}(u_n)\longrightarrow \infty$.
			\end{enumerate}
			
		\end{theorem}

		\begin{theorem}[\textbf{Nemytskii Operators}]\label{athnem} Let $f:\Omega\times\mathbb{R}\to\mathbb{R}$ be a Carath\'{e}odory function and $p:\Omega\to [1,\infty)$ be a measurable and bounded exponent. For each function $u\in L^{p(x)}(\Omega)$ we consider $\mathcal{N}_f(u):\Omega\to\mathbb{R},\ \mathcal{N}_f(u)(x)=f(x,u(x)),\ x\in\Omega$. 
			
			\noindent $\bullet$ If there is a non-negative function $g\in L^{q(x)}(\Omega)$, where $q:\Omega\to [1,\infty)$ is a measurable bounded exponent, and a constant $c\geq 0$ such that:
			\begin{equation}\label{aeqnem}
				|f(x,s)|\leq g(x)+c|s|^{\frac{p(x)}{q(x)}},\ \forall\ s\in\mathbb{R}\ \text{and}\ \text{for a.e.}\ x\in\Omega
			\end{equation}
			\noindent then $\mathcal{N}_f:L^{p(x)}(\Omega)\to L^{q(x)}(\Omega)$ is a continuous and bounded operator.

			\noindent $\bullet$ Conversely if it happens that $\mathcal{N}_f:L^{p(x)}(\Omega)\to L^{q(x)}(\Omega)$, then $\mathcal{N}_f$ is continuous and bounded and moreover there is a non-negative function $g\in L^{q(x)}(\Omega)$ and some constant $c\geq 0$ such that \eqref{aeqnem} holds.\footnote{For a complete proof see \cite[Theorem 1.16]{Fan2} or \cite[Theorem 4.1 and Theorem 4.2]{kovacik}.}
		\end{theorem}
		
		\begin{remark}\label{arenem} Theorem \ref{athnem} generalizez to Carath\'{e}odory functions of the type $f:\Omega\times\mathbb{R}^N\to\mathbb{R}$. For each $u=(u_1,u_2,\dots,u_N)\in L^{p_1(x)}(\Omega)\times L^{p_2(x)}(\Omega)\times\dots\times L^{p_N(x)}(\Omega)$, where $p_1,p_2,\dots,p_N:\Omega\to [1,\infty)$ are measurable and bounded exponents, we take the \textit{Nemytskii operator} $\mathcal{N}_f(u):\Omega\to\mathbb{R},\ \mathcal{N}_f(u)(x)=f\big (x,u(x)\big )=f\big (x,u_1(x),u_2(x),\dots,u_N(x)\big )$.
			
			\noindent $\bullet$ If there is a non-negative function $g\in L^{q(x)}(\Omega)$, where $q:\Omega\to [1,\infty)$ is a measurable bounded exponent, and some constants $c_1,c_2,\dots,c_N\geq 0$ such that:
			\begin{equation}\label{aeqnem2}
				|f(x,s_1,s_2,\dots,s_N)|\leq g(x)+c_1|s_1|^{\frac{p_1(x)}{q(x)}}+\dots+c_N|s_N|^{\frac{p_N(x)}{q(x)}},\ \forall\ s_1,s_2,\dots,s_N\in\mathbb{R}\ \text{and}\ \text{for a.e.}\ x\in\Omega
			\end{equation}
			\noindent then $\mathcal{N}_f:L^{p_1(x)}(\Omega)\times L^{p_2(x)}(\Omega)\times\dots\times L^{p_N(x)}(\Omega)\to L^{q(x)}(\Omega)$ is a continuous and bounded operator.
			
			\noindent $\bullet$ Conversely if it happens that $\mathcal{N}_f:L^{p_1(x)}(\Omega)\times L^{p_2(x)}(\Omega)\times\dots\times L^{p_N(x)}(\Omega)\to L^{q(x)}(\Omega)$, then $\mathcal{N}_f$ is continuous and bounded and moreover there is a non-negative function $g\in L^{q(x)}(\Omega)$ and some constants $c_1,c_2,\dots,c_N\geq 0$ such that \eqref{aeqnem2} holds.
		\end{remark}
			
		\begin{theorem}\label{athmgatfre} Let $\big (X,\Vert\cdot\Vert_X\big )$ and $\big (Y,\Vert\cdot\Vert_Y\big )$ be two real normed linear spaces. If the function $f:X\to Y$ is \textbf{G\^ateaux differentiable} in an open neighborhood $U\subset X$ of a point $x_0\in X$ and the mapping $\partial f: U\to\mathcal{L}(X,Y)$ is continuous then $f$ is \textbf{Fr\'{e}chet differentiable} at $x_0$ and $f'(x_0)=\partial f(x_0)$.\footnote{This result is taken from \cite[Proposition 3.2.15, page 133]{drabek}. It can also be found as \cite[Theorem 1.1.3, page 3]{chang}.}
		\end{theorem}
		
		\begin{theorem}\label{athmstru} Let $\big (X,\Vert\cdot\Vert_X\big )$ be a real reflexive Banach space and $M\subset X$ be a weakly closed subset of $X$. Suppose that the functional $\mathcal{J}:M\to\mathbb{R}\cup\{+\infty\}$ is:
			
			\begin{enumerate}
				\item[$\bullet$] \textbf{coercive on $M$}, i.e.: $\forall\ r>0,\ \exists\ R>0$ such that $\forall\ u\in M$ with $\Vert u\Vert_X\geq R$ it follows that $\mathcal{J}(u)\geq r$;
				
				\item[$\bullet$] \textbf{sequentially weakly lower semicontinuous on $M$}, i.e.: for any $u\in M$ and any sequence $(u_n)_{n\geq 1}\subset M$ such that $u_n\rightharpoonup u$ in $X$ there holds:
				
				\begin{equation*}
					\liminf\limits_{n\to\infty}\mathcal{J}(u_n)\geq \mathcal{J}(u).
				\end{equation*}
			\end{enumerate}
			
			\noindent Then $\mathcal{J}$ is bounded from below on $M$, meaning that $\displaystyle\inf_{u\in M} \mathcal{J}(u)>-\infty$, and $\mathcal{J}$ attains its infimum on $M$, i.e. $\exists\ u\in M$ such that $\mathcal{J}(u)=\displaystyle\inf_{v\in M} \mathcal{J}(v)$.\footnote{This is Theorem 1.2 that can be found at page 4 in \cite{struwe}.}
		\end{theorem}
		
		\begin{lemma}\label{alemhos} Let $a,b\in\mathbb{R}$ and $p>1$. Then $\lim\limits_{\theta\to 0}\dfrac{|a+b\theta|^p-|a|^p}{p\theta}=|a|^{p-2}ab$.
		\end{lemma}
		
		\begin{proof} If $a>0$ then for $|\theta|$ sufficiently small we also have that $a+b\theta>0$. Therefore, using \textit{L'H\^{o}spital Rule}, we get that:
			
			\begin{equation}
				\lim\limits_{\theta\to 0} \dfrac{|a+b\theta|^p-|a|^p}{p\theta}=\lim\limits_{\theta\to 0}\dfrac{(a+b\theta)^p-a^p}{p\theta}=\lim\limits_{\theta\to 0} (a+b\theta)^{p-1}b\stackrel{p>1}{=}a^{p-1}b=|a|^{p-2}ab.
			\end{equation}
			
			\noindent If $a<0$ then for $|\theta|$ sufficiently small we also have that $a+b\theta<0$. So:
			
			\begin{equation}
				\lim\limits_{\theta\to 0} \dfrac{|a+b\theta|^p-|a|^p}{p\theta}=\lim\limits_{\theta\to 0}\dfrac{(-a-b\theta)^p-(-a)^p}{p\theta}=\lim\limits_{\theta\to 0} (-a-b\theta)^{p-1}(-b)\stackrel{p>1}{=}(-a)^{p-1}(-b)=|a|^{p-2}ab.
			\end{equation}
			
			\noindent For $a=0$:
			
			\begin{equation}
				\lim\limits_{\theta\to 0} \dfrac{|a+b\theta|^p-|a|^p}{p\theta}=\lim\limits_{\theta\to 0}\dfrac{|b|^p\cdot |\theta|^p}{p\theta}=\dfrac{|b|^p}{p}\lim\limits_{\theta\to 0}|\theta|^{p-1}\operatorname{sgn}(\theta)\stackrel{p>1}{=}0=|a|^{p-2}ab.
			\end{equation}
		\end{proof}
		
		\begin{lemma}\label{alplap} For every $\mathbf{u},\mathbf{v}\in\mathbb{R}^N$ the following inequalities hold for some constants $c_1,c_2>0$:\footnote{See Lemma 3.3 and Lemma 3.4 from \cite{bystrom} where you can find the best constants $c_1,c_2$ and the cases when equality holds.}
			
			\begin{equation} \big | |\mathbf{u}|^{p-2}\mathbf{u}-|\mathbf{v}|^{p-2}\mathbf{v}\big |\leq c_1 |\mathbf{u}-\mathbf{v}|^{p-1},\ \text{if}\ p\in (1,2),
			\end{equation}
			
			\noindent and:
			
			\begin{equation}
				\big | |\mathbf{u}|^{p-2}\mathbf{u}-|\mathbf{v}|^{p-2}\mathbf{v}\big |\leq c_2\big (|\mathbf{u}|+|\mathbf{v}| \big )^{p-2}\cdot |\mathbf{u}-\mathbf{v}|,\ \text{if}\ p\in (2,\infty).
			\end{equation}
			
		\end{lemma}

	\bibliographystyle{apalike}
	\bibliography{doubly}
\end{document}